\documentclass[12pt,reqno]{amsart}
\usepackage{amsmath,amsthm,amscd,amsfonts,eepic}

\swapnumbers
\newtheorem{thm}[subsection]{Theorem}
\newtheorem{cor}[subsection]{Corollary}
\newtheorem{lem}[subsection]{Lemma}
\newtheorem{prop}[subsection]{Proposition}

\theoremstyle{definition}
\newtheorem{defn}[subsection]{Definition}
\newcommand{\quash}[1]{}
\newcommand{\ts}[1]{{\scriptstyle{#1}}}
\numberwithin{equation}{subsection}

\begin{document}

\centerline{\bf Chern classes of automorphic vectorbundles}
\medskip
\bigskip

\centerline{Mark Goresky\footnote{School of Mathematics, Institute for Advanced Study,
Princeton N.J.  Research partially supported by NSF grants \# DMS 9626616 and DMS
9900324.}
and William Pardon \footnote{Dept. of Mathematics, Duke University, Durham N.C.  Research
partially supported by NSF grant \# DMS 9504900.} }

\setcounter{section}{0}
\section{Introduction}
\subsection{}  Suppose $X$ is a compact $n$-dimensional complex manifold.  Each partition
$I=$ $\{i_1,i_2,$ $\ldots,i_r\}$ of $n$ corresponds to a Chern number $c^I(X) = \epsilon
(c^{i_1}(X)\cup c^{i_2}(X) \cup \ldots \cup c^{i_r}(X) \cap [X]) \in \mathbb Z$ where
$c^k(X)\in H^{2k}(X;\mathbb Z)$ are the Chern classes of the tangent bundle, $[X] \in
H_{2n}(X;\mathbb Z)$ is the fundamental class, and $\epsilon:H_0(X;\mathbb Z) \to \mathbb
Z$ is the augmentation.  Many invariants of $X$ (such as its complex cobordism class) may
be expressed in terms of its Chern numbers (\cite{Milnor2}, \cite{Stong}).  During the
last 25 years, characteristic classes of singular spaces have been
defined in a variety of contexts:  Whitney classes of Euler spaces \cite{Sullivan},
\cite{Halperin}, \cite{Akin}, Todd classes of singular varieties \cite{BFM}, Chern
classes of singular algebraic varieties \cite{Mac}, L-classes of stratified spaces with
even codimension strata \cite{GM1}, Wu classes of singular spaces \cite{Go1}, \cite{GP}
(to name a few).  However, these characteristic classes are invariably {\it homology}
classes and as such, they cannot be multiplied with each other.  In some cases it has
been found possible to ``lift'' these classes from homology to intersection homology,
where (some) characteristic numbers may be formed (\cite {BBF}, \cite{BW},
\cite{Go1},\cite{GP}, \cite{Totaro}). 

The case of locally symmetric spaces is particularly interesting.
Suppose $\Gamma$ is a torsion-free arithmetic group acting on a complex $n$-dimensional
Hermitian symmetric domain $D=G/K$, where $G$ is the group of real points of a semisimple
algebraic group $\mathbf G$ defined over $\mathbb Q$ with $\Gamma \subset \mathbf
G(\mathbb Q)$, and where $K\subset G$ is a maximal compact subgroup.  Then $X = \Gamma
\backslash D$ is a Hermitian locally symmetric space.  When $X$ is compact,  Hirzebruch's
proportionality theorem \cite{Hirzebruch} says that there is a number
$v(\Gamma)\in\mathbb Q$ so that for every partition $I = \{ i_1, i_2,\ldots , i_r \}$ of
$n$,  the Chern number satisfies $c^I(X) = v(\Gamma)c^I(\check{D}),$ where
$\check{D}=G_u/K$ is the compact dual symmetric space (and $G_u$ is a compact real form
of $\mathbf G(\mathbb C)$ containing $K$).
 
If $X = \Gamma \backslash D$ is noncompact,  it has a canonical
{\it Baily-Borel (Satake) compactification}, $\overline{X}$.  This is a (usually highly
singular) complex projective algebraic variety.  To formulate a proportionality theorem
in the noncompact case, one might hope that the tangent bundle $TX$ extends as a complex
vectorbundle over $\overline{X}$, but this is false.  In \cite{Mumford}, D. Mumford
showed that $TX$ has a particular extension $\overline{E}_{\Sigma} \to
\overline{X}_{\Sigma}$ over any toroidal resolution $\tau:\overline{X}_{\Sigma} \to
\overline{X}$ of the Baily-Borel compactification and that for any partition $I$ of $n$,
the resulting Chern numbers
\begin{equation}\label{eqn-Hirzebruch}
c^I(\overline{E}_{\Sigma}) = \epsilon(c^{i_1}(\overline{E}_{\Sigma}) \cup
c^{i_2}(\overline{E}_{\Sigma}) \cup \ldots \cup
c^{i_r}(\overline{E}_{\Sigma}) \cap [\overline{X}_{\Sigma}]) \end{equation}
satisfy the same equation, $c^I(\overline{E}_{\Sigma}) = v(\Gamma)c^I(\check{D}).$ (The
toroidal resolution $\overline{X}_{\Sigma}$ is constructed in \cite{AMRT}; it depends on
a choice $\Sigma$ of polyhedral cone decompositions of certain self-adjoint homogeneous
cones.) Mumford also showed that if $\Sigma'$ is a refinement of $\Sigma$ then there is a
natural morphism $f:\overline{X}_{\Sigma'} \to \overline{X}_{\Sigma}$ and that
$f^*(\overline{E}_{\Sigma}) \cong \overline{E}_{\Sigma'}$ (hence
$f^*c^i(\overline{E}_{\Sigma}) = c^i(\overline{E}_{\Sigma'})$).
Moreover, it is proven in \cite{Harris3} that the coherent sheaf
$\tau_*\overline{E}_{\Sigma}$ is independent of the choice of $\Sigma.$  One is therefore
led to suspect the existence of a closer relationship between the characteristic classes
of the vectorbundles $\overline{E}_{\Sigma}$ and the topology of the Baily-Borel
compactification $\overline{X}.$  In Theorem \ref{thm-A} and Theorem \ref{thm-toroidal}
we show that, at least for the variety $\overline{X}$, the original goal of constructing
Chern numbers can be completely realized:

\subsection*{Theorem} {\it  Every Chern class $c^i(X)$ has a canonical lift $\overline
c^i \in H^{2i}(\overline{X};\mathbb C)$ to the cohomology of the Baily-Borel
compactification.  Moreover, if $\tau:\overline{X}_{\Sigma} \to \overline{X}$ is any
toroidal resolution of singularities then
\begin{equation*} \tau^*(\overline{c}^i) = c^i(\overline{E}_{\Sigma})
\in H^{2i}(\overline{X}_{\Sigma};\mathbb C).  \end{equation*} }
\noindent
It follows (\S \ref{subsec-prop}) that the lifts $\overline c^i$ satisfy
(\ref{eqn-Hirzebruch}).  In \S \ref{sec-constructible} we show that the
homology image $\bar c^*\cap [\bar X] \in H_*(\bar X)$ lies in integral homology and
coincides with (MacPherson's) Chern class \cite{Mac} of the constructible function which
is $1$ on $X$ and is $0$ on $\overline{X} -X.$
\subsection{}
Moreover, a similar result holds for any automorphic vectorbundle.  Let $\lambda:K \to
GL(V)$ be a representation of $K$ on some finite dimensional complex vectorspace $V$. By
\cite{Mumford}, the automorphic vectorbundle $E' = (\Gamma \backslash G)\times_KV $ on
$X$ has a particular extension $\overline{E}'_{\Sigma}$ over any toroidal resolution
$\overline{X}_{\Sigma}$.  We show that each Chern class $c^i(E_{\Gamma})$ has a canonical
lift $\bar{c}^i(E') \in H^{2i}(\overline{X};\mathbb C)$ 
and that these lifts also satisfy the proportionality formula. Moreover,
$\tau^*(\bar c^i(E')) = c^i(\overline{E}'_{\Sigma})$  and the image of $\bar
c^i(E')$ in $\text{Gr}^W_{2i}(\overline{X};\mathbb C)$ (the top graded piece of
the weight filtration) is uniquely determined by this formula. 

\subsection{}
In \S \ref{sec-cohomology} we consider the subalgebra
$H^*_{\text{Chern}}(\overline{X};\mathbb C)$ of the cohomology of the Baily-Borel
compactification that is generated by the (above defined lifts of) Chern classes of
certain ``universal'' automorphic vectorbundles, and show that 
\subsection*{Theorem} {\it Suppose the Hermitian symmetric domain $D$ is a product of
irreducible factors $G_i/K_i$ (where $K_i$ is a maximal compact subgroup of $G_i$), and
that each $G_i$ is one of the following:  $Sp_n(\mathbb R)$, $U(p,q)$,
$SO(2n)$, or $SO(2,p)$ with $p$ odd or $p=2$.  Then there is a (naturally defined)
surjection $h:H^*_{\text{Chern}}(\overline{X};\mathbb
C) \to H^*(\check{D};\mathbb C)$ from this subalgebra to the cohomology of the compact
dual symmetric space.}  \medskip

This result is compatible with the few known facts about the cohomology of the
Baily-Borel compactification.  Charney and Lee \cite{CharneyLee} have shown, when
$D=Sp_{2n}(\mathbb R)/U(n)$ is the Siegel upper half space, and when $\Gamma =
Sp_{2n}(\mathbb Z)$ that the ``stable'' cohomology of $\overline{X}$ contains a
polynomial algebra which  coincides with the ``stable'' cohomology of the compact dual
symmetric space $\check{D}$ (which is the complex Lagrangian Grassmannian).  It is a
general fact (cf. \S \ref{sec-cohomology}) that the intersection cohomology
$IH^*(\overline{X};\mathbb C)$ contains a copy of $H^*(\check{D};\mathbb C).$
\subsection{}
Here are the main ideas behind the proof of Theorem \ref{thm-A}.  In \cite{Hi2},
Hirzebruch shows that the Chern classes of Hilbert modular varieties have  lifts to
the cohomology of the Baily-Borel compactification because the tangent bundle has a
trivialization in a neighborhood of each of the finitely many cusp points, cf.
\cite{Donnelly}.  If $X = \Gamma\backslash D$ is a $\mathbb Q$-rank 1 locally symmetric
space such that $\overline{X}$ is obtained from $X$ by adding finitely many cusps, then
the tangent bundle is not necessarily trivial near each cusp, but it admits a
connection which is flat near each cusp, and so the Chern forms vanish near each cusp,
hence the Chern classes lift to the cohomology of the Baily-Borel compactification. 
 A similar argument applies to arbitrary automorphic vectorbundles (cf. \cite{HZ1} \S
3.3.9). 

In the general $\mathbb Q$-rank 1 case, the singular set of the Baily-Borel
compactification $\overline{X}$ consists of finitely many disjoint smooth compact
manifolds (rather than finitely many cusp points).  If $Y$ denotes such a singular
stratum, then it admits a  neighborhood $\pi_Y:N_Y\to Y$ such that every slice
$\pi_Y^{-1}(y)\cap X$ is diffeomorphic to a neighborhood of a cusp similar to the kind
described above.  It is then possible to construct a connection $\nabla$ (on the tangent
bundle) which is ``flat along each fiber $\pi_Y^{-1}(y)$''.  (We call this the
``parabolic connection''; it is constructed in Section \ref{sec-induced}.)  Moreover,
within the neighborhood $N_Y$, each Chern form $\sigma^i(\nabla)$ is the pullback
$\pi_Y^*(\sigma_Y^i)$ of a certain differential form $\sigma_Y^i$ on $Y$.  Differential
forms with this ``$\pi$-fiber property'' form a complex whose cohomology is the
cohomology of $\overline{X}$, as discussed in Section \ref{sec-forms}.  So the Chern form
$\sigma^i(\nabla)$ determines a class $\overline c^i(\nabla) \in
H^{2i}(\overline{X};\mathbb C).$  (In fact, even the curvature form satisfies the
$\pi$-fiber condition.)  
\subsection{}
In higher rank cases, there are more problems.  If $Y_1 \subset \overline{Y}_2
\subset \overline{X}$ are singular strata of the Baily-Borel compactification, then it is
possible to define a ``parabolic'' connection in a neighborhood $N(Y_1)$ of $Y_1$ whose
curvature form has the $\pi$-fiber property relative to the tubular projection
$\pi_1:N(Y_1)\to Y_1.$ It is also possible to construct a ``parabolic'' connection in a
neighborhood $N(Y_2)$ of $Y_2$ whose curvature form has the $\pi$-fiber property relative
to the tubular projection $\pi_2:N(Y_2)\to Y_2.$  However these two connections do not
necessarily agree on the intersection $N(Y_1)\cap N(Y_2)\cap X$, nor do their curvature
forms.  When we patch these two connections together using a partition of unity, the
curvature form of the resulting connection fails to have the $\pi$-fiber property.
Nevertheless it is possible (as explained in Remark \ref{subsec-explanation}) to patch
together connections of this type so as to obtain a connection whose curvature form
$\Omega \in \text{End}(V)$ differs from a $\pi$-fiber differential form by a nilpotent
element $n \in \text{End}(V)$ which commutes with $\Omega$ (cf. \S
\ref{subsec-completion}).  (Here, $V$ is the representation of $K$ that gives rise to
the automorphic vectorbundle $E_{\Gamma} = (\Gamma \backslash G)\times_KV$ on $X = \Gamma
\backslash G /K.$)  This is enough to guarantee that the Chern forms of this ``patched''
connection are $\pi$-fiber differential forms (cf. Lemma \ref{lem-Springer}).  A standard
argument shows that the resulting cohomology class is independent of the choices that
were involved in the construction.

\subsection{}  A number of interesting questions remain.  We do not know whether the
results on Chern classes which are described in this paper for Hermitian symmetric spaces
may be extended to the ``equal rank'' case (when the real rank of $G$ and of $K$
coincide).  We do not know if the lifts $\overline c^i(E') \in
H^{2i}(\overline{X};\mathbb C)$ are integer or even rational cohomology
classes.   We do not know to what extent these lifts are uniquely determined by the
properties (\ref{cor-lift}), (\ref{eqn-toroidal}), (\ref{thm-constructible}).   We do not
know whether similar techniques can be applied to the Euler class of automorphic
vectorbundles (when such a class exists:  see \S \ref{sec-cohomology}).  We do not know
whether the surjection $h$ of Theorem \ref{thm-cohomology} admits a natural splitting.
We expect that $\bar c^*(E')=0$ whenever the automorphic vectorbundle
$E'$ arises from a representation $\lambda:K \to GL(V)$ which extends to a
representation of $G$.  If $\overline{E}^{\text{\tiny RBS}}_{\Gamma}$ denotes the
canonical extension (\cite{GT} \S9) of the automorphic vectorbundle $E$ over the
reductive Borel-Serre (\cite{Zucker} \S 4.2 p. 190, \cite{GHM} \S 8) compactification
$\nu:\overline{X}^{\text{\tiny RBS}} \to \overline{X}$ then it is likely that $\nu^*(\bar
c^*(E')) = c^*(\overline{E}^{\text{\tiny RBS}}).$  We expect these results to
have interesting applications to the study of the signature defect (\cite{Hi2} \S 3,
\cite{Donnelly}, \cite {Muller}, \cite{Stern1}) and to variations of weight 1 (and some
weight 2) Hodge structures (\cite{Gr1}, \cite{Gr2}).
\subsection{}
We would like to thank A. Borel, R. Bryant, D. Freed, R. Hain, M. Harris, R. MacPherson,
and L. Saper for valuable conversations.  We are grateful to S. Zucker for many
suggestions and comments on an earlier draft of this paper.  We are especially grateful
to an anonymous referee for pointing out a mistake in an earlier version of this paper,
and for his many helpful comments, suggestions, and corrections.  Both authors would like
to thank the Institute for Advanced Study in Princeton N.J. for its support and
hospitality while this paper was prepared.


\section{   Control Data}\label{sec-controldata}
\subsection{ }\label{subsec-controldata}
  A {\it weakly stratified space} $W$ is a compact Hausdorff space with a decomposition
into finitely many smooth manifolds $W=Y_1 \cup Y_2 \cup \ldots \cup Y_r$ (called the
strata of $W$) which satisfy the axiom of the frontier:	If $Y$ and $Z$ are strata and
if $Z\cap \bar Y \neq \phi$ then $Z \subset \bar Y$; we write $Z<Y$ and say that $Z$ is
{\it incident} to $Y$.  The {\it boundary} $\partial \bar Y = \bar Y - Y = \cup_{Z<Y}Z$
of the stratum $Y$ is the union of all strata incident to $Y$. If $W = \bar X$ is the
closure of a single stratum $X$ then we say that $X$ is the {\it nonsingular} part of $W$
and the other strata $Y<X$ are {\it boundary} or {\it singular} strata of $W$.  Fix a
positive real number $\epsilon > 0$.
\begin{defn}\label{def-control}
An {\it $\epsilon$-system of control data} on a weakly stratified space $W$ is a
collection $\{T_Y(\epsilon),\pi_Y,\rho_Y\}$ indexed by the boundary strata $Y\subset W$,
where
\begin{enumerate}
\item $T_Y(\epsilon)\subset W$ is an open subset of $W$ containing $Y$,
\item\label{item-frontier}
       $T_Y(\epsilon) \cap T_Z(\epsilon) = \phi$ unless $Y<Z$ or $Z<Y$.
\item\label{item-retraction}
      The {\it tubular projection} $\pi_Y:T_Y(\epsilon)\to Y$ is a retraction of
$T_Y(\epsilon)$ to $Y$ which is smooth on each stratum,
\item\label{item-comm1}
      $\pi_Z\pi_Y=\pi_Z$ whenever $Z<Y$ and both sides of the equation are
defined,
\item\label{item-rho} $\rho_Y:T_Y(\epsilon)\to [0,\epsilon )$ is a continuous ``distance
function'',  with $\rho_Y^{-1}(0)=Y,$ such that the mapping
$(\rho_Y,\pi_Y):T_Y(\epsilon) \to [0,\epsilon)\times Y$ is proper and its
restriction to each stratum is a submersion,
\item\label{item-comm2}
 $\rho_Z\pi_Y=\rho_Z$ whenever $Z<Y$ and both sides of the equation are defined,
\end{enumerate}\end{defn}
For $\tau \le \epsilon$, write $T_Y(\tau) = \rho_Y^{-1}([0,\tau)).$  By shrinking the
neighborhood $T_Y(\epsilon)$ and scaling $\rho_Y$ if necessary, we may assume that each
$\rho_Y$ is defined on a slightly larger neighborhood $T_Y(\epsilon')$ (where $\epsilon'
> \epsilon$) and that the ``boundary'' of $T_Y(\epsilon)$ is 
\[ \partial T_Y(\epsilon) = \overline{T_Y(\epsilon)}-T_Y(\epsilon) =
\overline{\rho_Y^{-1}(\epsilon)} = \rho_Y^{-1}(\epsilon)\cup \partial \overline{Y}.\]
Such neighborhoods are illustrated in the following diagram.

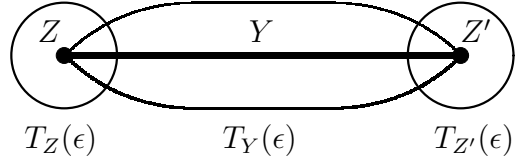
\begin{figure}[h!]\begin{center}
\centerline{\begin{picture}(200,80)
\put(40,40){\circle*{6}}
\linethickness{2pt}
\qbezier(40,40)(110,40)(190,40)
\linethickness{.5pt}
\qbezier(40,40)(60,60)(90,60)
\qbezier(90,60)(120,60)(140,60)
\qbezier(140,60)(170,60)(190,40)
\put(190,40){\circle*{6}}
\thicklines
\put(40,40){\circle{40}}
\put(190,40){\circle{40}}
\qbezier(40,40)(60,20)(90,20)
\qbezier(90,20)(120,20)(140,20)
\qbezier(140,20)(170,20)(190,40)
\put(30,45){$Z$}
\put(190,45){$Z'$}
\put(110,45){$Y$}
\put(25,5){$T_Z(\epsilon)$}
\put(100,5){$T_{Y}(\epsilon)$}
\put(180,5){$T_{Z'}(\epsilon)$}
\end{picture}}\caption{Tubular neighborhoods}\end{center}
\end{figure}

\noindent
If $W$ is the closure of a single stratum $X$ we extend this notation by
setting $T_X(\epsilon) =X$, $\pi_X(x) = x$ and $\rho_X(x)=0$ for all $x\in X.$

\subsection{}\label{subsec-existence}
Any compact real or complex algebraic or analytic variety admits a Whitney
stratification (\cite{Hardt1}, \cite{Hardt2}).  Any compact Whitney stratified subset of
a smooth manifold admits a system of control data (see \cite{Mather} or \cite{Gibson}
Thm. 2.6).   If $W$ is a compact Whitney stratified
set and if the mappings $\{\pi_Y\}$ are preassigned so as to satisfy Conditions
(\ref{item-retraction}) and (\ref{item-comm1}) above, then distance functions $\rho_Y$
may be found which are compatible with the mappings $\pi_Y.$

\section{   Partition of Unity  }\label{sec-partition}
\subsection{  }\label{subsec-bump}  
Throughout this paper we will fix a choice of a smooth nondecreasing function $s:\mathbb
R \to
[0,1]$ so that $s(x) = 0$ for all $x \le 1/2$ and $s(x) =1$ for all $x \ge 3/4.$  For any
$\epsilon > 0$ define $s_{\epsilon}(\rho) = s(\rho/\epsilon).$
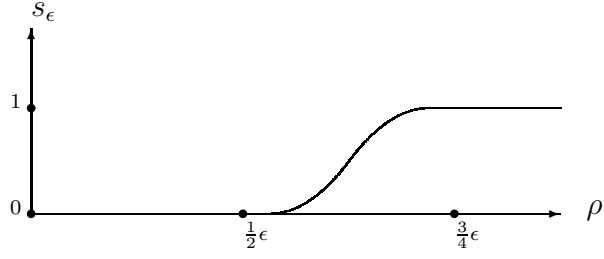
\begin{figure}[h!]\begin{center}
\centerline{\begin{picture}(200,80)
\put(0,0){\vector(1,0){200}}  \put(210,0){$\rho$}
\put(0,0){\vector(0,1){70}}   \put(0,75){$s_{\epsilon}$}
\put(0,0){\circle*{3}} \put(0,40){\circle*{3}}
\put(-8,0){$\scriptstyle{0}$}  
\put(-8,40){$\scriptstyle{1}$}
\qbezier(90,0)(105,0)(120,20)
\qbezier(120,20)(135,40)(150,40)
\qbezier(150,40)(175,40)(200,40)
\put(80,0){\circle*{3}} \put(160,0){\circle*{3}}
\put(80,-10){$\scriptstyle{\frac{1}{2}\epsilon}$}  
\put(160,-10){$\scriptstyle{\frac{3}{4}\epsilon}$}
\end{picture}}\end{center}\caption{The function $s_{\epsilon}(\rho).$}
\end{figure}

Fix $0 < \epsilon \le \epsilon_0.$  Let $W$ be a weakly stratified
space with an $\epsilon_0$ system of control data $\{T_Y(\epsilon_0), \pi_Y,
\rho_Y \}$.  For each stratum $Z\subset W$ define the modified distance function
$s_Z^{\epsilon}:T_Z (\epsilon)\to [0,1]$ by $s_Z^{\epsilon}(x)=s_{\epsilon}(\rho_Z(x))$.
Then $s_Z^{\epsilon}=0$ on $T_Z(\frac{\epsilon}2)$ and $s_Z^{\epsilon}=1$ near the edge
$\partial T_Z(\epsilon)$ of the tubular neighborhood $T_Z(\epsilon)$.

For each stratum $Y \subset W$ define a smooth function $t^Y_{\epsilon}:Y\to \mathbb R$
as follows:  If $y\in Y$ is not contained in the tubular neighborhood $T_Z(\epsilon )$ of
{\it any} stratum $Z<Y$ then set  $t^Y_{\epsilon}(y)=1.$ Otherwise, there is a unique
maximal collection of boundary strata $Z_1,Z_2,\ldots ,Z_r$ such that $y\in
T_{Z_1}(\epsilon )\cap T_{Z_2}(\epsilon )\cap\ldots\cap T_{Z_r}(\epsilon )$ and in this
case, by (\ref{def-control}) (Condition \ref{item-frontier}), these boundary strata form
a flag $Z_1<Z_2<\ldots <Z_r$ (after possibly relabeling the indices).  Define
\begin{equation}
t^Y_{\epsilon}(y)=s_{Z_1}^{\epsilon}(y).s_{Z_2}^{\epsilon}(y)\ldots
s_{Z_r}^{\epsilon}(y).
\end{equation}  
Then the function $t^Y_{\epsilon}:Y\to [0,1]$ is smooth and vanishes on 
\[\bigcup_{Z<Y}T_Z(\textstyle{\frac{\epsilon}{2}})\cap Y.\]  
Pull this up to a function $\pi_Y^{*}t^Y_{\epsilon}:T_Y(\epsilon )\to [0,1]$ by setting
$\pi_Y^{*}t^Y_{\epsilon}(x)=t^Y_{\epsilon}(\pi_Y(x)).$

\subsection{} 
For each stratum $Y\subset W$, the product
\begin{equation}
B_Y^{\epsilon}=(\pi_Y^{*}t^Y_{\epsilon}).(1-s_Y^{\epsilon}):T_Y(\epsilon
)\to [0,1] \end{equation}
is smooth, vanishes near $\partial T_Y(\epsilon),$ and also near $\partial \overline{Y}:$ 
\begin{equation}\label{eqn-Bpi00}
x \in \bigcup_{Z<Y} T_Z(\textstyle{\frac{\epsilon}{2}})\ \implies\
B_Y^{\epsilon}(x)=0.\end{equation}
Hence $B_Y^{\epsilon}$ admits an extension to $W$ which is defined by setting 
\begin{equation}\label{eqn-Bpi0}
B_Y^{\epsilon}(x)=0 \text{ if } x\notin T_Y(\epsilon ).\end{equation} 
\quash{
This extension is smooth on each stratum of $W$ and satisfies the following condition
whenever $Z < Y$:
\begin{equation}\label{eqn-Bpi1}
   B_Z^{\epsilon}\pi_Y(x) = B_Z^{\epsilon}(x) \text{ for all } x \in
T_Y(\epsilon).\end{equation}
Moreover $B_Y^{\epsilon}(\pi_Y(x)) \ge B_Y^{\epsilon}(x)$ for all $x\in T_Y(\epsilon_0).$
For any $x\in T_Y(\epsilon/2),$
\begin{equation}\label{eqn-Bpi2}
   B_Y^{\epsilon}\pi_Y(x)= B_Y^{\epsilon}(x)= t^Y_{\epsilon}(x)
\text{ and }
\end{equation}
\begin{equation}\label{eqn-Bpi3}
x\in T_Y(\epsilon/2) - \bigcup_{Z<Y}T_Z(\epsilon) \ \implies\ B_Y^{\epsilon}(x)=1.
\end{equation}
}  
This extension is smooth on each statum of $W$ and satisfies the following conditions
whenever $Z < Y:$ 
\begin{align}
\label{eqn-Bpi1}B_Z^{\epsilon}\pi_Y(x) = B_Z^{\epsilon}(x) &\text{ for all } x \in
T_Y(\epsilon)   \\
\label{eqn-Bpi1a}
B_Y^{\epsilon}(\pi_Y(x)) \ge B_Y^{\epsilon}(x) &\text{ for all }x\in T_Y(\epsilon)\\
\label{eqn-Bpi2}
   B_Y^{\epsilon}\pi_Y(x)= B_Y^{\epsilon}(x)= t^Y_{\epsilon}(x)
&\text{ for all } x \in T_Y(\epsilon/2) \\
\label{eqn-Bpi3}B_Y^{\epsilon}(x)=1 &\text{ for all }
x\in T_Y(\epsilon/2) - \bigcup_{Z<Y}T_Z(\epsilon).
\end{align}

\begin{lem}{}\label{lem-partition}
For every stratum $Y\subset W$ and for every point $y\in \overline{Y}$ we have
\begin{equation}\label{eqn-partition}
      B_Y^{\epsilon}(y)+\sum_{Z<Y}B_Z^{\epsilon}(y)=1\end{equation}
\end{lem}
\subsection{  Proof} It suffices to verify (\ref{eqn-partition}) for $y\in Y.$  Then
$B_Y^{\epsilon}(y) = t^Y_{\epsilon}(y).$  If $y$ is not in any tubular neighborhood
$T_Z(\epsilon )\cap Y$ (for $Z<Y)$ then $t^Y_{\epsilon}(y)=1$ and
$1-s_Z^{\epsilon}(y)=0$.  Otherise, let $\{Z_1,Z_2,\ldots Z_r\}$ be the  collection of
strata for which $Z_i<Y$ and $y\in T_{Z_i}(\epsilon)$.  By relabeling the indices, we may
assume that $Z_1<Z_2<\ldots Z_r<Y$ form a flag of strata.  The nonzero terms in the sum
(\ref{eqn-partition}) involve only the functions $s_1,s_2,\ldots,s_r$ (where $s_i =
s_{Z_i}^{\epsilon}$) and can be written:
\begin{equation}  (1-s_1)+s_1(1-s_2+s_2(\ldots +s_{r-1}(1-s_r+s_r)\ldots ))
                =1.\qed\end{equation}
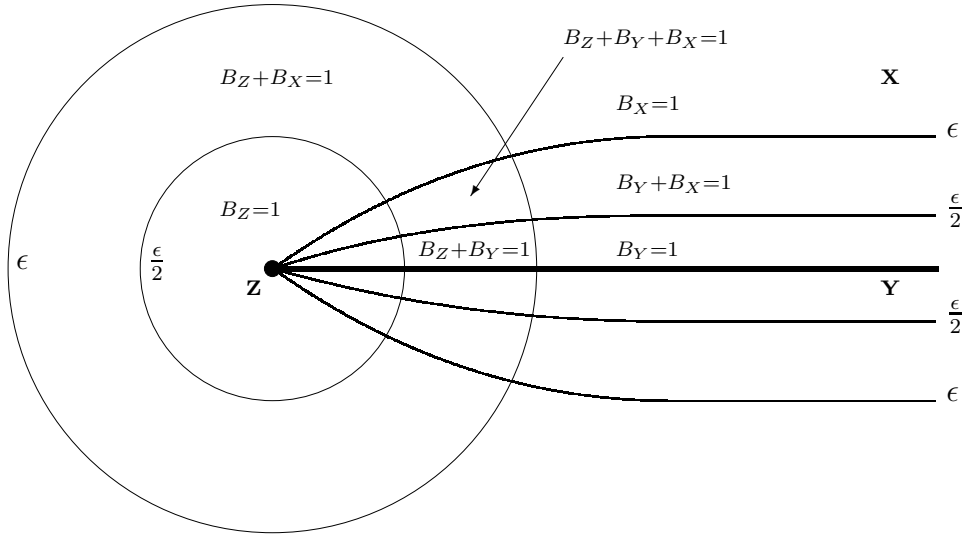
\begin{figure}[h!]\label{Fig-neighborhoods}\begin{center}
\centerline{\begin{picture}(250,180)(-10,-50)
\put(40,40){\circle*{6}}
\linethickness{2pt}
\qbezier(40,40)(110,40)(290,40)
\linethickness{.5pt}
\qbezier(40,40)(100,60)(180,60)
\qbezier(180,60)(200,60)(290,60)
\put(40,40){\circle{100}}
\qbezier(40,40)(110,20)(190,20)
\qbezier(190,20)(200,20)(290,20)
\put(40,40){\circle{200}}
\qbezier(40,40)(110,90)(190,90)
\qbezier(190,90)(200,90)(290,90)
\qbezier(40,40)(110,-10)(190,-10)
\qbezier(190,-10)(200,-10)(290,-10)
\put(30,30){$\mathbf{\scriptstyle{Z}}$}
\put(270,30){$\mathbf{\scriptstyle{Y}}$}
\put(270,110){$\mathbf{\scriptstyle{X}}$}
\put(95,45){$\scriptstyle{B_Z+B_Y=1}$}
\put(20,60){$\scriptstyle{B_Z=1}$}
\put(20,110){$\scriptstyle{B_Z+B_X=1}$}
\put(170,45){$\scriptstyle{B_Y=1}$}
\put(170,70){$\scriptstyle{B_Y+B_X=1}$}
\put(170,100){$\scriptstyle{B_X=1}$}
\put(150,125){$\scriptstyle{B_Z+B_Y+B_X=1}$}
\put(150,120){\vector(-2,-3){35}}
\put(295,60){${\frac{\epsilon}{2}}$}
\put(295,90){${\epsilon}$}
\put(295,20){$\frac{\epsilon}{2}$}
\put(295,-10){$\epsilon$}
\put(-7,40){${\frac{\epsilon}{2}}$}
\put(-57,40){${\epsilon}$}
\end{picture}}\caption{Partition of Unity for fixed $\epsilon$}\end{center}
\end{figure}

\subsection{}\label{subsec-family}
In \S \ref{subsec-patched} we will construct a connection on a modular variety by
induction, patching together connections which have been previously defined on
neighborhoods of boundary strata.  For each step of this induction we will need a
different partition of unity, which is obtained from (\ref{eqn-partition}) by shrinking
the parameter $\epsilon.$  The purpose of this subsection is to construct the family of
partitions of unity.

Let $W$ be a weakly stratified space with an $\epsilon_0 > 0$ system of control data,
$\{T_Y(\epsilon_0), \pi_Y, \rho_Y\}.$  Suppose each stratum $Y$ is a complex manifold and
define 
\[ \epsilon Y = \epsilon_0/2^{\dim_{\mathbb C}(Y)}. \]
(The complex structure is irrelevant to this construction and is only introduced so as to
agree with the later sections in this paper.)  By Lemma \ref{lem-partition} for every
point $x\in Y$ we have
\[ \sum_{Z \le Y} B_Z^{\epsilon Y}(x) = 1.\]
By (\ref{eqn-Bpi1}) and (\ref{eqn-Bpi2}) the same equation holds, in fact for all 
$x\in T_Y({\frac{\epsilon Y}{2}}).$  

\quash{ From (\ref{eqn-Bpi00}) and (\ref{eqn-Bpi0}) we find:
\begin{lem}
Suppose that $Y<X$ are strata of $W$, that $x\in T_Y(\frac{\epsilon X}{2})$ and that
$B_Y^{\epsilon Y}(x) \ne 0.$  Then $B_Y^{\epsilon X}(x) = 1. $ \qed 
\end{lem}
} 

Let $x\in W.$  Then there is a maximal collection of strata $Y_1,Y_2,\ldots,Y_r$ such
that 
\[ x\in T_{Y_1}(\epsilon_0)\cap T_{Y_2}(\epsilon_0) \cap \ldots \cap T_{Y_r}(\epsilon_0).
\]
These strata form a partial flag which (we may assume) is given by $Y_1<Y_2<\cdots <Y_r.$
Set
\[ B_n^{\epsilon m} = B_{Y_n}^{\epsilon Y_m} \text{ and } \pi_m(x) = \pi_{Y_m}(x).\]  
In the following lemma we assume $1 \le m,n,m',n' \le r,$ that $m \ge n$ and that $m' \ge
n'.$  
\begin{lem} \label{lem-vanishing}
 If $B_n^{\epsilon m}(\pi_n(x)) \ne 0$ then 
$B_{n'}^{\epsilon m'}(x)=0$ for all $n'<n$ and for all $m'>m.$
If $B_n^{\epsilon m}(x) \ne 0$ then
$B_{n'}^{\epsilon m'}(\pi_{m'}(x)) = 0$ for all 
$n' > n$ and for all $m' < m.$  
\end{lem}
\subsection{Proof}  If $B_n^{\epsilon m}(\pi_n)(x) \ne 0$ then $\pi_n(x) \notin
T_{n'}(\epsilon m/2) \supset T_{n'}(\epsilon m')$ by (\ref{eqn-Bpi00}).  Therefore
$B_{n'}^{\epsilon m'}(\pi_n(x)) = 0$ by (\ref{eqn-Bpi0}).  Hence $B_{n'}^{\epsilon
m'}(x) = 0$ by (\ref{eqn-Bpi1}).  The second statement is the contrapositive of the
first.  \qed 

\section{   $\pi$-Fiber Differential Forms }\label{sec-forms}

\subsection{  }
  Suppose $W$ is a stratified space with a fixed $\epsilon$-system of control data
$\{ T_Y(\epsilon),\pi_Y,\rho_Y\}$.  Define a $\pi$-{\it fiber differential form}
$\omega$ to be a collection $\omega =\{\omega_Y \in \mathcal{ A}^{*}(Y;\mathbb C)\}$ of
smooth differential forms (with complex coefficients) on the strata $Y$ of $W$,
which satisfy the following compatibility condition whenever $Z<Y$:
There exists a neighborhood $T(\omega )\subset T_Z(\epsilon)$ of $Z$ such that
\begin{equation}\label{eqn-pifiber}
\omega_Y|(Y\cap T(\omega ))=\pi_{YZ}^{*}(\omega_Z)|(Y\cap T(\omega
)).
\end{equation}
(Here, $\pi_{YZ}$ denotes the restriction $\pi_Z|Y \cap T_Z(\epsilon)$.)  We refer to
equation (\ref{eqn-pifiber}) as the {\em $\pi$-fiber condition}.  If $\omega
=\{\omega_Y\}$ is a $\pi$-fiber differential form, define its differential to be the
$\pi$-fiber differential form $d\omega =\{d\omega_Y\}$.  Let
$\mathcal{A}^{\bullet}_{\pi}(W)$ denote the complex of $ \pi$-fiber differential forms
and let $H^{*}_{\pi}(W)$ denote the resulting cohomology groups.  These differential
forms are analogous to the $\pi$-fiber cocycles in \cite{G}.  The following result is
proven in \cite{Verona}.

\begin{thm}\label{thm-mythesis}  The inclusion $\mathbb C\to
\mathcal{A}^{\bullet}_{\pi}(W)$ of the constant functions into the complex of $\pi$-fiber
differential forms induces an isomorphism
\begin{equation}  H^i(W;\mathbb C)\cong H^i_{\pi}(W)\end{equation}
for all $i$.  The restriction $i^*[\omega] \in H^*(\overline{Y})$ of the cohomology
class $[\omega]$ represented by a closed, $\pi$-fiber differential form $\omega =
\{\omega_Y\}_{Y\subset W}$ to the closure $\overline{Y}$ of a single stratum is given by
the $\pi$-fiber differential form $\{\omega_Z\}_{Z \subset \overline{Y}}$ (where
$i:\overline{Y} \to W$ denotes the inclusion).\end{thm}
\quash{
\subsection{  Outline of Proof  } 
  Since the $\pi$-fiber condition is local, the $\pi$-fiber differential forms
$\mathcal{A}^i_{\pi}(W)$ are the global sections of a sheaf \ ${\mathbf{
\mathcal{A}}}^i_{\pi}$ of $\pi$-fiber differential forms (which is obtained by first
restricting  $\pi$-fiber forms to open sets and then sheafifying the resulting
presheaf).  This sheaf is a module over the sheaf of $\pi$-fiber functions and hence is
fine.  Let $Z$ be a stratum in $W$, fix $x\in Z$ and let $B\subset Z$ be a small open
ball containing the point $x$.	Let
\begin{equation}  U=\pi_Z^{-1}(B)=\coprod_{Y\geq Z}\pi_{YZ}^{-1}(B)\end{equation}
be the resulting ``basic'' open neighborhood of $x\in W$ (\cite{Mather}, \cite{GM1} \S
1.1, \cite{GM2} \S 1.1, \cite{GM3} \S 1.4 ).  If $\omega=\{\omega_Y\}$ is a differential
form which satisfies the $\pi$-fiber condition (\ref{eqn-pifiber}) in the whole
neighborhood $U$, then it is completely determined by the smooth differential form
$\omega_Z\in \mathcal{ A}^{*}(B).$ Hence, the Poincar\'e lemma for the ball $ B\subset Z$
implies the Poincar\'e lemma for the stalk $\mathcal{A}^{*}_{\pi ,x}$ at $x$ of the
complex of $\pi$-fiber differential forms.  This shows that the inclusion of the constant
sheaf $\mathbb C_W\to {\mathbf{ \mathcal{A}}}^{*}_{\pi}$ into the complex of sheaves of $
\pi$-fiber differential forms is a quasi-isomorphism and hence induces an isomorphism on
hypercohomology.  (This corrects the proof in \cite{Verona}.) The statement about the
restriction to $H^*(\overline{Y})$ follows
immediately.  \qed } 

\subsection{}\label{subsec-resolve}
Suppose the stratified space $W$ is the closure of a single stratum $X$. Then a smooth
differential form $\omega_X \in \mathcal A^I(X)$ is the $X$-component of a $\pi$-fiber
differential form $\omega \in \mathcal A^i_{\pi}(W)$ if and only if for each stratum $Y$
there exists a neighborhood $U_Y$ of $Y$ such that for each point $p \in U_Y \cap X$ and
for every tangent vector $v\in T_pX$ the following condition holds:
\begin{equation}\label{eqn-contraction}
\text{If } d\pi(p)(v)=0 \text{ then } i_v\omega = 0
\end{equation}
where $i_v$ denotes the contraction with $v$.  If (\ref{eqn-contraction}) holds, then the
$\pi$-fiber form $\omega$ is uniquely determined by $\omega_X.$

Now suppose that $W = \overline{X}$ is a compact subanalytic Whitney stratified subset of
some (real) analytic manifold, and that $\tau:\widetilde{W} \to W$ is a (subanalytic)
resolution of singularities (cf. \cite{Hironaka1}, \cite{Hironaka2}).  This means that
$\widetilde{W}$ is a smooth compact subanalytic manifold, the mapping $\tau$ is
subanalytic, its restriction $\tau^{-1}(X) \to X$ to $\tau^{-1}(X)$ is a diffeomorphism,
and $\tau^{-1}(X)$ is dense in $\widetilde{W}.$ Let $\omega_X \in \mathcal{A}^i(X)$ be a
$\pi$-fiber differential form on $W$.

\begin{lem}{}\label{lem-integrate}
Suppose the differential form $\tau^{*}(\omega_X)$ is the restriction of a smooth closed
differential form $\widetilde{\omega} \in \mathcal{A}^i(\widetilde{W}).$  Let
$[\widetilde{\omega}] \in H^i(\widetilde{W};\mathbb R)$ and $[\omega_X]\in H^i(W;\mathbb
R)$ denote the corresponding cohomology classes.  Then $[\widetilde{\omega}] =
\tau^*([\omega_X]).$ \end{lem}

\subsection{Proof}  The cohomology classes $[\widetilde{\omega}]$ and $[\omega_X]$ are
determined by their integrals over subanalytic cycles by \cite{Hardt1}, \cite{Hardt2}.
Any subanalytic cycle $\xi \in C_i(\widetilde{W};\mathbb R)$ may be made transverse
(within its homology class) to the ``exceptional divisor'' $\tau^{-1}(W-X)$ (\cite{GM3}
\S 1.3.6). Then
\[ \int_{\xi}\widetilde{\omega} = \int_{\xi \cap \tau^{-1}(X)}\widetilde{\omega}=
\int_{\tau(\xi)\cap X}\omega_X = \int_{\tau(\xi)}\omega_X \qed \]


\section{   Homogeneous Vectorbundles   }\label{sec-vb}
\subsection{   }
If $M$ is a smooth manifold and $E \to M$ is a smooth vectorbundle, let $\mathcal
A^i(M,E)$  denote the space of smooth differential $i$-forms with values in $E$.
Throughout this section, $K$ denotes a closed subgroup of a connected Lie group $G$ with
Lie algebras $\mathfrak k\subset \mathfrak g$ and with quotient $D=G/K.$  We fix a
representation $\lambda:K \to GL(V)$ on some finite dimensional (real or complex)
vectorspace $V.$ Write $\lambda':\mathfrak k \to \text{End}(V)$ for its derivative at the
identity, and note that its derivative at a general point $g\in K$ is given by 
\begin{equation}\label{eqn-Lg}
d\lambda (h)((L_h)_{*}(\dot k ))=\lambda (h)\lambda'(\dot k )\in\text{\rm End}
(V)\end{equation} for any $\dot k \in \mathfrak k$, where $L_h:K \to K$ is multiplication
from the left by $h \in K.$ The quotient mapping $q:G \to D=G/K$ is a principal
$K$-bundle.  The {\it fundamental vertical vectorfields} $Y_{\dot k}(g) = L_{g*}(\dot
k)$ (for $\dot k \in \mathfrak k$) determine a canonical trivialization, $\ker (dq) \cong
G \times \mathfrak k.$

\subsection{}\label{subsec-VB}
The representation $\lambda:K \to GL(V)$ determines an associated
homogeneous vectorbundle $E = G \times_K V$ on $D = G/K,$ which consists of equivalence
classes $[g,v]$ of pairs $(g,v) \in G \times V$ under the equivalence relation $(gk,v)
\sim (g, \lambda(k)v)$ for all $g\in G$, $k\in K$, $v\in V.$  It admits the homogeneous
$G$ action given by $g' \cdot [g,v] = [g'g,v].$  Smooth sections $\tilde s$ of $E$ may be
identified with smooth mappings $s:G \to V$ such that 
\begin{equation}\label{eqn-rightk}
s(gk)= \lambda(k^{-1})s(g) \end{equation} 
by $\tilde s(gK) = [g,s(g)]\in E.$ Then (\ref{eqn-rightk}) implies
\begin{equation}\label{eqn-vert}
ds(g)(L_{g*}(\dot k)) = -\lambda'(\dot k)s(g) \end{equation}
for all $g\in G$ and $\dot k \in \mathfrak k.$ 

Similarly we identify smooth differential forms $\tilde\eta \in \mathcal A^i(D,E)$ (with
values in the vectorbundle $E$) with  differential forms $\eta \in \mathcal
A^i_{\text{bas}}(G,V)$ which are ``basic,'' meaning they are both {\it $K$-equivariant}
($R_k^*(\eta) = \lambda(k)^{-1}\eta$ for all $k \in K$) and {\it horizontal} ($i(Y_{\dot
k})\eta = 0 \text{ for all } \dot k \in \mathfrak k$).  Here, $i(Y)$ denotes the interior
product with the vectorfield $Y$ and $R_k(g) = gk$ for $g\in G.$

A connection $\nabla$ on $E$ is determined by a connection $1$-form $\omega \in \mathcal
A^1(G, \text{End}(V))$ which satisfies $\omega(L_{g*}({\dot k})) = \lambda'(\dot k)$ and
$R_k^*(\omega) = \text{Ad}(\lambda(k^{-1}))\omega $ for any $\dot k\in\mathfrak k$, $k
\in K$, and $g\in G.$  The covariant derivative $\nabla_X$ with respect to a vectorfield
$X$ on $D$ acts on sections $s:G \to V$ satisfying (\ref{eqn-rightk}) by 
\begin{equation}\label{eqn-covariant}
\nabla_Xs(g) = ds(g)(\tilde X(g)) +  \omega_g(\tilde X(g))(s(g))
\end{equation} where $\tilde X$ is any lift of $X$ to
a smooth vectorfield on $G$.  We write $\nabla = d+\omega.$  
The curvature form  $\Omega \in \mathcal A^2(D, \text{End}(E))$ takes values in the
vectorbundle $\text{End}(E) = G \times_K \text{End(V)}$ and it will be identified
with the ``basic'' 2-form $\Omega \in \mathcal A^2_{\text{bas}}(G,\text{End}(V))$ which
assigns to tangent vectors $X,Y\in T_gG$ the endomorphism 
\begin{equation}\label{eqn-structure}
\Omega(X,Y) = d\omega(X,Y) + [\omega(X),\omega(Y)]
\end{equation}
If the Lie bracket is extended in a natural way to Lie algebra-valued
$1-$forms, then it turns out (\cite{Berline} \S 1.12) that $[\alpha,\alpha](X,Y) =
2[\alpha(X),\alpha(Y)]$ so we may express the curvature form as $\Omega = d\omega +
\frac{1}{2}[\omega,\omega].$

The connection $\nabla = d+\omega$ is $G$-invariant iff $L_g^*(\omega)=\omega$, in which
case it is determined (on the identity component $G^0$) by its value $\omega_0:\mathfrak
g \to \text{End}(V)$ at the identity.  Using \cite{Wang} or \cite{KN} Chapt. II Thm.
11.5. it is easy to verify the following result. \begin{prop}\label{prop-classify}
Suppose $G$ is a connected Lie group and $K$ is a closed subgroup.  Then the
$G$-invariant connections on the homogeneous vectorbundle $E = G\times_KV$ are given by
linear mappings $\omega_0:\mathfrak g \to \text{End}(V)$ such that
\begin{enumerate}
\item $\omega_0(\dot k) = \lambda'(\dot k)$ for all $\dot k \in \mathfrak k$ and
\item $\omega_0([\dot g, \dot k]) = [\omega_0(\dot g), \lambda'(\dot k)]$ for all $\dot g
\in
\mathfrak g$ and all $\dot k \in \mathfrak k.$ \end{enumerate} 
Moreover the curvature $\Omega\in \mathcal A^2(G, \text{End}(V))$ of such a connection is
the left-invariant ``basic'' differential form whose value $\Omega_0$ at the identity is
given by
\begin{equation}\label{eqn-invcurv}
\Omega_0(\dot g,\dot h) = [\omega_0(\dot g),\omega_0(\dot h)] - \omega_0([\dot g,\dot h])
\end{equation} 
for any $\dot g,\dot h \in \mathfrak g.$  The connection is flat iff $\omega_0$ is a Lie
algebra homomorphism.\qed \end{prop}

\subsection{Example}  Suppose the representation $\lambda:K \to GL(V)$ is the restriction
of a representation $\tilde\lambda:G \to GL(V)$.  Then we obtain a  {\it flat} connection
with $\omega_0(\dot g) = \tilde\lambda'(\dot g)$  for $\dot g \in \mathfrak g$.  

\subsection{Example}\label{ex-Nomizu}  A connection in the principal bundle $G \to D$ is
given by a connection 1-form $\theta\in\mathcal A^1(G,\mathfrak k)$ such that
$R_k^*(\theta) = Ad(k^{-1})(\theta)$ and $\theta(Y_{\dot k}) = \dot k$ (for any $k \in K$
and any fundamental
vectorfield $Y_{\dot k}$).  It determines a connection  $\nabla = d + \omega$ in the
associated bundle $E = G\times_KV$ by $\omega(X) = \lambda'(\theta(X)).$  The principal
connection $\theta$ is $G$-invariant iff $L_g^*(\theta)=\theta$ in which case it is
determined by its value $\theta_0:\mathfrak g \to \mathfrak k$ at the identity.
Conversely, by \cite{Nomizu}, any linear mapping $\theta_0:\mathfrak g \to \mathfrak k$
determines a $G$-invariant principal connection iff $\ker(\theta_0)$ is preserved under
the adjoint action of $K$.  If $K$ is a maximal compact subgroup of $G$ then the Cartan
decomposition $\mathfrak g = \mathfrak k \oplus \mathfrak p$  determines
a  canonical $G$ invariant connection in the principal bundle $G \to G/K,$ and hence a
connection on $E$ which we refer to as the {\it Nomizu connection.}   It is given by
$\omega_0(\dot k + \dot p) = \lambda'(\dot k)$ for $\dot k \in \mathfrak k$ and $\dot p
\in \mathfrak p$. By (\ref{eqn-invcurv}), its curvature is given by
\begin{equation}\label{eqn-Nomizu2}
\Omega(L_{g*}(\dot g_1),L_{g*}(\dot g_2)) = - \lambda'([\dot p_1,\dot p_2])
\end{equation}
(where $\dot g_i = \dot k_i + \dot p_i \in \mathfrak k \oplus \mathfrak p$), 
since $[\mathfrak k, \mathfrak p] \subset \mathfrak p$ and $[\mathfrak p, \mathfrak p]
\subset \mathfrak k.$

\quash{
\subsection{   } Let $\lambda:K\to GL(V)$ be a representation.
Suppose $\rho:H \to GL(V)$ is a representation of another Lie group $H$ on the
same vectorspace $V$.  Suppose that $\rho$ and $\lambda$ commute in the sense
that for any $h\in H$ and for any $k \in K$ we have
\begin{equation}  \rho(h)\lambda(k) = \lambda(k)\rho(h)\in GL(V).\end{equation}
Then $H$ acts on $E=G\times_KV$ by vectorbundle automorphisms which cover the
identity mapping on $D=G/K$ with $h.[g,v] = [g,\rho(h)v].$
\begin{lem}\label{lem-invariance} Suppose the representations $\lambda:K\to GL(V)$ and
$\rho:H \to GL(V)$ commute.  Then for any principal connection $\theta\in\mathcal
A^1(G,\mathfrak k)$ on $G\to G/K$ the resulting connection $\nabla$ on
$E=G\times_KV$ is $H$-invariant. \end{lem}
\subsection{  Proof }  The group $H$ acts on sections $s:G\to V$ of $E$ by
$h\cdot s = \rho(h)s.$  Now compute
\begin{equation}  \begin{split}
\nabla_X(h\cdot s) &= d_X(\rho(h)s) + \lambda'(\theta(X))(\rho(h)s) \\
  &= \rho(h)(d_Xs + \lambda'(\theta(X))s) = h\cdot\nabla_Xs\ \qed \end{split}
\end{equation}
} 

\subsection{  }\label{subsec-automorphy}
There is a further description of homogeneous vectorbundles on $D$ which are 
topologically trivial.   Let $E=G\times_KV$ be a homogeneous vectorbundle
corresponding to a representation $\lambda :K\to GL(V)$ of $K$ .  A (smooth)
{\it automorphy factor} $J:G\times D\to GL(V)$ for $E$ (or for $\lambda$) is a 
(smooth) mapping such that \begin{enumerate}
\item\label{item-aut1} $J(gg',x)=J(g,g'x)J(g',x)$ for all $g,g'\in G$ and for
all $x\in D$
\item\label{item-aut2} $J(k,x_0)=\lambda (k)$ for all $k\in K.$ \end{enumerate}
It follows (by taking $g=1$) that $J(1,x)=I$.  The automorphy
factor $J$ is determined by its values $J(g,x_0)$ at the basepoint:  any
smooth mapping $j:G\to GL(V)$ such that $j(gk) = j(g)\lambda(k)$ (for all
$k\in K$ and all $g\in G$) extends in a unique way to an automorphy factor
$J:G\times D \to GL(V)$ for $E$, namely 
\begin{equation}\label{eqn-aut3}J(g,hx_0) = j(gh)j(h)^{-1}.\end{equation}

An automorphy factor $J$, if it exists, determines a (smooth) trivialization
\begin{equation}\label{eqn-Jtrivialization}
\Phi_J:G\times_KV\to (G/K)\times V\end{equation}
by $[g,v]\mapsto (gK,J(g,x_0)v)$.  This trivialization is $G$-equivariant with respect to
the following {\it $J$-automorphic action} of $G$ on $(G/K)\times V$:
\begin{equation}\label{eqn-Jaction}  g.(x,v)=(gx,J(g,x)v).\end{equation}
Conversely, any smooth trivialization $\Phi :E\cong (G/K)\times
V$ of the homogeneous vectorbundle $E$ determines a unique automorphy factor
$J$ such that $\Phi =\Phi_J$.  A trivialization of $E$ (if one exists) allows one to
identify smooth sections $s$ of $E$ with smooth mappings $r:D\to V$.  If the
trivialization is given by an automorphy factor $J$ and the smooth section $s$ of $E$ is
given by a $K$-equivariant mapping $s:G \to V$ as in (\ref{eqn-rightk}) then the
corresponding smooth mapping is
\begin{equation}\label{eqn-J1}  r(gK)=J(g,x_0)s(g)\end{equation}
which is easily seen to be well defined.  Sections $s$ which are invariant under $\gamma
\in G$ correspond to functions $r$ such that $r(\gamma x) = J(\gamma,x)r(x)$ for all
$x\in D.$

If $D= G/K$ is Hermitian symmetric of noncompact type then the {\it canonical automorphy
factor} (\S\ref{subsec-canon-automorphy}) $J_0:G \times D \to \mathbf K(\mathbb C)$
determines an automorphy factor $J = \lambda_{\mathbb C} \circ J_0$ for every homogeneous
vectorbundle $E = G \times_K V,$ where $\lambda_{\mathbb C}: \mathbf K(\mathbb C) \to
\text{GL}(V)$ is the complexification of $\lambda.$

\subsection{}\label{subsec-twoauto}
 If $J_1,J_2$ are two automorphy factors for $E$, then the
mapping $\phi:D\times V \to D\times V$ which is given by $\phi(gx_0,v) =
(gx_0, J_1(g,x_0)J_2(g,x_0)^{-1}v)$ is a well defined $G$-equivariant
isomorphism of trivial bundles, where $G$ acts on the domain via the
$J_1$-automorphic action and $G$ acts on the target via the $J_2$-automorphic
action.

\subsection{  }(The following fact will be used in \S \ref{subsec-proof-induced}.)
Suppose that $J:G\times D \to GL(V)$ is an automorphy factor for $E$.  Let
$\nabla=d+\omega$ be a connection on $E.$  The trivialization
$\Phi_J$ of $E$ (\ref{eqn-Jtrivialization}) determines a connection
$\nabla^J = d + \eta$ on the trivial bundle $D\times V$ (with $\eta\in
\mathcal{A}^1(D,\text{End}(V))$) as follows: if $r(gK) = J(g,x_0)s(g)$ as in
(\ref{eqn-J1}) then
$\nabla^J_{q_*X}r(gK) =  J(g,x_0) \nabla_Xs(g)$ (for any $X \in T_gG$).  It
follows from (\ref{eqn-Lg}) that the connection 1-forms are related by
\begin{equation}\label{eqn-related}  
\eta(q_*(X)) = J(g,x_0)\omega(X)J(g,x_0)^{-1} - (d_XJ(g,x_0))
J(g,x_0)^{-1}.\end{equation}

\section{   Lemmas on curvature  }\label{sec-curvature}
\subsection{ }
Suppose $E = G \times_K V$ is a homogeneous vectorbundle over $D = G/K$
arising from a representation $K \to GL(V)$ of a closed subgroup $K$ of a Lie group
$G$ on a complex vectorspace $V$.  If $\theta\in \mathcal A^1(G,\text{End}(V))$ is a Lie
algebra-valued 1-form, denote by $[\theta,\theta]$ the Lie algebra valued 2-form
$(X,Y)\mapsto 2[\theta(X),\theta(Y)]$ (cf \cite{Berline} \S 1.12).  The proof of the
following lemma is a direct but surprisingly tedious computation.

\begin{lem}\label{lem-patch} Suppose $\{f_1,f_2,\ldots,f_n\}$ form a smooth partition of
unity on $D$, that is $0 \le f_i(x) \le 1$  and $\sum_{i=1}^nf_i(x) =1$ for all $x\in D.$
Let $\nabla_i = d+\omega_1$ be connections on $E$ with curvature forms $\Omega_i$ for
$1 \le i \le n.$  Let $\nabla = \sum_{i=1}^n f_i \nabla_i$ be the connection with
connection form $\omega = \sum_{i=1}^n f_i \omega_i.$  Then the curvature $\Omega$ of
$\nabla$ is given by
\[
\Omega = \sum_{i=1}^n f_i\Omega_i -{\textstyle\frac{1}{2}}
\sum_{i<j}f_if_j[\omega_i-\omega_j,
\omega_i-\omega_j] + \sum_{i=1}^{n-1} df_i\wedge(\omega_i-\omega_n).  \qed \]
\end{lem}
(Even if the $\nabla_i$ are flat and the $f_i$ are constant, the
connection $\nabla$ is not necessarily flat.)

\subsection{  }\label{subsec-characteristicforms}
Let $\mathcal E$ be a complex vectorspace and let $f:\mathcal E \to \mathbb C$ be a
homogeneous polynomial of degree $k$.  The {\it polarization} of $f$ is the unique
symmetric $k$-linear form $P:\mathcal E\times \mathcal E\times\ldots \times \mathcal E
\to \mathbb C$ such that $f(x) = P(x,x,\ldots,x)$ for all $x\in \mathcal E.$  If
$N\subset \mathcal E$ is a vectorsubspace such that $f(x+n)=f(x)$ for all $x\in \mathcal
E$ and all $n\in N$ then the polarization $P$ satisfies
\begin{equation}\label{eqn-polarization}
P(x_1+n_1,x_2+n_2,\ldots,x_k+n_k) = P(x_1,x_2,\ldots,x_k) \end{equation}
for all $x_1,\ldots,x_k\in \mathcal E$ and all $n_1,\ldots,n_k \in N.$

Now let $K\to GL(V)$ be a representation on a complex vectorspace $V$ as above, and let
$f:\mathcal E=\text{End}(V) \to \mathbb C$ be a homogeneous polynomial of degree $k$,
which is invariant under the adjoint action $K \to GL(\mathfrak k) \to GL(\text{End}(V))$
of $K$. Then the  polarization $P$ of $f$ is also $\text{Ad}_K$-invariant.  If $\nabla$
is a  connection on $E = G\times_KV$ with curvature form $\Omega \in \mathcal A^2(G,
\text{End}(V))$
then the characteristic form associated to $f$ is
\begin{equation}  \sigma^f_{\nabla}(X_1,Y_1,\ldots,X_k,Y_k) = P(\Omega(X_1,Y_1),
\ldots,\Omega(X_k,Y_k))\in \mathcal{A}^{2k}(G,\mathbb C).\end{equation}
It is ``basic'' (cf. \S \ref{subsec-VB}) and hence descends uniquely to a differential
form on $D$ which we also denote by $\sigma^f_{\nabla} \in \mathcal A^{2k}(D,\mathbb C).$
Throughout this paper we shall be concerned only with homogeneous polynomials
$f:\text{End}(V)
\to \mathbb C$ which are invariant under the full adjoint action of $GL(V)$ and which we
shall refer to as a {\it {\rm Ad}-invariant polynomials}. 
When $f(x)$ is the $i$-th elemen\-tary symmetric function in the eigenvalues
of $x$, the resulting characteristic form is the {\it $i$-th Chern form} and it will be
denoted $\sigma^i(\nabla).$

\begin{lem}\label{lem-Springer}  
Let $V$ be a complex vectorspace, $H =  \text{\rm GL}(V)$ and $\mathfrak h = 
\text{\rm End}(V).$  Let $x,n\in \mathfrak h$ and suppose that $[x,n] = 0$ and that $n$
is nilpotent.  Then for any
{\rm $Ad$}-invariant polynomial $f:V \to \mathbb C$ we have:
\[ f(x+n) = f(x). \] \end{lem}

\subsection{Proof}  Let $\mathfrak t \subset \mathfrak h$ be the Lie algebra of a Cartan
subgroup $T \subset H$ and let $W$ denote its Weyl group.  The {\it adjoint quotient}
 mapping $\chi:\mathfrak h \to \mathfrak t/W$ associates to any $a\in \mathfrak h$ the
$W$-orbit $C(a_s)\cap \mathfrak t$ where $a = a_s + a_n$ is the Jordan
decomposition of $a$ into its semisimple and nilpotent parts (with
$[a_s,a_n]=0$), and where $C(a_s)$ denotes the conjugacy class of $a_s$ in
$\mathfrak h,$ cf. \cite{Springer}.  (The value $\chi(a)$ may be interpreted as the
coefficients of the characteristic polynomial of $a \in \text{End}(V).$)  Since $x$ and
$n$ commute, $(x+n)_s = x_s + n_s = x_s,$ hence $\chi(x+n) = \chi(x).$  But every
Ad-invariant polynomial $f:\mathfrak h \to \mathbb C$ factors through $\chi$, hence
$f(x+n)=f(x).$ \qed


\section{  Hermitian Symmetric Spaces  }\label{sec-hermitian}
\subsection{   }\label{subsec-parabolics}
Throughout the remainder of this paper, algebraic groups will be denoted by boldface type
and the associated group of real points will be denoted in Roman.  Fix a reductive
algebraic group $\mathbf G$ which is defined over $\mathbb Q$ and which (for convenience
only) is assumed to be connected and simple over $\mathbb Q.$    Suppose $G$ acts
as the identity component of the group of automorphisms of a Hermitian symmetric
space $D$.  A choice of basepoint $x_0 \in D$ determines a Cartan involution $\theta:G
\to G$, a maximal compact subgroup $K =G^{\theta}=  \text{Stab}_G(x_0)$ and a
diffeomorphism $G/K \cong D.$  Let $D^{*}$ denote the {\it Satake partial
compactification of $D$}, in other words, the union of $D$ and all its rational boundary
components, with the Satake topology (cf. \cite{BB}).  The closure $D_1^*$ in $D^*$ of a
rational boundary component $D_1$ is again the Satake partial compactification of $D_1.$

If $D_1\subset D^{*}$ is a rational boundary component of $D$, its normalizer
$P$ is the group of real points of a rationally defined maximal proper parabolic
subgroup $\mathbf P \subset \mathbf G$.  Let $\mathcal{U}_P$ denote the unipotent radical
of $P$ and let $L(P) = P/\mathcal{U}_P$ be the Levi quotient with projection $\nu_P:P\to
L(P)$.  It is well known that $L(P)$ is an almost direct product (commuting product with
finite intersection) of two subgroups, $L(P) = G_hG_{\ell}$ (we include  possible compact
factors in the $G_{\ell}$ factor), where the ``hermitian
part'' $G_h$ is semisimple and defined over $\mathbb Q.$  (It may be trivial).  The
choice of basepoint $x_0\in D$ determines a unique $\theta$-stable lift $L_P(x_0)\subset
P$ of the Levi quotient \cite{BorelSerre} (which, from now on, we shall use without
mention), as well as basepoints $x_1 \in D_1$ in the boundary component $D_1$ such that
$\text{Stab}_{G_h}(x_1) = K\cap G_h.$  We obtain a decomposition
\begin{equation}\label{eqn-decompose1}
P = \mathcal U_P G_h G_{\ell}. \end{equation}
The group $P$ acts on the boundary component $D_1$ through the projection
$\nu_h:P \to L(P) \to G_h/(G_h\cap G_{\ell})$ which also determines a diffeomorphism
$G_h/K_h \cong D_1$ (where $K_h = K \cap G_h$).  This projection $\nu_h$ also gives
rise to a $P$-equivariant {\it canonical projection}
\begin{equation}\label{eqn-pi}  
\pi: D \to D_1\end{equation}
by $\pi(ug_hg_{\ell}K_P) = g_hK_h.$

The ``linear part'' $G_{\ell}$ is reductive and contains the 1-dimensional
$\mathbb Q$-split torus $\mathbf{S_P}(\mathbb R)$ in the center of the Levi
quotient. If $\mathfrak z \subset \mathfrak N_P$ denotes the center of the Lie
algebra $\mathfrak N_P$ of the unipotent radical of $P$, then the adjoint action of
$G_{\ell}$ on $\mathfrak z$ has a unique open orbit $C(P)$ which is a self adjoint
homogeneous cone.

\subsection{   }
Let $\mathbf{P_0} \subset \mathbf G$ be a fixed minimal rational parabolic subgroup and
define the standard parabolic subgroups to be those which contain $\mathbf{P_0}.$  Let
$\mathbf{S_0} \subset L(\mathbf{P_0})$ be the greatest $\mathbb Q$-split torus in the
center of $L(\mathbf{P_0})$ and let $\Phi=\Phi(\mathbf{S_0},G)$ be the (relative) roots
of $G$ in $\mathbf{S_0}$ with positive roots $\Phi^+$ consisting of those roots which
appear in the unipotent radical $\mathcal{U}_{P_0}$ and resulting simple roots $\Delta.$
Each simple root $\alpha$ corresponds to a vertex in the (rational) Dynkin diagram for
$\mathbf G$ and also to a maximal standard parabolic subgroup $\mathbf P$ such that
$\mathbf {S_P} \subset \ker(\beta)$ for each $\beta \in \Delta - \{ \alpha \}.$
\quash{
\begin{equation}  \mathbf{S_P} = \left( \bigcap_{{\beta \in \Delta -
\{\alpha\}
}} \ker (\beta)\right)^0\subset \mathbf{S_0}
\end{equation}
} 

\subsection{  Two maximal parabolic subgroups  }\label{subsec-twoparab}
\begin{figure}[h!]\begin{center}
\begin{picture}(400,130)
\thicklines

\put(85,120){$\alpha_1$}  \put(295,120){$\alpha_2$}

\put(0,100){\circle{15}} \put(30,100){\circle{15}} \put(60,100){\circle{15}}
\put(90,100){\circle*{15}} \put(120,100){\circle{15}} \put(150,100){\circle{15}}
\put(180,100){\circle{15}} \put(210,100) {\circle{15}} \put(240,100){\circle{15}}
\put(270,100){\circle{15}}  \put(300,100){\circle*{15}} \put(330,100){\circle{15}}
\put(360,100){\circle{15}} 

\put(7,102.5){\line(1,0){16}}
\put(7, 97.5){\line(1,0){16}}

\put(37.5,100){\line(1,0){15}} 
\put(67.5,100){\line(1,0){15}}
\put(97.5,100){\line(1,0){15}}
\put(127.5,100){\line(1,0){15}}
\put(157.5,100){\line(1,0){15}}
\put(187.5,100){\line(1,0){15}}
\put(217.5,100){\line(1,0){15}}
\put(247.5,100){\line(1,0){15}}
\put(277.5,100){\line(1,0){15}}
\put(307.5,100){\line(1,0){15}}
\put(337.5,100){\line(1,0){15}}

\put(0,50){\circle{15}} \put(30,50){\circle{15}} \put(60,50){\circle{15}}
\put(120,50){\circle{15}} \put(150,50){\circle{15}}
\put(210,50) {\circle{15}} \put(180,50){\circle{15}}
\put(240,50){\circle{15}} \put(270,50){\circle{15}}
\put(330,50){\circle{15}} \put(360,50){\circle{15}} 

\put(7,52.5){\line(1,0){16}} \put(7,47.5){\line(1,0){16}} 

\put(37.5,50){\line(1,0){15}} 
\put(127.5,50){\line(1,0){15}}
\put(157.5,50){\line(1,0){15}}
\put(187.5,50){\line(1,0){15}}
\put(217.5,50){\line(1,0){15}} \put(247.5,50){\line(1,0){15}}
\put(337.5,50){\line(1,0){15}}

\put(25,25){$G_{1h}$}    \put(185,25){$G'_{\ell}$}
\put(340,25){$G_{2\ell}$}
\end{picture}\vskip-.2in
\caption{Dynkin diagrams for $G$ and $P$}\label{fig-Dynkin1}\end{center}
\end{figure}
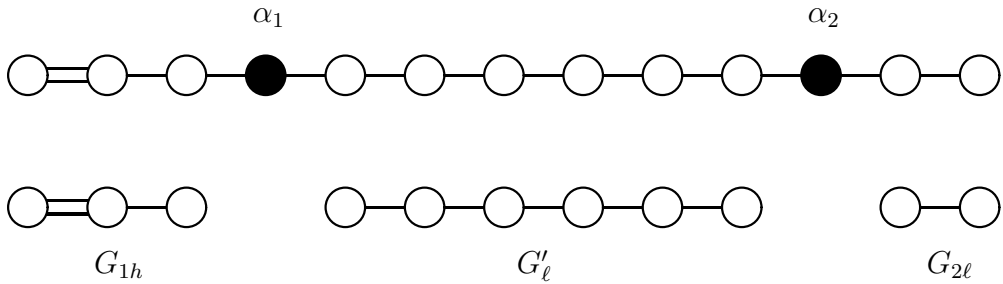

For simplicity, let us assume that $\mathbf G$ is (almost) simple over $\mathbb Q.$ 
The (rational) Dynkin diagram for $\mathbf G$ is linear, of type BC, and determines a
canonical ordering among the maximal standard rational
parabolic subgroups with $\mathbf{P_1} \prec \mathbf{P_2}$ iff $D_1 \prec D_2$ (meaning
that $D_1 \subset {D_2^*}\subset D^{*}$) where $D_i$ is the rational boundary component
fixed by ${P_i}.$  Write $P_1 = \mathcal{U}_1 G_{1h}G_{1\ell}$  and $P_2 = \mathcal{U}_2
G_{2h}G_{2\ell}$ as in (\ref{eqn-decompose1}).  If $\mathbf{P_1} \prec \mathbf{P_2}$ then
$G_{1h}\subset G_{2h}$ and $G_{1\ell} \supset G_{2\ell}.$      Let $\mathbf{P}=
\mathbf{P_1} \cap \mathbf{P_2}.$  In Figure \ref{fig-Dynkin1}, $\Delta-\{\alpha_1\}$ and
$\Delta-\{\alpha_2\}$ denote the simple roots corresponding to $P_1$ and $P_2$
respectively. 

Then we have a commutative diagram
\begin{equation}  \begin{CD}
\mathcal{U}_1 & \ \subset\  & \mathcal{U}_P &\ \subset\  & P &\ \subset\  & P_1 \\
&&@VVV @VVV  @VV{\nu_1}V \\
&&\mathcal{U}_{\overline{P}} &\ \subset\  & \overline P &\ \subset\  & L(P_1)&
=G_{1h}G_{1\ell} \\
\end{CD} \end{equation}
where $\overline{P} =\nu_1(P)\subset L(P_1)$ is the image of $P$.  Then
$\overline{ P} =
G_{1h}P_{\ell}$ where $P_{\ell}\subset G_{1\ell}$ is a parabolic subgroup of
$G_{1\ell}$ whose Levi factor decomposes as a commuting, almost direct product
$L(P_\ell)=G_{\ell}'G_{2\ell}.$  Writing $\overline{\mathcal{U}}$ for the lift of
$\mathcal{U}_{\overline P}= \mathcal U_{P_{\ell}}$ we conclude that $P$ has a
decomposition
\begin{equation}\label{eqn-twopa}
P = \mathcal{U}_1 G_{1h}  (\overline{\mathcal{U}} G_{\ell}' G_{2\ell})
= \mathcal U_P G_{1h} G_{P\ell}\end{equation}   
with $\mathcal U_P = \mathcal U_1\overline{\mathcal U}$, $G_{P\ell} =
G_{\ell}'G_{2\ell}$ and $P_{\ell} =\overline{\mathcal{U}} G_{\ell}'G_{2\ell} \subset
G_{1\ell}.$  Similarly, we have a diagram
\begin{equation}  \begin{CD}
\mathcal{U}_2 & \ \subset\  & \mathcal{U}_P & \ \subset\  & P & \ \subset\ & P_2 \\
&& @VVV @VVV @VV{\nu_2}V  \\
&& \mathcal{U}_{\overline{\overline{P}}} &\ \subset\  & \overline{\overline{P}}
&\ \subset\  & L(P_2)& = G_{2h}G_{2\ell} \end{CD} \end{equation}
where $\overline{\overline{P}} = \nu_2(P)\subset L(P_2)$ is the image of $P$.
Then $\overline{\overline{P}} =P_hG_{2\ell}$ with $P_h\subset G_{2h}$ a
parabolic subgroup of $G_{2h}$ whose Levi factor decomposes as a product
$L(P_h) = G_{1h}G_{\ell}'.$  Writing $\overline{\overline{\mathcal{U}}}$ for the
canonical lift of $\mathcal U_{\overline{\overline{P}}}=\mathcal{U}_{P_h}$ we obtain
another decomposition,
\begin{equation}\label{eqn-twopb}
P = \mathcal{U}_2 (\overline{\overline{\mathcal{U}}} G_{1h}G_{\ell}')G_{2\ell}
= \mathcal U_P G_{1h} G_{P\ell}\end{equation}
with $P_h = \overline{\overline{\mathcal{U}}}G_{1h}G'_{\ell}.$ 

\quash{
Henceforth we shall refer to $P_h \subset G_{2h}$ as the {\it relative parabolic
subgroup} for the pair $D_1 \subset D_2^*$, and to $G'_{\ell} \subset G_{2h} \cap
G_{1\ell}$ as the {\it relative linear factor}. 
}  

Similarly, an arbitrary standard parabolic subgroup $ Q$ may be expressed in a unique way
as an intersection $Q=P_{1}\cap P_{2} \cap \ldots \cap P_{m}$ of maximal standard
parabolic subgroups, with $\mathbf {P_{1}} \prec \mathbf{P_2} \prec \ldots \prec
\mathbf{P_m}.$   In this case we write $Q^{\flat}=P_1$ and $Q^{\sharp} = P_m.$  If $P_1 =
\mathcal U_1 G_{1h}G_{1\ell}$ (with projection $\nu_1:P_1 \to G_{1h}G_{1\ell}$) then the
Levi factor $L(Q)$ decomposes as an almost direct product of $m+1$ factors
\begin{equation}\label{eqn-manyLevi}
L(Q) = G_{1h} (G'_1 G'_2 \ldots G'_{m-1} G_{m\ell}) = G_{1h}G_{Q\ell}\end{equation}
where $G_{1h}$ is the hermitian part of $L(P_1)$, and $G_{Q\ell}$ consists of the
remaining
factors, including $G_{m\ell},$ the linear
part of $L(P_m).$  Each factor in $G_{Q\ell}$ acts as an automorphism group of a certain
symmetric cone in the boundary of the cone $C(P_1).$  The projection
\[ \overline{Q} = \nu_{1}(Q) = G_{1h} \mathcal U_{P_1Q} G_{Q\ell} \]
is parabolic in $L(P_1)$ with unipotent radical $\mathcal U_{P_1Q} \subset G_{1\ell}$
which also has a lift depending on the choice of basepoint.  In summary we obtain a
decomposition 
\begin{equation}\label{eqn-decompose2} Q = \mathcal U_1  G_{1h}  \mathcal
U_{P_1Q} G_{Q\ell}.\end{equation}

\section{Cayley transform}
\subsection{}As in \S \ref{sec-hermitian}, suppose that $\mathbf G$ is defined over
$\mathbb Q$ and simple over $\mathbb Q$, that $G = \mathbf G(\mathbb R)$, and that $K$ is
a maximal compact subgroup of $G$ with $D = G/K$ Hermitian.  The following proposition is
the key technical tool behind our construction of a connection which is flat
along the fibers of $\pi.$  The proof follows from the existence of a ``canonical
automorphy factor for $P_1$'' as defined by M. Harris \cite{Harris2} \cite{HZ1} (1.8.7).
See also the survey in \cite{Zucker2}.  In this section we will approximately follow
\cite{Harris2} and derive these results from known facts about the Cayley transform
\cite{WK}, \cite{Satake} Chapter III.

\begin{prop}\label{prop-Cayley}
Let $P_1 = \mathcal U_{1} G_{1h}G_{1\ell}$ be a maximal rational parabolic subgroup of
$\mathbf G.$  Set $K_1 = K\cap P_1= K_{1h}K_{1\ell}.$  Let $\lambda:K \to GL(V)$ be a
representation of $K$.   Then the restriction $\lambda|K_1$ admits a natural extension
$\lambda_1:K_{1h}G_{1\ell} \to \text{GL}(V).$
\end{prop}

\quash{ It will be necessary to briefly review this material since we need some finer
properties of $\lambda$ in Proposition \ref{prop-finer}.} 

\subsection{Proof}\label{subsec-canon-automorphy}
The group $K$ is the set of real points of an algebraic group $\mathbf K$
defined over $\mathbb R.$ As in
\cite{Helgason} VIII\S 7, \cite{Satake} II \S 4, \cite{AMRT} III \S 2 when the Cartan
decomposition $\mathfrak g = \mathfrak k \oplus \mathfrak p$ is complexified it gives
rise to two abelian unipotent subgroups $P^{+}$ and $P^{-}$ of $\mathbf G(\mathbb C)$
such that the complex structure on $\mathfrak g(\mathbb C)$ acts with eigenvalue $\pm i$
on $\text{Lie}(P^{\pm}).$  The natural mapping $P^{+}\mathbf K({\mathbb C}) P^{-}\to
\mathbf G(\mathbb C)$ is injective and its image contains $G=\mathbf G(\mathbb R)$.  Let
$j:P^{+}\mathbf K(\mathbb C)P^{-} \to \mathbf K(\mathbb C)$ denote the projection to the
middle factor.  The group $P^{+}$ is the unipotent radical of the maximal parabolic
subgroup $P^{+}\mathbf K(\mathbb C)$ and hence is normal in $P^{+}\mathbf K(\mathbb C)$;
similarly $P^{-}$ is normal in $\mathbf K(\mathbb C)P^{-}.$  It follows that, for all $h
\in P^{+}\mathbf K(\mathbb C)$, for all $h'\in \mathbf K(\mathbb C)P^{-}$ and for all
$g\in P^{+}\mathbf K(\mathbb C)P^{-}$ we have 
\begin{equation}\label{eqn-hgh} j(hgh') = j(h)j(g)j(h'). \end{equation}
In particular, $j(gk)=j(g)k$ whenever $k\in \mathbf K(\mathbb C)$, so by equation
(\ref{eqn-aut3}), $j$ determines a unique automorphy factor (the ``usual''
canonical automorphy factor) $J_0:G\times D \to \mathbf K(\mathbb C)$ such that
$J_0(g,x_0) = j(g)$ (for all $g\in G$).

Let $c_1 \in \mathbf G(\mathbb C)$ be the \quash{(inverse of) the} Cayley element
\cite{Satake},\cite{WK}, \cite{BB}, \cite{Harris2}, \cite{Zucker2}  which corresponds to
$\mathbf {P_1}$.  That is, $c_1$ is a lift to $\mathbf G(\mathbb C)$ of the standard
choice of Cayley element in $\text{Ad}(\mathfrak g)(\mathbb C)$.  We follow Satake's
convention, rather than that of \cite{WK} (which would associate $c_1^{-1}$ to $\mathbf
P_1.$)  The element $c_1$ satisfies the following
properties, \cite{Satake} (Chapt. III (7.8),(7.9),(2.4)), \cite{Harris2}, \cite{Zucker2}:
\begin{enumerate}
\item  $c_1, c^{-1}_1 \in P^{+} \mathbf K(\mathbb C) P^{-}$ and $c_1G \subset P^+\mathbf
K(\mathbb C) P^- ;$
\item $c_1$ commutes with $G_{1h};$
\item $c_1 \mathcal U_1 K_{1h} G_{1\ell} c_1^{-1} \subset P^{+}\mathbf K(\mathbb C);$
\item $c_1 K_{1h}G_{1\ell}c_1^{-1} \subset \mathbf K(\mathbb C),$ and
hence
\item $j(c_1g) = j(c_1gc_1^{-1})j(c_1) = c_1gc_1^{-1}j(c_1)$ for all $g\in
K_{1h}G_{1\ell}.$ \end{enumerate}

Harris then defines the {\it canonical automorphy factor} (for $\mathbf{P_1}$) to be the
automorphy factor $J_1:G \times D \to \mathbf K(\mathbb C)$ which is determined by its
values at the basepoint, $J_1(g,x_0) = j(c_1)^{-1}j(c_1g).$ This is well defined because
$j(c_1)^{-1}j(c_1gk)= j(c_1)^{-1}j(c_1g)k$ by (\ref{eqn-hgh}), see (\ref{eqn-aut3}). So
we may define the {\it canonical extension}
\begin{equation}\label{eqn-lambda1}
\lambda_1(g) = \lambda_{\mathbb C}J_1(g,x_0) = \lambda_{\mathbb C}(j(c_1)^{-1}
j(c_1g)) \end{equation}
for any $g\in K_{1h}G_{1\ell},$ where $\lambda_{\mathbb C}$ is the complexification of
$\lambda.$  Then $\lambda_1(k) = \lambda(k)$ for any $k \in K_1.$
Moreover $\lambda_1$ is a homomorphism:  if $k_{1h}g_{1\ell}$ and $k'_{1h}g'_{1\ell}$ are
elements of $K_{1h}G_{1\ell}$ then
\begin{align*}
J_1(k_{1h}g_{1\ell} k'_{1h}g'_{1\ell},x_0) &= j(c_1)^{-1}(c_1k_{1h}
g_{1\ell}c_1^{-1}\cdot c_1k'_{1h} g'_{1\ell} c_1^{-1}) j(c) \\
&= j(c_1)^{-1} c_1k_{1h}g_{1\ell}c_1^{-1}j(c_1)\cdot j(c_1)^{-1}c_1 k'_{1h} g'_{1\ell}
c_1^{-1} j(c_1) \\
&= J_1(k_{1h}g_{1\ell},x_0) J_1(k'_{1h} g'_{1\ell},x_0).  \end{align*}
Verification that $J_1(k_{1h}^{-1}g_{1\ell}^{-1},x_0) = J_1(k_{1h}g_{1\ell},x_0)^{-1}$ is
similar.  We remark, following \cite{Harris2} that modifying $c_1$ by any element $d\in
P^{+}\mathbf K(\mathbb C)$ will not affect the values of $J_1(k_{1h}g_{1\ell},x_0).$ \qed

Now suppose $\mathbf{P_1} \prec \mathbf{P_2}$ are rational maximal parabolic subgroups of
$\mathbf G$ with $P_i = \mathcal U_i G_{ih}G_{i\ell}$ ($i = 1,2$) and, as in
(\ref{eqn-twopb}),  
\[ P_1 \cap P_2 = \mathcal U_2 (\overline{\overline{\mathcal U}} G_{1h} G'_{\ell})
G_{2\ell} \text{ with } P_h = \overline{\overline{\mathcal U}} G_{1h} G'_{\ell} \subset
G_{2h}. \]

\begin{prop}\label{prop-finer}
Let $\lambda_i:K_{ih}G_{i\ell} \to \text{GL}(V)$ be the canonical extensions of
$\lambda|K_i= K \cap P_i$ ($i = 1,2$) and let
$\lambda_{21}:K_{1h}G'_{\ell} \to \text{GL}(V)$ be the canonical extension of
$\lambda|K_{2h}$ corresponding to the canonical automorphy factor $J_{21}$ for $P_h
\subset G_{2h}.$  Then
\[ \lambda_1(g_{2\ell}) = \lambda_2(g_{2\ell}) \text{ and }  
 \lambda_1(g'_{\ell}) = \lambda_{21}(g'_{\ell})\] 
 for all $g_{2\ell} \in G_{2\ell}$ and all $g'_{\ell} \in G'_{\ell}.$
\end{prop}
\subsection{Proof}
Let $c_1,c_2\in \mathbf G(\mathbb C)$ and $c_{21} \in \mathbf{G_{2h}}(\mathbb C)$ be the
Cayley elements for $ P_1,P_2 \subset G$ and $P_h \subset G_{2h}$ respectively.  By
\cite{Satake} Chapter III (9.5),
\[c_{21} = c_1c_2^{-1} = c_2^{-1}c_1. \]
The lift $G_{2h}\subset G$ is stable under the Cartan involution on $G$ so the
corresponding decomposition $P_{2h}^+ \mathbf{K_{2h}}(\mathbb C) P_{2h}^-$ coincides
with $(G_{2h} \cap P^+)(G_{2h} \cap \mathbf K(\mathbb C)) (G_{2h} \cap P^-)$ and in
particular $j_{2h}:G_{2h} \to \mathbf{K_{2h}}(\mathbb C)$ is the restriction of $j:G \to
\mathbf K(\mathbb C).$ Let $c_{21} = c_{21}^{+}c_{21}^0c_{21}^{-}$ be the resulting
decomposition. Then $c_2$ commutes with each of the
factors $c_{21}^*.$ It follows from (\ref{eqn-hgh}) that
\[ j(c_1) = j(c_{21}c_2)= j(c_{21}^+c_{21}^0c_2 c_{21}^-) = j(c_{21}^0 c_2)
= j(c_{21})j(c_2)\]
and similarly $j(c_1) = j(c_2)j(c_{21}).$  Since $g_{2\ell}\in G_{2\ell}$ also commutes
with $c_{21}\in G_{2h}(\mathbb C)$ and with $j(c_{21})\in K_{2h}(\mathbb C)$, we find
\begin{align*}
\lambda_1(g_{2\ell})\lambda_2(g_{2\ell})^{-1} &= \lambda_{\mathbb C} \left(
j(c_1)^{-1} c_1 g_{2\ell} c_1^{-1} j(c_1) j(c_2)^{-1} c_2g_{2\ell}^{-1} c_2^{-1} j(c_2)
\right)
\\
&= \lambda_{\mathbb C} \left( j(c_1)^{-1}c_1g_{2\ell}c_1^{-1} j(c_{21})c_2 g_{2\ell}^{-1}
c_2^{-1} j(c_2)\right) \\
&= \lambda_{\mathbb C} \left( j(c_1)^{-1} c_1 g_{2\ell} c_1^{-1}c_2 g_{2\ell}^{-1}
c_2^{-1} j(c_{21})j(c_2)\right) \\
&= \lambda_{\mathbb C} \left( j(c_1)^{-1}c_1g_{2\ell} c_{21}^{-1}g_{2\ell}^{-1} c_2^{-1}
j(c_1)\right) \\
&= \lambda_{\mathbb C} \left( j(c_1)^{-1}c_1c_{21}^{-1}c_2^{-1}j(c_1)\right) = 1.
\end{align*}
Similarly if $g'_{\ell} \in G'_{\ell}$ then using (4) above, 
\begin{align*}
\lambda_{21}(g'_{\ell})  &= \lambda_{\mathbb C} \left(
j(c_{21})^{-1}j(c_{21}g'_{\ell}c_{21}^{-1})j(c_{21})\right) \\
&= \lambda_{\mathbb C} \left(j(c_{21})^{-1} j(c_2^{-1}c_1 g'_{\ell} c_1^{-1}c_2)j(c_{21})
\right) \\
&= \lambda_{\mathbb C} \left( j(c_{21})^{-1} j(c_2)^{-1}
j(c_1g'_{\ell}c_1^{-1})j(c_2)j(c_{21})\right) \\
&= \lambda_{\mathbb C} \left( j(c_1)^{-1} j(c_1 g'_{\ell}c_1^{-1})j(c_1)\right) =
\lambda_1(g'_{\ell}). \qed \end{align*}


\section{    Baily-Borel Satake compactification  }\label{sec-BB}
\subsection{  }\label{subsec-BB1}
As in \S \ref{sec-hermitian}, suppose that $\mathbf G$ is defined over $\mathbb Q$ and
simple over $\mathbb Q$, that $G = \mathbf G(\mathbb R)$, and that $K$ is a maximal
compact subgroup of $G$ with $D = G/K$ Hermitian.  Let $D^*$ be the Satake partial
compactification of $D$, consisting of $D$ together with its rational boundary
components $D_P$, one for each (proper) maximal rational parabolic subgroup $\mathbf P
\subset \mathbf G$; with the Satake topology \cite{BB}.  The action of $\mathbf G(\mathbb
Q)$ on $D$ extends continuously to an action of $\mathbf G(\mathbb Q)$ on $D^*.$  Let
$\Gamma \subset \mathbf G(\mathbb Q)$ be a neat arithmetic subgroup and let $q:D^*\to
\overline{X} = \Gamma \backslash D^*$ denote the quotient mapping.  Then $\overline{X}$
is the  Baily-Borel compactification  of $X$ and it admits the structure of a complex
projective algebraic variety with a canonical stratification with a single stratum for
every $\Gamma$-conjugacy class of rational boundary
components as follows.   Let $D_1\subset D^{*}$ be a rational boundary component with
normalizing maximal parabolic subgroup $P= \mathcal{ U}_PG_hG_{\ell}.$  Let $\nu_h:P\to
G_h' = G_h/(G_h\cap G_{\ell})$ and let $\nu_{\ell}: P\to G_{\ell}' = G_{\ell}/(G_h\cap
G_{\ell}).$  The closure
$D_1^{*}$ of $D_1$ in $D^{*}$ is the Satake partial compactification of $D_1$.
The group $P$ acts on $D_1$ through its projection to $G_h'$ and the group
$\Gamma_P = \Gamma \cap P$ acts on $D_1^{*}$ through its projection $\Gamma_h
= \nu_h(\Gamma_P)$ to $G_h'.$  Then $X_1 = \Gamma_P \backslash D_1^{*} =
\Gamma_h \backslash D_1^{*}$ is a stratum of $\bar X.$  Its closure
$\bar X_1 = \Gamma_h\backslash D_1^{*}$ in $\bar X$ is the Baily-Borel
compactification   of $X_1.$  The stratum $X_1$ is also the image of the
(infinitely many) rational boundary components $D'_1$ which are
$\Gamma$-conjugate to $D_1.$

\subsection{  } 
Let $D_1\subset D^*$ be a rational boundary component which projects to $X_1$.  We will
say that a neighborhood $\tilde U \subset D^*$ is {\it a $\Gamma$-parabolic neighborhood}
of $D_1$ if the following holds:  if $x_1,x_2\in \tilde U$ and $\gamma \in \Gamma$
satisfy $x_2 = \gamma x_1$ then $\gamma \in \Gamma \cap P.$  If $X_1 \subset \bar X$ is a
stratum in the Baily-Borel compactification of $X$, we say that a neighborhood $U\subset
\bar X$ of $X_1$ is parabolic if for some (and hence for any) boundary component $D_1
\subset D^*$ with $q(D_1) = X_1$, there is a $\Gamma$-parabolic neighborhood $\tilde
U\subset D^*$ of $D_1$ such that $U = q(\tilde U).$  This means that the covering
$\Gamma_P \backslash D^* \to \Gamma \backslash D^*$ is one to one on $\tilde U$, and we
have a commutative diagram
\begin{equation*}
\begin{CD}
D^* &\ \supset\ &\ \tilde U \ &\ \supset\ &\ D_1 \\
@V{q}V{\Gamma}V @VV{\Gamma_P}V @VV{\Gamma_h}V \\
\overline{X}\ &\ \supset\ &\ U\ &\ \supset\ &\ X_1
\end{CD}.
\end{equation*}
\begin{lem}Each stratum $X_1\subset\bar X$ has a
fundamental system of neighborhoods, each of which is $\Gamma$-parabolic.
\end{lem}
\subsection{Proof}\cite{Saper2} Let $D^{BS}$ be the Borel-Serre partial compactification
of $D$ together with its ``Satake'' topology \cite{BorelSerre}.  It is a manifold with
corners, having one corner $e(P)$ for each rational parabolic subgroup $\mathbf P$.
According to \cite{ZuckerSatake} the identity mapping $D\to D$ has a unique continuous
extension  $\mathfrak w: D^{BS} \to D^*,$ and it is surjective.  If $\mathbf P$ is
standard then $\mathfrak w(e(P)) = D_{P^{\flat}}.$  Let $D^{\dagger}$ denote the quotient
topology on the underlying set $|D^*|$ which is induced by $\mathfrak w.$  Then
$D^{\dagger} \to D^*$ is a continuous bijection and the quotient mapping $\Gamma
\backslash D^{\dagger} \to \Gamma \backslash D^*$ is a homeomorphism.  In \cite{Saper}
Theorem 8.1, Saper constructs a basis of parabolic neighborhoods $U_P$ of each corner
$e(P)$ in $D^{BS}$.  If $\mathbf P$ is maximal then $\mathfrak w(U_P)$ is open in $D^*$
as may be shown by verifying the condition at the bottom of page 264 in \cite{AMRT}. \qed 

We remark that the image of $U_P$ is $\Gamma$-parabolic and is open in $D^{\dagger}$ by
construction, and that the topology $D^{\dagger}$ may be substituted for the Satake
topology $D^*$ throughout this paper.

\subsection{  }\label{subsec-compatible}
Fix a standard rational boundary component $D_1\subset D^*$ normalized by a standard
maximal rational parabolic subgroup $P_1.$  In the Satake topology (or in the topology
$D^{\dagger}$) there is a natural
neighborhood $T(D_1) = \bigcup \{D_2|\ D_1 \prec D_2 \prec D \}$ consisting of
the union of all rational boundary components (including $D$, the nonproper boundary
component) whose closures contain $D_1$.  The projection $\pi:D \to D_1$ (\ref{eqn-pi})
has a unique continuous extension $T(D_1) \to D_1$ to this neighborhood.  Its restriction
to each intermediate boundary component $D_2$  coincides with the canonical
projection $D_2 \to D_1$ which is obtained by considering $D_2$ to be the symmetric space
corresponding to the Hermitian part $G_{2h}$ of the Levi factor of $P_2$ and by
considering $D_1\subset D_2^{*}$ to be the rational boundary component preserved by the
parabolic subgroup $P_h\subset G_{2h}$ (notation as in \S \ref{subsec-twoparab}).  
It follows that
$\pi_1 (x) = \pi_1 \pi_2 (x)$ 
for all $x\in T(D_1)\cap T(D_2).$  (The above union can be quite large:  if $D_{20}$ is a
standard boundary component normalized by a standard parabolic subgroup $\mathbf{P_2}
\succ \mathbf{P_1}$ then the boundary components $D_2 \subset T(D_1)$ which are conjugate
to the standard one $D_{20}$ are in one to one correspondence with elements of
$\mathbf{P_1}(\mathbb Q)/(\mathbf{P_1}(\mathbb Q) \cap \mathbf{P_2}(\mathbb Q))$ cf.
(\ref{eqn-twopa}).)

\begin{lem}\label{lem-BBdata}For $\epsilon_0 > 0$ sufficiently small, 
the Baily-Borel compactification $\overline X$ admits an $\epsilon_0$-system of control
data $\{T_Y(\epsilon_0),$ $\pi_Y,\rho_Y\}$ such that for each stratum $Y \subset \partial
\overline{X}$ and for any choice of boundary component $D_1 \subset D^*$ with $q(D_1)=Y$
we have \begin{enumerate}
\item  The neighborhood $T_Y(\epsilon_0)\subset \overline{X}$ is a parabolic neighborhood
of $Y$; it is the image, say, of some $\Gamma$-parabolic neighborhood $\tilde
U(\epsilon_0)\subset D^*$ of $D_1$, and
\item $\pi_Y(q(x)) = q(\pi(x))$ for all $x\in \tilde U(\epsilon_0)$ (where $\pi:D \to
D_1$ is the canonical projection).
\end{enumerate}
\end{lem}

\subsection{Proof}
For any $g\in P$ and $x\in D$ we have $\pi(gx) = \nu_h(g)\pi(x) \in D_1$.  It follows
that the projection function $\pi$ passes to the quotient $\Gamma_P \backslash D^*$,
where it may be restricted to a parabolic neighborhood $U$ of $Y = q(D_1)$; write
$\pi_Y:U\to Y$ for the result.   If $P$ and $P'$ are $\Gamma$ conjugate maximal rational
parabolic subgroups corresponding to conjugate boundary components $D_1$ and $D'_1$ then
the projections $T(D_1)\to D_1$ and $T(D'_1) \to D'_1$ are compatible with conjugation,
which shows that the resulting projection $\pi_Y:U\to Y$ is independent of the choice of
lift $D_1\subset D^*$ of the stratum $Y\subset \bar X$.   The tubular neighborhood
$T_Y(\epsilon_0)$ may be chosen inside $U$.  The compatibility between these projections
follows from (\ref{subsec-compatible}).  As mentioned in \S \ref{subsec-existence}, by
further shrinking the tubular neighborhoods if necessary, control data may be found for
which the tubular projections agree with these $\pi_Y.$ \qed


\section{ Parabolically Induced Connection  }\label{sec-induced}
\subsection{   }
As in \S \ref{sec-hermitian}, \ref{sec-BB} we suppose that $\mathbf G$ is semisimple,
defined over $\mathbb Q$ and simple over $\mathbb Q$; that $G = \mathbf G(\mathbb R)$,
and $K\subset G$ is a maximal compact subgroup with $D = G/K$ Hermitian symmetric.  Fix
$\Gamma \subset \mathbf G(\mathbb Q)$ a neat arithmetic subgroup.  Let $\lambda:K \to
GL(V)$ be a representation of $K$ on some complex vectorspace $V$ and denote by $E =
G\times_K V$ the associated homogeneous vectorbundle on $D$.

Let $D_1$ be a rational boundary component of $D$ with canonical projection
$\pi:D\to D_1.$  Let $P$ be the maximal parabolic subgroup of
$G$ which preserves $D_1.$  Write $P=\mathcal{U}G_hG_{\ell}$ as in \S
\ref{subsec-parabolics} and let $K_h = K\cap G_h$ and $K_{\ell} = K \cap
G_{\ell}$ be the corresponding maximal compact subgroups.  Let $\mathfrak g_h =
\mathfrak k_h \oplus \mathfrak p_h$ and $\mathfrak g_{\ell} = \mathfrak k_{\ell} \oplus
\mathfrak p_{\ell}$ denote the corresponding Cartan decompositions.

The restriction of $\lambda$ to $K_h$ determines a homogeneous vectorbundle $E_1 = G_h
\times_{K_h}V$ over $D_1$.  By Proposition \ref{prop-Cayley} the representation
$\lambda|\ K_hK_{\ell}$ admits an extension to a representation $\lambda_1:K_hG_{\ell}
\to GL(V).$  This extension determines an action of $P$ on $E_1$ which is given by
\begin{equation}\label{eqn-Paction}
ug_hg_{\ell}.[g_h',v] = [g_hg_h',\lambda_1(g_{\ell})v]. \end{equation}
We obtain a vectorbundle mapping (which covers $\pi$), 
\begin{equation}\label{eqn-Phi}
\tilde\Phi:E=P\times_{K_P}V \to G_h\times_{K_h}V = E_1\end{equation}
by $\tilde\Phi([ug_hg_{\ell},v]) = [g_h,\lambda_1(g_{\ell})v].$  Then $\tilde\Phi$
induces an isomorphism, 
\begin{equation}\label{eqn-Phiiso}  \Phi:E\cong \pi^{*}(E_1);\quad [g,v] \mapsto 
(gK_P,\tilde\Phi([g,v])) \in D\times E_1\end{equation}
 of $P$-homogeneous vectorbundles (where $g\in P$ and $v \in K$).

\subsection{  Definition }\label{def-parab}  Let $\nabla_1 = d+\omega_1$ be a
connection on $E_1.$  The {\it parabolically induced connection} $\nabla = d+\omega$ on
$E$ is defined to be the pullback $\nabla = \Phi^{*}(\nabla_1)$ of $\nabla^1$ under the
isomorphism $\Phi.$  It is the unique connection whose covariant derivative
(\ref{eqn-covariant}) satisfies
\begin{equation}\label{eqn-pullback}
\nabla_v(\Phi^{*}(s)) = \Phi^{*}((\nabla_1)_{\pi_*v}s)
\end{equation}
for any section $s$ of $E_1$ and for any tangent vector $v\in T_xD.$  
\quash{  Parabolic induction
is an affine operation.  Suppose $\nabla_1 = \sum_{i=1}^r B_i\nabla'_i$ is a convex
combination of connections on $E_1$, where $\sum_{i=1}^rB_i=1$ is a partition of unity on
$D_1.$  Then
\[ \Phi^*(\nabla_1) = \sum_{i=1}^r \pi^*(B_i) \Phi^*(\nabla'_i).  \]
} 
\begin{prop}\label{prop-induced}  Suppose $\nabla_1=d+\omega_1$ is a
connection on $E_1=G_h\times_{K_h}V$.  Let $\nabla =d+\omega$ denote the
parabolically induced connection on $E=P\times_{K_P}V$.  Then
\begin{equation}\label{eqn-induced}
\omega (L_{g*}(\dot {u}+\dot {g}_h+\dot {g}_{\ell}))=\lambda'_1(\dot{g}_{\ell})
+ Ad(\lambda_1(g_{\ell}^{-1}))(\omega_1(L_{g_h*}(\dot {g}_h))\end{equation}
for any $g=ug_hg_{\ell}\in P$  and any $\dot {u}+\dot {g}_h+\dot
{g}_{\ell}\in \text{Lie}(\mathcal U_P)\oplus \mathfrak g_h\oplus \mathfrak g_{\ell}$.
\end{prop}

\subsection{  Proof  }\label{subsec-proof-induced}
  Let $J_1:G_h\times D_1\to GL(V)$ be an automorphy factor for $E_1$,
corresponding to a trivialization $E_1 \cong D_1 \times V$.  Composing this
with the isomorphism $\Phi:E \to \pi^{*}(E_1)$ determines an automorphy
factor $J:P\times D \to GL(V)$ with
\begin{equation}\label{eqn-autP}
J(ug_hg_{\ell},x_0) = J_1(g_h,x_1)\lambda_1(g_{\ell})\end{equation}
where $x_1 = \pi(x_0)\in D_1$ denotes the basepoint in $D_1.$  To simplify
notation we will write $j(g)$ and $j_1(g_h)$ rather than $J(g,x_0)$ and
$J_1(g_h,x_1)$.

By (\ref{eqn-related}), the connection $\nabla_1$ in $E_1$ determines a connection
$\nabla^{J_1} = d+\eta_1$ in the $J_1$-trivialization $E_1 \cong D_1\times V$ with 
\begin{equation}  \eta_1(q_{1*}(X_h))=j_1(g_h)\omega_1(X_h)j_1(g_h)^{-1}-d_{X_h}
(j_1(g))\circ j_1(g)^{-1}\end{equation}
for any $g_h\in G_h$ and any $X_h\in T_{g_h}G_h$, where $q_1:G_
h\to D_1$ denotes the projection.

The parabolically induced connection $\nabla$ in $E$ determines a connection $\nabla^J=
d+\eta$ in the $J$-trivialization $E\cong D\times V$ with 
\begin{equation}  \eta (q_{*}X)=j(g)\omega (X)j(g)^{-1}-d_X(j(g))\circ
j(g)^{-1}.\end{equation}
By (\ref{eqn-pullback}) and (\ref{eqn-covariant}) the connection forms $\eta$ and
$\eta_1$ are related by $\eta (q_{*}(X))=\eta_1(\pi_{*}q_{*}(X))$  for any $X\in T_gG.$
Take $X = L_{g*}(\dot u + \dot g_h + \dot g_{\ell}) \in T_gG$ and let $X_h =
L_{g_h*}(\dot g_h) \in T_{g_h}G_h$ denote its projection to $G_h.$  Then we
have
\begin{align*}
j_1(g_h)\omega_1(X_h)j_1(g_h)^{-1}=&j_1(g_h)\lambda_1 (g_{\ell})\omega
(X)\lambda_1 (g_{\ell})^{-1}j_1(g_h)^{-1}\\
 &-j_1(g_h)d_X(\lambda_1 (g_{\ell}))\lambda_1 (g_{\ell})^{-1}j_1(g_h)^{-1}
\end{align*}   or, using (\ref{eqn-Lg})
\begin{equation*}  \omega_1(L_{g_h*}(\dot g_h)) =
\lambda_1(g_{\ell})\omega(X)\lambda_1(g_{\ell}^{-1}) - \lambda_1(g_{\ell})
\lambda_1'(\dot g_{\ell}) \lambda_1(g_{\ell}^{-1})\qed\end{equation*}

\quash{
\subsection{Remark}  If, in (\ref{eqn-autP}), we take $J_1 = \lambda \circ J_h :G_h
\times D_1 \to GL(V)$ to be the (usual) canonical automorphy factor for $G_h$, then the
resulting automorphy factor $J$ of (\ref{eqn-autP}) is closely
related to but does not necessarily coincide with Harris' ``canonical automorphy
factor for $P$'' (Proposition \ref{prop-Cayley}) which we will denote by
$J^H.$  Both automorphy factors agree with $J_h$ on the boundary component.
Our $J(g,x_0)$ is trivial for all $g\in \mathcal U_P$ while $J^H(g,x_0)$ is
trivial for $g\in \text{Center}(\mathcal U_P).$  The automorphy factor $J^H$ is
defined on $G\times D \to GL(V)$ whereas our $J:P\times D\to GL(V)$ does not
necessarily extend to $G\times D.$  Nevertheless, the $J$ and $J^H$
automorphic actions of $P$ on the trivial vectorbundle $P/K_P \times V$ are
equivalent    (cf. \S \ref{subsec-twoauto}).}  


\begin{cor}\label{cor-parab}  
Suppose $\omega_1 \in \mathcal A^1(G, End(V))$ commutes with the adjoint action
of $\lambda_1:G_{\ell} \to GL(V).$  Then the curvature form $\Omega$ of the parabolically
induced connection $\nabla = \Phi^*(\nabla_1) = d+\omega$ satisfies the $\pi$-fiber
condition,
\[ \Omega = \pi^*(\Omega_1) \]
where $\Omega_1\in \mathcal A^2_{\text{bas}}(G_h,\text{End}(V))$ is the curvature form of
$\nabla_1 = d+\omega_1.$ \end{cor}

\subsection{Proof}
Let us compute $\Omega(X,Y)$ where $X = L_{g*}(\dot x)$ and $Y = L_{g*}(\dot y),$ where
$\dot x, \dot y \in \text{Lie}(P)$ and where $g=ug_hg_{\ell}\in P.$
Set $\dot x = \dot u_X + \dot g_{Xh} + \dot g_{X\ell}\in \text{Lie}(\mathcal U_P) \oplus
\mathfrak g_h \oplus \mathfrak g_{\ell}$ (and similarly for $\dot y$).  By
(\ref{eqn-induced}), 
\begin{equation}\label{eqn-cor-parab} \omega(X) =  \omega_1(L_{g_h*}(\dot g_{Xh}))  +
\lambda_1'(\dot g_{X\ell}) \end{equation}
and similarly for $\omega(Y).$  Set $X_h = L_{g_h*}(\dot g_{Xh})$ and $Y_h =
L_{g_h*}(\dot g_{Yh}). $  Using the structure equation (\ref{eqn-structure}) and
the fact that $\text{Lie}(\mathcal U_P)$ is an ideal in $\text{Lie}(P)$ gives
\begin{align*} 
\Omega(X,Y) &= X(\omega(Y)) - Y(\omega(X)) -\omega([X,Y]) + [\omega(X),\omega(Y)] \\ 
       &=  X(\omega_1(L_{g_h*}(\dot g_{Yh}))) + X(\lambda'_1(\dot g_{Y\ell}))
               -Y(\omega_1(L_{g_h*}(\dot g_{Xh}))) - Y(\lambda'_1(\dot g_{X\ell})) \\
            &\quad - \omega_1(L_{g_h*}([\dot g_{Xh}, \dot g_{Yh}])) 
             - \lambda'_1([\dot g_{X\ell}, \dot g_{Y\ell}]) \\
            &\quad +[\omega_1(L_{g_h*}(\dot g_{Xh})), \omega_1(L_{g_h*}(\dot g_{Yh}))] 
             + [ \lambda'_1(\dot g_{X\ell}), \lambda'_1(\dot g_{Y\ell})] \\
   &= X_h \omega_1(L_{g_h*}(\dot g_{Yh}))) - Y_h \omega_1 (L_{g_h*}(\dot g_{Xh}))) - 
    \omega_1([X_h,Y_h]) + [\omega_1(X_h), \omega_1(Y_h)]\\
   & = \Omega_1(X_h,Y_h). \qed
\end{align*}

\subsection{} \label{subsec-nuh}
If $\nabla_1 = d+\omega_1$ is a
connection on $E_1$ which is invariant under a subgroup $\Gamma_h \subset G_h$, then 
by (\ref{eqn-induced}) the induced connection $\Phi^*(\nabla_1)$ is invariant under the
group $\nu_h^*(\Gamma_h) \subset P$ which is obtained by first projecting $\Gamma_h$ to
$G'_h = G_h/(G_h\cap G_{\ell})$ then taking the pre-image under the projection $\nu_h: P
\to G'_h$ (cf. \S \ref{subsec-parabolics}, \S \ref{subsec-BB1}).

As in \S \ref{sec-BB} let $\Gamma \subset G$ be a neat arithmetic subgroup with $X =
\Gamma \backslash D.$  Write $X_1 = \Gamma_h \backslash D_1$ for the stratum in
$\overline{X} = \Gamma \backslash D^*$ corresponding to the boundary component $D_1.$
The homogeneous vectorbundles $E \to D$ and $E_1 \to D_1$ pass to automorphic
vectorbundles $E' \to X$ and $E'_1 \to X_1$ respectively.  Parabolic induction then
passes to an operation on these vectorbundles as follows.
Suppose $\nabla'_1$ is a connection on $E'_1.$  It pulls back to a $\Gamma_h$-invariant
connection $\nabla_1$ on $E_1 \to D_1.$  The parabolically induced connection $\nabla =
\Phi^*(\nabla_1)$
is invariant under $\Gamma_P =\Gamma \cap P \subset \nu_h^*(\Gamma_h)$ so it passes to a
connection on $\Gamma_P\backslash E \to \Gamma_P \backslash X.$  Since
$T_{X_1}(\epsilon_0)$
is a $\Gamma$-parabolic neighborhood of $X_1$ in $\overline{X}$ this defines a connection 
$\nabla' = \Phi^*_{XX_1}(\nabla'_1)$ on the restriction $E'|(X \cap
T_{X_1}(\epsilon_0)).$

\quash{Since $\overline{X}_1 = \Gamma_h \backslash D^*_1$ is the Baily-Borel
compactification of
$X_1$,} This procedure may be applied to any pair of strata, say, $X_2 < X_1$ of
$\overline{X}.$  Thus, if $\nabla'_2$ is any connection on the automorphic vectorbundle
$E'_2\to X_2$ defined by $\lambda$, then we obtain a parabolically induced connection
\[ \nabla'_1 = \Phi^*_{X_1X_2}(\nabla'_2) \]
on $E'_1 | (X_1 \cap T_{X_2}(\epsilon_0)).$  However, if $X_3 < X_2 < X_1$ are strata of
$\overline{X}$ and if $\nabla'_3$ is a connection on $E'_3 \to X_3$ then the
parabolically induced connection $\Phi^*_{X_1X_3}(\nabla'_3)$ does not necessarily agree
with the connection $\Phi^*_{X_1X_2}\Phi^*_{X_2X_3}(\nabla'_3),$ even in the neighborhood
$X_1 \cap T_{X_2}(\epsilon_0) \cap T_{X_3}(\epsilon_0)$ where they are both defined, cf.
Proposition \ref{prop-multipleinduction}.



\subsection{}\label{subsec-decompositions}  
We will also need the following more technical result concerning parabolic
induction for the proof of the main theorem.  Suppose $P_1  \prec P_2 \prec \ldots \prec
P_r$ are standard maximal parabolic subgroups with corresponding rational boundary
components $D_1 \prec {D}_2 \prec \cdots \prec {D}_r.$  Let $Q = P_1 \cap P_2\cap \ldots
\cap P_r.$  Then $Q^{\flat} = P_1 = \mathcal U_1 G_{1h}G_{1\ell}$ so by
(\ref{eqn-decompose2}), $Q$ decomposes as $Q = \mathcal U_QG_{1h}G_{Q\ell}.$

Let $\lambda$ be a representation of $K$ with resulting homogeneous vectorbundles $E_i
\to D_i$ ($1 \le i \le r$).  By Proposition \ref{prop-Cayley}, $\lambda$ extends to a
representation $\lambda_1$ of $K_{1h}G_{1\ell}.$  For $1 \le i < j \le r$ let
$\Phi_{ji}:E_j \to \pi_{ji}^*E_i$ denote the vectorbundle isomorphism of (\ref{eqn-Phi})
which covers the canonical projection $\pi_{ji}:D_j \to D_i.$  Let $\Phi_{r}:E \to
\pi_r^*(E_r)$ be the vectorbundle isomorphism (\ref{eqn-Phi}) which covers the canonical
projection $\pi_r:D \to D_r.$

\begin{prop}\label{prop-multipleinduction}  
Suppose $\nabla_1 = d + \omega_1$ is a connection on $E_1$ and suppose that the
connection form $\omega_1 \in \mathcal A^1(G_{1h}, \text{End}(V))$ commutes with the
adjoint action of $\lambda_1(G_{1\ell}).$  Let $\nabla = d+\omega$ denote the connection
\[ \nabla = \Phi_{r}^* \Phi_{r,r-1}^*\ldots \Phi_{21}^*(\nabla_1).\]
Then 
\begin{equation}\label{eqn-heart}
 \omega(L_{g*}(\dot u_Q + \dot g_{1h} + \dot g_{Q\ell})) 
= \omega_1(L_{g_{1h}*}\dot g_{1h})) + \lambda'_1(\dot g_{Q\ell}) \end{equation}
for any $g=u_Qg_{1h}g_{Q\ell} \in Q$ and any $\dot u_Q + \dot g_{1h} + \dot
g_{Q\ell} \in \text{Lie}(Q) = \text{Lie}(\mathcal U_Q) + \mathfrak g_{1h} + \mathfrak
g_{Q\ell}.$
\end{prop}

\subsection{Proof}
First we determine the connection form of the connection
\[ \nabla_r = \Phi_{r,r-1}^* \Phi_{r-1,r-2}^* \ldots \Phi_{21}^*(\nabla_1) =
d+\omega_r\]
on the vectorbundle $E_r \to D_r.$  Set $P_r = \mathcal U_r G_{rh}G_{r\ell}$.   The
images under the projection $\nu_{rh}:P \to G_{rh}$ of  $Q \subset P_1 \cap P_r$ are
parabolic subgroups $\overline{Q} \subset \overline{P}_1 \subset G_{rh}.$  In fact,
$\overline{Q}$ is the parabolic subgroup corresponding to the flag of rational boundary
components $D_1 \prec D_2 \prec \cdots \prec D_{r-1}$ of $D_r.$ By (\ref{eqn-twopb})
there are compatible decompositions
\begin{alignat*}{2}
P_r &= \mathcal U_r G_{rh} G_{r\ell} & \\
P_1 \cap P_r &= \mathcal U_r (\overline{\overline{\mathcal U}} G_{1h} G'_{\ell})
G_{r\ell} &\text{ with } \overline{P}_1 &= \overline{\overline{\mathcal U}} G_{1h}
G'_{\ell} \\
Q &= \mathcal U_r \overline{\overline{\mathcal U}} G_{1h} (\mathcal U_{\overline{Q}\ell}
G_{\overline{Q}\ell}) G_{r\ell} &\text{ with } \overline{Q} &=
\overline{\overline{\mathcal U}} \mathcal U_{\overline{Q}\ell}G_{1h}G_{\overline{Q}\ell}.
\end{alignat*}
Corresponding to the maximal parabolic subgroup $\overline{P}_1\subset G_{rh}$ the
representation $\lambda|K_{1h}K'_{\ell}$ has a canonical extension
\[ \lambda_{r1}:K_{1h}G'_{\ell} \to \text{GL}(V).\]
According to Proposition \ref{prop-finer}, $\lambda_{r1}|G'_{\ell} = \lambda_1|G'_{\ell}$
which, by assumption, commutes with $\omega_1.$  So by induction,
for any $\bar q = u_{\overline{Q}}g_{1h}g_{\overline{Q}\ell} \in
\overline{Q}$ and for any $\dot{\overline{q}} = \dot u_{\overline{Q}} + \dot g_{1h} +
\dot g_{\overline{Q}\ell} \in \text{Lie}(\overline{Q}),$ 
\[ \omega_r(L_{\bar q*}(\dot u_{\overline{Q}} + \dot g_{1h} + \dot g_{\overline{Q}\ell}))
=\omega_1(L_{g_{1h}*}(\dot g_{1h})) + \lambda'_1(\dot g_{\overline{Q}\ell}).\]
Moreover $\nabla = \Phi_r^*(\nabla_r) = d+\omega$ so by Proposition \ref{prop-induced},
for any $g = u_rg_{rh} g_{r\ell} \in P_r,$
\[ \omega(L_{g*}(\dot u_r + \dot{g}_{rh} + \dot g_{r\ell}))
= \omega_r(L_{\bar q*}(\dot g_{rh})) + \lambda'_1(\dot g_{r\ell})\]
for any $\dot u_r + \dot g_{rh} + \dot q_{r\ell} \in \text{Lie}(P_r).$  Taking $g_{rh} =
\bar q =u_{\overline{Q}}g_{1h}g_{\overline{Q}\ell}\in \overline{Q}$ and 
\[\dot g_{rh} = \dot{\bar q} = \dot u_{\overline{Q}} + \dot g_{1h} +
\dot g_{\overline{Q}\ell}\in \text{Lie}(\overline{Q})\] 
gives equation (\ref{eqn-heart}):
\[ \omega(L_{g*}(\dot u_r + \dot u_{\overline{Q}} + \dot g_{1h} + \dot
g_{\overline{Q}\ell} + \dot g_{r\ell}))
= \omega_1(L_{g_{1h}*}(\dot g_{1h})) + \lambda'_1(\dot g_{\overline{Q}\ell} + \dot
g_{r\ell}). \qed \]


\section{   The Patched Connection  }\label{sec-maintheorem}

\subsection{   }\label{subsec-main1}
As \S \ref{sec-hermitian}, \ref{sec-BB}, suppose that $D=G/K$ is a Hermitian symmetric
space which is irreducible over $\mathbb Q$, and that $\Gamma\subset \mathbf G(\mathbb
Q)$ is a neat arithmetic group.  Let $\overline{X} = \Gamma \backslash D^*$ denote the
Baily Borel compactification of $X=\Gamma\backslash D$ with projection $q:D^* \to
\overline{X}.$  By lemma \ref{lem-BBdata}, for any sufficiently small $\epsilon_0 > 0$
there exists an $\epsilon_0$-system of control data $\{T_Y(\epsilon_0),\pi_Y,\rho_Y\}$
(which we now fix) on $\bar {X},$ so that $\pi_Y$ is obtained from the canonical
projection $D \to D_1$ whenever $q(D_1) = Y$ and so that $T_Y(\epsilon)$ is a
$\Gamma$-parabolic neighborhood of $Y$ in $\overline{X}.$  Applying \S
\ref{subsec-family} to this system of control data yields a partition of unity on each
stratum $Y$ of $\overline{X}$,
\begin{equation}  B_Y^{\epsilon Y}(y)+\sum_{Z<Y}B_Z^{\epsilon Y}(y)=1\end{equation} for
all $y\in Y,$ where $\epsilon Y = \epsilon_0/2^{\dim Y}.$

A choice of representation $\lambda:K \to GL(V)$ on some complex vectorspace $V$
determines homogeneous vectorbundles $E = G\times_KV$ on $D$ and $E_1 = G_h
\times_{K_h}V$ on $D_1$ which pass to automorphic vectorbundles $E'=\Gamma\backslash E$
on $X$ and  $E'_Y \to Y$ on $Y$.  Here, $D_1$ is a rational boundary component (with
$q(D_1)=Y$), normalized by some maximal parabolic subgroup $P =
\mathcal{U}_PG_hG_{\ell};$ and $K_h = K \cap G_h$; cf. \S \ref{subsec-parabolics},
\ref{subsec-BB1}.  The Nomizu connections $\nabla_D^{\text{Nom}}$ (on $E$) and
$\nabla_1^{\text{Nom}}$ (on $E_1$) pass to  connections $\nabla_X^{\text{Nom}}$ on $E'\to
X$ and $\nabla_Y^{\text{Nom}}$ on $E'_Y \to Y$ respectively.  We use an inductive
procedure to define the {\it patched connection} $\nabla_Y^{\sf p}$ on the vectorbundle
$E'_Y \to Y,$ for any stratum $Y \le X$ as follows.  If $Y \subset \overline{X}$ is a
minimal stratum set $\nabla_Y^{\sf p} = \nabla_Y^{\text{Nom}}.$ Now suppose that the
patched connection $\nabla^{\sf p}_Z$ has been constructed on every stratum $Z < Y.$

\subsection{  Definition }\label{subsec-patched}
 The patched connection $\nabla^{\sf p}_Y$ on $E'_Y\to Y$ is the connection
\begin{equation}\label{eqn-patched}
  \nabla^{\sf p}_Y = B_Y^{\epsilon Y}\nabla_Y^{\text{Nom}}+ \sum_{Z<Y} B_Z^{\epsilon Y}
\Phi_{YZ}^*(\nabla_Z^{\sf p})\end{equation}
(where the sum is taken over all strata $Z<Y$ in the Baily Borel compactification  of
$X$).  
\subsection{Remarks}\label{subsec-explanation}
The idea behind this construction may be explained when there are two singular strata $Z
< Y < X.$  A simpler candidate for a connection on $X$
whose Chern forms might satisfy the $\pi$-fiber condition is
\begin{equation}\label{eqn-simpler} \nabla'_X = B_Z^{\epsilon X}
\Phi_{XZ}^*\nabla_Z^{\text{Nom}} + B_Y^{\epsilon X}
\Phi_{XY}^*\nabla_Y^{\text{Nom}} + B_X^{\epsilon X}\nabla_X^{\text{Nom}}.\end{equation}
In the region $T_Y(\epsilon X/2)$ only the first two terms contribute to $\nabla'_X$.
Both connection $\Phi_{XZ}^*\nabla_Z^{\text{Nom}}$ and $\Phi_{XY}^*\nabla_Y^{\text{Nom}}$
have curvature forms which satisfy the $\pi$-fiber condition with respect to $Y$.
However the curvature form of (and even the Chern forms of) any affine combination of
these fails to satisfy the $\pi$-fiber condition.  (cf. Figure \ref{Fig-neighborhoods}:
this occurs in the region where $B_Z+B_Y =1.$)  The remedy is to create a connection on
$X$ for which no nontrivial affine combination of $\Phi_{XZ}^*\nabla_Z^{\text{Nom}}$ and
$\Phi_{XY}^*\nabla_Y^{\text{Nom}}$ ever occurs.  Replacing $\nabla_Y^{\text{Nom}}$ by
$\nabla_Y^{\sf p}$ in (\ref{eqn-simpler}) gives
\[ \nabla_X^{\sf p} = B_Z^{\epsilon X}\Phi_{XZ}^*\nabla_Z^{\text{Nom}} + B_Y^{\epsilon X}
B_Z^{\epsilon Y}\Phi_{XY}^*\Phi_{YZ}^*\nabla_Z^{\text{Nom}} + B_Y^{\epsilon X}
B_Y^{\epsilon Y} \Phi_{XY}^* \nabla_Y^{\text{Nom}} + B_X^{\epsilon X}
\nabla_X^{\text{Nom}}\]
Within the region $T_Y(\epsilon X/2)$ only the first three terms appear:  the first term
alone appears in $T_Z(\epsilon X/2)$; the first and second terms appear in $T_Z(\epsilon
X) - T_Z(\epsilon X/2)$; the second term alone appears in $T_Z(\epsilon Y/2) -
T_Z(\epsilon X)$; the second and third terms appear in $T_Z(\epsilon Y) - T_Z(\epsilon
Y/2)$ and the third term alone appears outside $T_Z(\epsilon Y).$  In the region
$T_Z(\epsilon Y) - T_Z(\epsilon Y/2),$ 
\[ \nabla^{\sf p}_X = \Phi_{XY}^*(B_Z^{\epsilon Y}
\Phi_{YZ}^*\nabla_Z^{\text{Nom}} + B_Y^{\epsilon Y} \nabla_Y^{\text{Nom}}).\]
So $\nabla^{\sf p}_X$ is parabolically induced from a connection on $Y$, and by Corollary
\ref{cor-parab} its curvature form satisfies the $\pi$-fiber condition relative to $Y$.
In the region $T_Z(\epsilon X) - T_Z(\epsilon X/2)$, 
\[ \nabla^{\sf p}_X = (B_Z^{\epsilon X}\Phi_{XZ}^* + B_Y^{\epsilon X}B_Z^{\epsilon Y}
\Phi_{XY}^*\Phi_{YZ}^*)(\nabla_Z^{\text{Nom}}).\]
In this region, the curvature form still does not satisfy the $\pi$-fiber condition
however we show in \S \ref{subsec-completion} that the difference  $(\Phi_{XZ}^* -
\Phi_{XY}^*\Phi_{YZ}^*)\omega_Z^{\text{Nom}}$ is nilpotent and commutes with the
curvature form.  This turns out to be enough (Lemma \ref{lem-Springer}) to imply that the
Chern forms of $\nabla_X^{\sf p}$ satisfy the $\pi$-fiber condition with respect to $Y$.
\subsection{} Returning to the general case, suppose $\mathbf{Z} = Z(1) < Z(2) < \cdots <
Z(r)$ is a chain of strata in $\overline{X}.$  Write $\mathbf Z \le Y$ if $Z(r) =Y$, that
is, if the chain ends in $Y$.  Write $Y \le \mathbf Z$ if $Z(1) = Y$, that is, if the
chain begins at $Y$.  Suppose $\mathbf Z$ is such a chain of strata and suppose $x\in
T_{Z(i)}(\epsilon_0)$ for $1 \le i \le r.$  Denote by 
\begin{align*}\epsilon i &= \epsilon(Z(i)) = \epsilon_0/2^{\dim Z(i)}\\  
B_i^{\epsilon} &= B_{Z(i)}^{\epsilon}\\
\Phi_{ji}^* &= \Phi_{Z(j)Z(i)}^* \text{ for }j > i\\
\pi_i(x) &= \pi_{Z(i)}(x). \end{align*}
Define
\begin{align}\label{eqn-BZ}
B_{\mathbf Z}(x) &=  B_{r-1}^{\epsilon r}(\pi_r(x)) \ldots B_{2}^{\epsilon 3}(\pi_3(x))
B_1^{\epsilon 2}(\pi_2(x)) \\
\Phi_{\mathbf Z}^* &= \Phi_{r,r-1}^*\ldots \Phi_{32}^* \Phi_{21}^*.
\end{align}
If $\mathbf Z = \{Y\}$ consists of a single element, set $B_{\mathbf Z}(x)=1$ and
$\Phi_{\mathbf Z}^* = \text{Id}.$  The following lemma is easily verified by induction.

\begin{lem}\label{lem-chains}
The patched connection may be expressed as follows,
\begin{equation}\label{eqn-patched2} \nabla_Y^{\sf p}(x) = \sum_{\mathbf Z} B_{\mathbf
Z}(x) \Phi_{\mathbf Z}^* \left( B_{Z(1)}^{\epsilon
Z(1)}(\pi_{Z(1)}(x))\nabla_{Z(1)}^{\text{\rm Nom}}\right) \end{equation}
where the sum is over all chains of strata $\mathbf Z \le Y$ ending in $Y$.  \qed
\end{lem} 
\noindent (One checks that, although the projection functions $x\mapsto \pi_{Z(i)}(x)$
are not everywhere defined, they occur in (\ref{eqn-patched2}) with coefficient $0$
unless $x$ lies in the region of definition.) 

\quash{
\subsection{Proof}
Using induction, replace $\nabla_W^{\sf p}$ in (\ref{eqn-patched}) by the corresponding
expression (\ref{eqn-patched2}) to obtain
\begin{align}\label{eqn-intermediate}
 \nabla_Y^{\sf p}(x) = &B_Y^{\epsilon Y}(x) \nabla_Y^{\text{Nom}}(x) \\
 &+ \sum_{W<Y}\sum_{\mathbf Z'
\le W}B_W^{\epsilon Y}(x)B_{\mathbf Z'}(\pi_W(x))\Phi_{YW}^*\Phi_{W\mathbf
Z'}^*(B_{Z'(1)}^{\epsilon Z'(1)}(\pi_{Z'(1)}(x)\nabla_{Z'(1)}^{\text{Nom}}).\end{align}
In (\ref{eqn-patched2}) the trivial chain $\mathbf Z = \{Y\}$ gives rise to a term
$B_Y^{\epsilon Y}(\pi_Y(x))\nabla_Y^{\text{Nom}}.$  Each nontrivial chain may be
expressed as
$\mathbf Z = \mathbf Z' < Y,$ that is, $\mathbf Z = \{ Z'(1) < \cdots < Z'(r)=W <Y$ for
some stratum $W<Y.$  Then $B_{\mathbf Z}(x) = B_W^{\epsilon Y}(x) B_{\mathbf
Z'}(\pi_W(x))$ and $\Phi_{\mathbf Z}^* = \Phi_{YW}^*\Phi_{\mathbf Z'}^*.$  When
substituted into (\ref{eqn-patched 2}) these give (\ref{eqn-intermediate}), which
completes the inductive step in the argument.  \qed
}  

Definition \ref{subsec-patched} constructs a patched connection $\nabla_X^{\sf
p}$ on each of the automorphic vectorbundles $E'_Y \to Y.$  The proof of the following
theorem will appear in \S\ref{sec-proof}.
\begin{thm}\label{thm-A}
  The Chern forms $\{\sigma^j(\nabla^{\sf p}_Y)\}_{Y \le X}$of the patched connection
$\nabla^{\sf p}_X$ constitute a closed $\pi$-fiber differential form.
\end{thm}

\begin{cor}\label{cor-lift}  For each $j$, the Chern form $\sigma^j(\nabla^{\sf p}_{X})
\in \mathcal A^{2j}_{\pi}(\overline{X};\mathbb C)$ of the patched connection determines a
lift
\begin{equation}
\bar c^j(E')=[\sigma^j(\nabla^{\sf p}_{X})] \in H^{2j}(\overline{X};\mathbb C)
\end{equation}
of the Chern class $c^j(E')\in H^{2j}(X;\mathbb C)$ which is independent of the choices
that were made in its construction.  For any stratum closure $i:\overline{Y}
\hookrightarrow \overline{X}$ the restriction $i^*\bar c^j(E')$ is equal to the Chern
class $\bar c^j(E'_Y)\in H^{2j}(\overline{Y};\mathbb C)$ of the automorphic vectorbundle
$E'_Y \to Y.$ \end{cor}

\subsection{Proof of Corollary \ref{cor-lift}}  The restriction map
$i^*:H^{2j}(\overline{X}) \to H^{2j}(X)$ associates to any $\pi$-fiber differential form
$\omega \in \mathcal A^{2j}_{\pi}(\overline{X})$ the cohomology class $[\omega_X]$ of the
differential form $\omega_X \in \mathcal A^{2j}(X)$ on the nonsingular part.  Hence
$i^*(\bar c^j(E')) = c^j(E') \in H^{2j}(X;\mathbb C)$ since the latter is independent of
the connection.

The patched connection $\nabla^{\sf p}_X$ depends on the choice of a pair (partition of
unity, control data which is subordinate to the canonical projections $\{\pi_Z\}$) (see
\S \ref{sec-partition} and \S  \ref{sec-controldata}).  It is tedious but standard to
check that two such choices are connected by a smooth 1-parameter family of choices
(partition of unity, control data subordinate to $\{\pi_Z\}$).  The resulting patched
connections $\nabla_0^{\sf p}$ and $\nabla_1^{\sf p}$ are therefore connected by a smooth
1-parameter family of patched connections $\nabla_t^{\sf p}$, each of whose Chern forms
is a $\pi$-fiber differential form.  So the usual argument (e.g. \cite{KN} Chapt. XII
lemma 5; \cite{Milnor}) produces a differential $2j-1$ form $\Psi$ such that
$\sigma^j(\nabla_1^{\sf p}) - \sigma^j(\nabla_0^{\sf p}) = d\Psi.$  It is easy to see
that $\Psi\in \mathcal A^{2j-1}_{\pi}(\overline{X};\mathbb C)$
is a $\pi$-fiber differential form.  Consequently the $\pi$-fiber cohomology
classes coincide:  $[\sigma^j(\nabla_1^{\sf p})] = [\sigma^j(\nabla_0^{\sf p})] \in
H^{2j}(\overline{ X};\mathbb C).$  The second statement follows from the analogous
statement in Theorem \ref{thm-mythesis}\qed

\subsection{Remarks}\label{subsec-remark}  Theorem \ref{thm-A} and Corollary
\ref{cor-lift} extend to the case that $\mathbf G$ is semisimple over
$\mathbb Q.$  The restriction map $H^{2j}(\overline{X};\mathbb C) \to H^{2j}(X;\mathbb
C)$ factors as follows,
\begin{equation*}
H^{2j}(\overline{X};\mathbb C) \to IH^{2j}(\overline{X};\mathbb C) \to
H_{2n-2j}(\overline{X};\mathbb C) \to H_{2n-2j}(\overline{X},\partial\overline{X};\mathbb
C) \cong H^{2j}(X;\mathbb C)\end{equation*}
where $2n=\dim_{\mathbb R}(X)$, and where $\partial\overline{X}$ denotes the singular set
of the Baily-Borel compactification $\overline{X}.$  The Chern class $c^j(E') \in
H^{2j}(X;\mathbb C)$ lives in the last group.  For any toroidal resolution of
singularities $\tau: \overline{X}_{\Sigma} \to \overline{X}$, the pushdown
\[ c_{n-j}(E')=\tau_*(c^j(\bar E'_{\Sigma})\cap
[\overline{X}_{\Sigma}]) \in H_{2n-2j}(\overline{X};\mathbb Z)\] of the Chern class of
Mumford's canonical extension (\cite{Mumford}) $\bar E'_{\Sigma}$ of $E'$ gives a
canonical lift of  $c^j(E')$ to the homology of the Baily-Borel
compactification. In \S\ref{sec-constructible} (in the case of the tangent bundle) we
identify this with the (homology) Chern
class of the constructible function $\mathbf 1_X.$  In \cite{BBF} it is shown that every
algebraic homology class (including $c_{n-j}(E')$) admits a (non-canonical) lift to
middle intersection homology with rational coefficients.

\section{Proof of Theorem \ref{thm-A}  }\label{sec-proof}
\subsection{Preliminaries}\label{subsec-preliminaries}
  As in \S \ref{subsec-BB1}, let $q:D^{\star} \to \overline{X}$
denote the projection.  If $\epsilon_0 >0$ is sufficiently small, then for each stratum
$Z$ of $\overline{X}$ the preimage 
\[q^{-1}(T_Z(\epsilon_0)) = \coprod_{q(D_1)=Z} U_{D_1}(\epsilon_0)\] 
is a disjoint union of $\Gamma$-parabolic neighborhoods $U_{D_1}(\epsilon_0)$ of
those boundary components $D_1$ such that $q(D_1) = Z.$  For such a boundary
component define \begin{equation*}
\chi_{D_1}(x) = \begin{cases} 1 \text{ if } x \in U_{D_1}(\epsilon_0) \\
                                   0 \text{ otherwise}
\end{cases}\end{equation*}
to be the characteristic function of $U_{D_1}.$ 

Fix a stratum $Y$ and a choice $D_2$ of rational boundary component such that
$q(D_2)=Y.$.   Denote by $P_{2}=\mathcal U_{2} G_{2h}G_{2\ell}$ the
rational maximal parabolic subgroup of $G$ which normalizes $D_2$ and by $\nu_h:P_2
\to G_{2h}$ the projection as in \S \ref{subsec-BB1}.  By Proposition \ref{prop-Cayley}
the representation $\lambda|K \cap P_2$ extends to a representation $\lambda_2$ of
$K_{2h}G_{2\ell}.$  Set $\Gamma_h =\nu_h(\Gamma \cap P_2).$

The partition of unity $B_Y^{\epsilon Y} + \sum_{Z < Y}B_Z^{\epsilon Y} = 1$ on $Y$ pulls
back to a $\Gamma_h$-invariant  locally finite partition of unity on $D_2,$
\[  B_{D_2}^{\epsilon Y} +  \sum_{D_1 \prec D_2}  B_{D_1}^{\epsilon Y} =1\]
where the sum is over all rational boundary components $D_1 \prec D_2$, where
$B_{D_1}^{\epsilon Y} = q^*(B_Z^{\epsilon Y})\chi_{D_1}$ (and similarly for
$B_{D_2}^{\epsilon Y}$, however $\chi_{D_2}=1$ on $D_2$). 
The patched connection $\nabla^{\sf p}_Y$ on the vectorbundle $E'_Y \to Y$ pulls
back to a $\Gamma_h$-invariant connection $\nabla^{\sf p}_{2}=q^*(\nabla^{\sf p}_Y)$ on
the homogeneous vectorbundle $E_{2} = G_{2h} \times_{K_{2h}} V.$
This connection may also be described as the affine locally finite combination
\begin{equation}\label{eqn-bigpatch}
\nabla^{\sf p}_{2} = B_{D_2}^{\epsilon Y} \nabla_{2}^{\text{Nom}}
+ \sum_{D_1 \prec D_2}  B_{D_1}^{\epsilon Y} \Phi^*_{21}(\nabla_{1}^{\sf p})
\end{equation}
where, for each rational boundary component $D_1 \prec D_2$ the obvious notation holds:
$\nabla^{\sf p}_1$ is the patched connection on $E_1 \to D_1$ and $\Phi_{21}: E_2 \to
\pi_{21}^*(E_1)$ is the vectorbundle isomorphism which is obtained from
(\ref{eqn-Phiiso}) upon replacing $G$ by $G_{2h}.$

Denote by $\omega^{\sf p}_{2} \in \mathcal A^1(G_{2h}, \text{End}(V))$  the connection
form of $\nabla^{\sf p}_{2}.$  The curvature form $\Omega^{\sf p}_Y \in \mathcal A^2(Y,
\text{End}(E'_Y))$ of $\nabla^{\sf p}_Y$ coincides with the curvature form $\Omega^{\sf
p}_{2} \in \mathcal A^2_{\text{bas}} (G_{2h}, \text{End}(V))$ of $\nabla^{\sf p}_{2}$
under the canonical isomorphism
\begin{equation}\label{eqn-basicgamma} \mathcal A^2(Y,\text{End}(E'_Y)) \cong \mathcal
A^2(D_2, \text{End}(E_{2}))^{\Gamma_h} \cong \mathcal
A^2_{\text{bas}}(G_{2h},\text{End}(V))^{\Gamma_h} \end{equation}
where the superscript $\Gamma_h$ denotes the $\Gamma_h$-invariant differential forms.

\begin{prop}\label{prop-commutes}
Let $D_2$ be a rational boundary component of $D = G/K.$  
Then the connection form $\omega^{\sf p}_{2} \in \mathcal A^1(G_{2h},
\text{End}(V))$ and the curvature form $\Omega^{\sf p}_{2}\in \mathcal
A^2_{\text{bas}}(G_{2h}, \text{End}(V))$ commute with the adjoint action of
$\lambda_2(G_{2\ell}) \subset GL(V).$
\end{prop}

\subsection{Proof}  The proof uses a double induction over boundary components $D_2$ in
$D^*.$  However, so as to avoid the horribly complicated notation which would arise in
the proof, we rephrase the double induction as follows:
\begin{enumerate}
\item By induction we assume the proposition has been proven for every rational boundary
component $D_2'$ of any Hermitian symmetric domain $D' = G'/K$ for which
$\text{dim}(D') < \dim(D).$  (The case $\dim(D')=0$ is trivial.)
\item  For our given domain $D$, we assume the proposition has been proven for every
rational boundary component $D_1$ of $D$ such that $\dim(D_1) < \dim (D_2).$  (The case
$D_1 = \phi$ is trivial.) \end{enumerate}

To prove Proposition \ref{prop-commutes} for $D_2 \subset D^*$,
it suffices to verify that $\lambda_2(G_{2 \ell})$ commutes with the
connection form of each of the connections appearing in the linear combination
(\ref{eqn-bigpatch}).  The connection form $\omega^{\text{Nom}}_{2}$ of the Nomizu
connection $\nabla^{\text{Nom}}_{2}$ is given by
\[ \omega^{\text{Nom}}_{2}(L_{g_{2h}*})(\dot g_{2h}) = \lambda'(\dot k_{2h})
= \lambda'_2(\dot k_{2h})\]
for any $g_{2h} \in G_{2h}$ and for any $\dot g_{2h} \in \mathfrak g_{2h}$, where
$\dot g_{2h} = \dot k_{2h} + \dot p_{2h}$ is its Cartan decomposition.  This commutes
with $\lambda_2(G_{2\ell})$ since $G_{2 h}$ and $G_{2\ell}$ commute.  

Now consider any boundary component $D_1 \prec D_2$ which appears in the sum
(\ref{eqn-bigpatch}). Let $P_{1} = \mathcal U_{1}G_{1h}G_{1\ell}$ be the rational
parabolic subgroup which normalizes $D_1.$  Decompose the intersection
\[P = P_{1} \cap P_{2} = 
\mathcal U_{2}(\overline{\overline{\mathcal U}}G_{1h} G'_{\ell}) G_{2\ell}\]
according to (\ref{eqn-twopb}), setting $P_h = \overline{\overline{\mathcal U}}G_{1h}
G'_{\ell}\subset G_{2h}.$  Let $\nabla^{\sf p}_1$ be the patched connection on $E_1 \to
D_1.$  According to Proposition \ref{prop-induced}, the
connection form $\omega_{21}$ of the parabolically induced connection $\Phi^*_{21}
(\nabla_{1}^{\sf p})$ is given by
\[ \omega_{21}(L_{g*}(\dot{\bar{\bar u}}) + \dot g_{1 h} + \dot g'_{\ell}) =
\lambda'_{21}(\dot g'_{\ell}) + Ad(\lambda_{21}(g'_{\ell})^{-1})(\omega_{1}^{\sf
p}(L_{g_{1h*}} (\dot g_{1h}))).\]
Here, $g = \bar{\bar u} g_{1h} g'_{\ell} \in P_h$ and $\dot{\bar{\bar u}}
+\dot g_{1h} + \dot g'_{\ell} \in \text{Lie}(P_h)$ and $\lambda_{21}:K_{1h}G'_{\ell} \to
\text{GL}(V)$ is as in Proposition \ref{prop-Cayley}.   Since $\dim Y < \dim X$ we
may apply the first induction hypothesis and conclude that the adjoint action of
$\lambda_{21}(G'_{\ell})$ commutes with the connection form $\omega_{1}^{\sf p} \in
\mathcal A^1(G_{1h}, \text{End}(V)).$ Hence, using Proposition \ref{prop-finer}, 
\[\omega_{21}(L_{g*}(\dot{\bar{\bar u}}) + \dot g_{1 h} + \dot g'_{\ell}) =
\lambda'_1(\dot g'_{\ell}) + \omega_{1}^{\sf p}(L_{g_{1h*}}(\dot g_{1h})).
\]
The group $G'_{\ell}$ commutes with $G_{2\ell}$ so the first term
$\lambda'_1(\dot g'_{\ell})$ commutes with $\lambda_1(G_{2\ell}).$  By the second
induction hypothesis, the connection form $\omega_{1}^{\sf p} \in \mathcal
A^1(G_{1h}, \text{End}(V))$ also commutes with $\lambda_1(G_{1\ell}).$  But
$\lambda_1(G_{1\ell}) \supset \lambda_1(G_{2 \ell})= \lambda_2(G_{2\ell})$ by Proposition
\ref{prop-finer} again, which completes the proof
that the connection form of the patched connection commutes with $\lambda_2(G_{2\ell}).$
\qed

\subsection{}
Let $x\in X$ be a point near the boundary of $\overline{X}.$  Then there is a maximal
collection of strata $Y_1, Y_2,\ldots, Y_t$ such that $x \in T_{\epsilon_0}(Y_i)$ for
each $i$, and we may assume they form a partial flag, 
\begin{equation}\label{eqn-flag}
Y_1 < Y_2 < \cdots < Y_t = X.\end{equation}
Let $W\le X$ be the largest stratum in this collection such that $B_W^{\epsilon
W}(\pi_W(x)) \ne 0.$  Such a stratum exists since $B_{Y_1}^{\epsilon
Y_1}(\pi_{Y_1}(x))=1.$  (Choosing $W$ in this way guarantees that, if $W \ne X$, then for
every stratum $Z>W$ in this partial flag, at the point $\pi_Z(x)$ the connection
$\nabla^{\sf p}_Z$ is an affine combination of connections induced from smaller strata
and contains no contribution from $\nabla_Z^{\text{Nom}},$ because $B_Z^{\epsilon
Z}(\pi_Z(x))=0.$)

\begin{prop}\label{prop-patchedconnection}
At the point $x\in X$, $\sum_{W \le \mathbf S \le X}B_{\mathbf S}\ts{(x)}=1$ and
\quash{the connection $\nabla_Y^{\sf p}=d + \omega_Y^{\sf p}$  is given by}
\begin{equation}\label{eqn-patchedconnection}
\nabla_X^{\sf p}\ts{(x)} 
= \left( \sum_{W \le \mathbf S \le X}B_{\mathbf S}\ts{(x)} \Phi_{\mathbf S}^*  \right)
\left( \nabla_W^{\sf p}\ts{(\pi_W(x))}  \right)
\end{equation}
where the sum is over sub-chains $\mathbf S$ in the partial flag (\ref{eqn-flag}) which
begin at $W$ and end at $X$.
\end{prop}

\subsection{Proof}
By Lemma \ref{lem-chains} applied to $\nabla^{\sf p}_W$, we need to show that
\begin{equation}\label{eqn-bigproduct}\nabla_X^{\sf p}(x) =
 \left( \sum_{W \le \mathbf S \le X}B_{\mathbf S}\ts{(x)} \Phi_{\mathbf S}^*  \right) 
\left( \sum_{\mathbf R \le W}B_{\mathbf R}\ts{(\pi_W(x))}\Phi_{\mathbf R}^* 
\left( B_{R(1)}^{\epsilon R(1)} \nabla_{R(1)}^{\text{\rm Nom}}
\ts{(\pi_{R(1)}(x))})  \right)\right).
\end{equation}
By Lemma \ref{lem-chains}, $ \nabla_X^{\sf p}\ts{(x)}$ is a sum over chains $\mathbf Z
\le X$ of terms 
\[ B_{\mathbf Z}\ts{(x)} \Phi_{\mathbf Z}^*B_{Z(1)}^{\epsilon Z(1)}
\nabla_{Z(1)}^{\text{Nom}}\ts{(\pi_{Z(1)}(x))}.\]
For any $\epsilon\le \epsilon_0$, $B_Z^{\epsilon}(x)=0$ unless $Z$ occurs in the
collection $\{Y_1,Y_2,\ldots, Y_t=X\}.$   By assumption the term $B_{Z(1)}^{\epsilon
Z(1)}\pi_{Z(1)}(x)$ also vanishes unless $Z(1) \le W.$  Therefore each chain $\mathbf Z =
Z(1) < \cdots < X$ appearing in the sum may be assumed to occur as a sub-chain of $Y_1 <
Y_2 < \cdots <Y_t=X$, and  we may also assume the chain begins at $Z(1) \le W.$  We claim
that if such a chain $\mathbf Z$ occurs with nonzero coefficient, then the stratum $W$
must appear in the chain.  For if not, then  $Z(k) < W < Z(k+1)$ for some $k.$  But the
term $B_{\mathbf Z}(x)$ contains a factor
\[ B_{Z(k)}^{\epsilon Z(k+1)}(\pi_{Z(k+1)}(x)).\]
Since $B_W^{\epsilon W}(\pi_W(x)) \ne 0$, this factor vanishes by Lemma
\ref{lem-vanishing}, which proves the claim.  Summarizing, every chain $\mathbf Z$ which
occurs with nonzero coefficient in the sum may be described as
\[ R(1) < R(2) < \cdots < R(r) = W = S(1) < S(2) < \cdots < S(s)=X.\]  The contribution 
to $\nabla_X^{\sf p}(x)$ in (\ref{eqn-patched2}) from such a chain is the product of
\[ B_{S(s-1)}^{\epsilon X}\ts{(\pi_{X}(x))} \ldots B_W^{\epsilon
S(2)}\ts{(\pi_{S(2)}(x))} \Phi_{XS(s-1)}^* \ldots \Phi_{S(2)W}^* \]
with
\[ B_{R(r-1)}^{\epsilon W}\ts{(\pi_W(x))} \ldots B_{R(1)}^{\epsilon R(2)}
\ts{(\pi_{R(2)}(x))}
\Phi_{WR(r-1)}^*\ldots \Phi_{R(2)R(1)}^* \]
applied to
\[ B_{R(1)}^{\epsilon R(1)} \nabla_{R(1)}^{\text{Nom}}\ts{(\pi_{R(1)}(x))}.\]
However this product is exactly a single term in (\ref{eqn-bigproduct}) and every such
product occurs exactly once, which verifies (\ref{eqn-patchedconnection}).    Since the
coefficients in (\ref{eqn-patched2}) sum to 1, so also does $\sum_{W \le \mathbf S \le
X}B_{\mathbf S}\ts{(x)}. $\qed
\quash{
the first factor is a single term in 
\[ \sum_{W \le \mathbf S \le Y} B_{\mathbf S}(y) \Phi_{\mathbf S}^*\]
while the second factor is a single term in 
\[ \sum_{\mathbf R \le W} B_{\mathbf R}\ts{(\pi_W(y))}\Phi_{\mathbf R}^* \]
and every such product appears exactly once.  \qed
} 

\subsection{}
We must prove that the Chern forms of the patched connection satisfy the $\pi$-fiber
condition near each stratum of $\overline{X}.$  Let $Y$ be such a stratum and let  $x
\in X \cap T_Y(\epsilon X/2).$  We will verify the $\pi$-fiber condition relative to the
stratum $Y$ at the point $x.$  The point $x\in X$ lies in an intersection of
$\epsilon_0$-tubular neighborhoods of a maximal collection of strata, which (we may
assume) form a partial flag $Z_1 < Z_2 < \cdots < X.$ Let $W$ be the largest stratum in
this chain such that $B_W^{\epsilon W}(\pi_W(x)) \ne 0.$  Then $W \le Y$ by
(\ref{eqn-Bpi00}).  Consider the subchain lying between $W$ and $Y$, which we shall
denote by
\[ W = Y_1 < Y_2 < \cdots < Y_t = Y.  \]
Fix corresponding boundary components $D_1 \prec {D}_2 \prec \cdots
\prec {D}_t$ and let $P_1\prec P_2\prec \ldots \prec P_t$ be their normalizing maximal
parabolic subgroups. Set $P= P_1 \cap P_2 \cap \ldots \cap P_t.$    Let $\nabla^{\sf
p}_W$ be the patched connection on $E'_W \to W.$ According to Proposition
\ref{prop-patchedconnection},
\begin{equation}\label{eqn-recall} \nabla_X^{\sf p}(x) = \left( \sum_{W \le \mathbf S \le
X} B_{\mathbf S}(x) \Phi_{\mathbf S}^* \right) \nabla_W^{\sf p}(\pi_W(x)).
\end{equation}
However the only chains $\mathbf S= \{S(1)<\cdots < S(s)\}$ which occur with nonzero
coefficient in this sum satisfy
\begin{equation}\label{eqn-Scondition}
\{S(1),S(2),\ldots, S(s-1)\} \subset \{ Y_1,Y_2,\ldots, Y_t\},\ S(1)=W=Y_1,\
S(s)=X\end{equation}
for the following reason.  Suppose a chain $\mathbf S$ contains a stratum larger than $Y$
(but not equal to $X$). Let $Z$ be the largest such stratum occurring in $\mathbf S$.
Then the first factor in $B_{\mathbf S}(x)$ is $B_Z^{\epsilon X}(x),$ cf. (\ref{eqn-BZ}).
By assumption, $x \in T_Y(\epsilon X/2).$  So by (\ref{eqn-Bpi00}) (with $\epsilon =
\epsilon X$ and where the roles of $Y$ and $Z$ are reversed), $B_Z^{\epsilon X}(x)=0.$

If $W=Y$ then (\ref{eqn-recall}) becomes $\nabla^{\sf p}_X = \Phi_{XY}^*(\nabla^{\sf
p}_Y)$ so by Corollary \ref{cor-parab} and Proposition \ref{prop-commutes} the curvature
form $\Omega_X^{\sf p}$ of $\nabla_X^{\sf p}$ satisfies the $\pi$-fiber condition with
respect to $Y$.  So the same is true of every Chern form which proves Theorem \ref{thm-A}
in this case.  

Therefore we may assume that $W < Y$.  As in \S \ref{subsec-preliminaries} the connection
$\nabla^{\sf p}_W$ on $E'_W \to W$ pulls back to a $\Gamma_h$-invariant connection
$\nabla^{\sf p}_{1}$ on $E_1 \to D_1$ and the connection $\nabla^{\sf p}_X$ on $E' \to
X$ pulls back to a $\Gamma$-invariant connection $\nabla^{\sf p}$ on $E \to D.$  As in
(\ref{eqn-basicgamma}) we identify the curvature form $\Omega^{\sf p}_X$ of
$\nabla^{\sf p}_X$ with the curvature form $\Omega^{\sf p}$ of $\nabla^{\sf p}$ under the
canonical isomorphism
\[ \mathcal A^2(X,\text{End}(E')) \cong \mathcal A^2(D, \text{End}(E))^{\Gamma} \cong
\mathcal A^2_{\text{bas}}(G, \text{End}(V))^{\Gamma}.\]
Similarly identify the curvature $\Omega^{\sf p}_1$ of $\nabla^{\sf p}_1$ with the
curvature $\Omega^{\sf p}_W$ of $\nabla^{\sf p}_W.$  Choose a lift $\tilde x \in D$ of
$x$, which lies in the intersection 
\[U_{D_1}(\epsilon_0) \cap U_{D_2}(\epsilon_0) \cap \ldots \cap U_{D_t}(\epsilon_0)\]
of $\Gamma$-parabolic neighborhoods of the boundary components $D_1 \prec D_2 \prec
\ldots \prec D_t.$  Let $\pi:D \to D_1$ denote the canonical projection.

\begin{lem}\label{lem-errorterm}
For any tangent vectors $U,V \in T_{\tilde x}D,$
\[ \Omega^{\sf p}(U,V) = \pi^*\Omega^{\sf p}_1(U,V) + n \in \text{End}(V)
\]  
for some nilpotent element $n \in \lambda'_1(\mathfrak g_{1\ell}).$  (Here, $P_1 =
\mathcal U_{P_1}G_{1h}G_{1\ell};$  $\mathfrak g_{1\ell} = \text{Lie}(G_{1\ell});$
and $\lambda_1$ is the extension (Proposition \ref{prop-Cayley}) of the representation
$\lambda|K\cap P_1$.)
\end{lem}

\subsection{Proof}
We will compute the curvature $\Omega_X^{\sf p}$ of $\nabla_X^{\sf p}.$  Each chain
$\mathbf S$ satisfying (\ref{eqn-Scondition}) corresponds also to a chain of rational
boundary components $D_1 = D_{\mathbf S(1)} \prec D_{\mathbf S(2)} \prec \ldots \prec
D_{\mathbf S(s)} = D.$  Let $\widetilde{B}_{\mathbf S} = q^*B_{\mathbf S}$ 
denote the pullback to $D.$ It follows from Proposition \ref{prop-patchedconnection} that
$\sum_{\mathbf S}\widetilde B_{\mathbf S}(\tilde x) = 1$ (where the sum is over all
chains $\mathbf S$ which appear in (\ref{eqn-recall})). Choose any
ordered labeling of the chains $\mathbf S$ which appear in (\ref{eqn-recall}),
say $\mathbf S_1, \mathbf S_2, \ldots, \mathbf S_M.$ By \S \ref{subsec-controldata} (4), 
\[\sum_{i=1}^M B_{\mathbf S_{i}}(x) \pi_{\mathbf S_i}^*(\Omega_W^{\sf p}) =
\sum_{i=1}^M B_{\mathbf S_{i}}(x)\pi_{XW}^*\Omega_W^{\sf p}= \pi_{XW}^*
\Omega_W^{\sf p}.\] Let $U,V \in T_xX.$  Pulling back the equation (\ref{eqn-recall}) to
$D$ and using  Lemma \ref{lem-patch} gives: 
\begin{equation*}\label{eqn-bigcurvature}\begin{split}
\Omega^{\sf p}(U,V) = &\pi^*(\Omega_1^{\sf p})(U,V)+ \sum_{i=1}^{M-1} d\tilde B_{\mathbf
S_i}\wedge(\Phi_{\mathbf S_i}^*\omega_1^{\sf p} -
\Phi_{\mathbf S_M}^*\omega_1^{\sf p})(U,V)\\
 & - \sum_{i<j}\tilde B_{\mathbf S_i}(x)\tilde B_{\mathbf S_j}(x)
\left[ \Phi_{\mathbf S_i}^*\omega_1^{\sf p}(U) - \Phi_{\mathbf S_j}^*\omega_1^{\sf p}(U),
       \Phi_{\mathbf S_i}^*\omega_1^{\sf p}(V) - \Phi_{\mathbf S_j}^*\omega_1^{\sf p}(V)   
\right]  
\end{split}\end{equation*}
where $\Phi_{\mathbf S}^*\omega^{\sf p}_1$ denotes the connection form of $\Phi_{\mathbf
S}^*\nabla^{\sf p}_1.$  Let us compute this connection form.  
Suppose that $\mathbf S$ is a chain satisfying (\ref{eqn-Scondition}).
For $ 1 \le j \le s-1$ let $P_{S(j)}$ be the corresponding normalizing maximal parabolic
subgroup and set $Q = P_{S(1)} \cap P_{S(2)} \cap \ldots \cap P_{S(s-1)}.$  
Then $P \subset Q$ and $P^{\flat} = Q^{\flat} = P_1$ which implies (as in \S
\ref{subsec-twoparab}) that $Q \subset P \subset P_1$ have compatible decompositions,
\begin{align*}
 P_1 &= \mathcal U_1 G_{1h} G_{1\ell}  \\ 
 Q &= \mathcal U_1 G_{1h} (\mathcal U_{P_1Q} G_{Q\ell}) \text{ with } \mathcal
U_Q = \mathcal U_1 \mathcal U_{P_1Q} \\
P &= \mathcal U_1 G_{1h} \mathcal U_{P_1Q}(\mathcal U_{QP}G_{P\ell}) \text{ with }
\mathcal U_P = \mathcal U_1 \mathcal U_{P_1Q} \mathcal U_{QP}.
\end{align*}
Here, $\mathcal U_{QP}G_{P\ell}$ is the parabolic subgroup of $G_{Q\ell}$ determined by
$P \subset Q$.  We also note that
\begin{equation}\label{eqn-UP1P}
\mathcal U_{P_1P} = \mathcal U_{P_1Q} \mathcal U_{QP} \end{equation}
is the unipotent radical of the parabolic subgroup $\nu_{1\ell}(P) \subset G_{1\ell}$
determined by $P$ (where $\nu_{1\ell}: P_1 \to G_{1\ell}$ is the projection).
Let $\mathfrak N_{P_1P}$ denote its Lie algebra.

By Proposition \ref{prop-commutes}, the connection form $\omega_1^{\sf p}(\pi_W(x))\in
\mathcal A^1(G_{1h}, \text{End}(V))$ commutes with the adjoint action of
$\lambda_1(G_{1\ell}) \subset GL(V)$.  So we may apply  Proposition
\ref{prop-multipleinduction} to determine $\Phi_{\mathbf S}^*(\omega^{\sf p}_1).$  Let $g
= u_P g_{1h} g_{P\ell} \in P = \mathcal U_P G_{1h} G_{P\ell}$ and let
\[ \dot g = \dot u_1 + \dot g_{1h} + \dot u_{P_1Q} + \dot u_{QP} + \dot g_{P\ell} \in
\text{Lie}(\mathcal U_1 G_{1h} \mathcal U_{P_1Q} \mathcal U_{QP} G_{P\ell}). \]
Apply Proposition \ref{prop-multipleinduction} using $\dot u_Q = \dot u_1 + \dot
u_{P_1Q}\in \text{Lie}(\mathcal U_Q)$ and $\dot g_{Q\ell} = \dot u_{QP} + \dot g_{P\ell}
\in \mathfrak g_{Q\ell}$
to find:
\[ \Phi_{\mathbf S}^*(\omega_{1}^{\sf p})(L_{g*}(\dot g)) 
= \omega_1^{\sf p}(L_{g_{1h}*}\dot g_{1h})) + \lambda'_1(\dot u_{QP})+
\lambda'_1(\dot g_{P\ell}) \ .\]
Moreover, $\lambda'_1(\dot u_{QP}) \in
\lambda'_1(\mathfrak N_{P_1P}) \subset \lambda'_1(\mathfrak g_{1\ell}).$

Now suppose that $\mathbf R$ is another chain in the sum (\ref{eqn-recall}) which makes a
nonzero contribution $\Phi_{\mathbf R}^*\nabla_1^{\sf p} = d+\Phi_{\bf
R}^*(\omega_{1}^{\sf p})$ to the connection $\nabla^{\sf p},$ say, $W = R(1) < R(2) <
\cdots < R(r-1) < X.$  Then $Q' = P_{R(1)} \cap P_{R(2)} \cap \ldots \cap P_{R(r-1)}$ and
$P$ also have compatible decompositions:
\begin{align*}
Q' &= \mathcal U_1 G_{1h} (\mathcal U_{P_1Q'} G_{Q'\ell}) \text{ with } \mathcal
U_{Q'} = \mathcal U_1 \mathcal U_{P_1Q'} \\
P &= \mathcal U_1 G_{1h} \mathcal U_{P_1Q'}(\mathcal U_{Q'P}G_{P\ell})\text{ with }
\mathcal U_P = \mathcal U_1 \mathcal U_{P_1Q'} \mathcal U_{Q'P}.
\end{align*}
 The same element $\dot g \in \text{Lie}(P)$ decomposes as
\[ \dot g = \dot u_1 + \dot g_{1h} + \dot u_{P_1Q'} + \dot u_{Q'P} + \dot g_{P\ell}.\]
So the same argument gives
 $\Phi_{\mathbf R}^*(\omega_{1}^{\sf p})(L_{g*}(\dot g)) = \omega_1^{\sf
p}(L_{g_{1h}*}(\dot g_{1h})) +\lambda'_1(\dot u_{Q'P})+
\lambda'_1(\dot g_{P\ell}) .$  We conclude that:
\begin{equation}\label{eqn-nilpotentdifference}
 (\Phi_{\bf S}^* \omega_{1}^{\sf p} - \Phi_{\bf R}^*\omega_{1}^{\sf p})(L_{g*}(\dot g)) =
\lambda'_1(\dot u_{QP} - \dot u_{Q'P}) \in\lambda'_1(\mathfrak N_{P_1P}) \subset
\lambda'_1(\mathfrak g_{1\ell}).\end{equation}
 Consequently each term in the sum (\ref{eqn-bigcurvature}) (except for the first) lies
in $\mathfrak N_{P_1P}.$  \qed 

\subsection{Completion of the proof} \label{subsec-completion} 
Using Lemma \ref{lem-errorterm}, at the point $\tilde x$ we may write
\[ \Omega^{\sf p}(U,V) = \pi^*(\Omega^{\sf p}_1)(U,V) + n\]
where $n \in \lambda'_1(G_{1\ell})$ is nilpotent and in fact lies in $\mathfrak
N_{P_1P}.$
Moreover, by Proposition \ref{prop-commutes}, the curvature $\Omega^{\sf p}_1$ commutes
with $n.$  If $f:\text{End}(V) \to \mathbb C$ is any
$Ad$-invariant polynomial, it follows from Lemma \ref{lem-Springer} that $f(\Omega^{\sf
p}(U,V)) = f(\Omega^{\sf p}_1(\pi_*U,)\pi_*V)$ hence also $f(\Omega_X^{\sf
p}(U,V)) =f(\Omega_W^{\sf p}(\pi_{XW*}U, \pi_{XW*}V)).$  So it follows from
(\ref{eqn-polarization}) that the corresponding characteristic form satisfies the
$\pi$-fiber condition relative to $W$.  This completes the proof of Theorem \ref{thm-A}. 
\qed

\section{Toroidal compactification}\label{sec-toroidal}

\subsection{}
Throughout this section we assume that $G = \mathbf G(\mathbb
R)$ is the set of real points of a connected semisimple algebraic group $\mathbf G$
defined over $\mathbb Q$, that $D = G/K$ is a Hermitian symmetric space,
$\Gamma \subset \mathbf G(\mathbb Q)$ is a neat arithmetic group, $X = \Gamma
\backslash G /K$ is the corresponding locally symmetric space with Baily-Borel
Satake compactification $\overline{X} = \Gamma \backslash D^*.$   Fix a
representation $\lambda:K \to GL(V)$ on some complex vectorspace $V$ and let
$E = G\times_K V$ be the corresponding homogeneous vectorbundle on $D$, and
$E' = \Gamma \backslash E$ the automorphic vectorbundle on $X$.  Choose a
system of control data on the Baily-Borel compactification $\overline{X}$ and
a partition of unity as in \S \ref{sec-partition}, and let
$\nabla_X^{\sf p}$ denote the resulting patched connection on $E' \to X.$
For each $i$, the Chern form $\sigma^i(\nabla_X^{\sf p})\in
\Omega^{\bullet}_{\pi}(\overline{X};\mathbb C)$ (cf \S \ref{subsec-characteristicforms})
is a $\pi$-fiber differential form on $\overline{X}$ so it determines a cohomology class
$\bar c^i(E') =[\sigma^i(\nabla_X^{\sf p})] \in
H^{2i}(\overline{X};\mathbb C).$

We also fix a nonsingular toroidal compactification $\overline{X}_{\Sigma}.$
This corresponds to a $\Gamma$-compat\-ible collection of
simplicial polyhedral cone decompositions $\Sigma_F$ of certain self adjoint homogeneous
cones.  These compactifications were constructed in \cite{AMRT} and are reviewed in 
\cite{Harris3}, \cite{HZ1}, \cite{Chai}, \cite{Namikawa}.  In \cite{Mumford}, D. Mumford
shows that the automorphic vectorbundle $E'\to X$ admits a
canonical extension $\bar E'_{\Sigma}$ over the toroidal compactification
$\overline{X}_{\Sigma}.$  In \cite{Harris3} Theorem 4.2, M. Harris shows that Mumford's 
canonical extension coincides with Deligne's canonical extension \cite{Deligne}
(for an appropriately chosen flat connection with unipotent monodromy).

The identity mapping $X\to X$ has a unique continuous extension, $\tau :
\overline{X}_{\Sigma} \to \overline{X}$ of $X$, and this is a resolution of
singularities.

\begin{thm}\label{thm-toroidal}  The patched connection $\nabla_X^{\sf p}$ on $E' \to
X$ extends to a smooth connection $\overline{\nabla}^{\sf p}_{\Sigma}$ on $\bar
E'_{\Sigma} \to \overline{X}_{\Sigma}.$
Moreover for each $i$,
\begin{equation}\label{eqn-toroidal}
\tau^*\bar c^i(E')=\tau^*([\sigma^i(\nabla_X^{\sf p})]) =
[\sigma^i(\overline{\nabla}^{\sf p}_{\Sigma})]=
c^i(\bar E'_{\Sigma}) \in H^{2i}(\overline{X}_{\Sigma};\mathbb C). \end{equation}
\end{thm}
\noindent
The proof will appear in Section \ref{sec-pfthmtor}.   S. Zucker has pointed out that it
follows from mixed Hodge theory that the image of $\bar c^i(E')$ in
$\text{Gr}^W_{2i}H^{2i}(\overline{X};\mathbb C)$ is uniquely determined by
(\ref{eqn-toroidal}).

\subsection{Proportionality theorem}\label{subsec-prop}
Fix representations $\lambda_j:K\to GL(V_j)$ for $j=1,2,\ldots,r$ and fix nonnegative
integers $I = (i_1,i_2,\ldots i_r)$ with $i_1 + i_2 + \ldots +
i_r =n= \dim(D).$  For $j = 1,2,\ldots,r$ let $E'_j = \Gamma \backslash G\times_KV_j
\to X$ be the resulting automorphic vectorbundle on $X$ and let $\check{E}_j = 
G_u \times_{K}V_j$ be the corresponding vectorbundle on the
compact dual symmetric space $\check D = G_u/K$ (where $G_u$ is a compact real form of
$\mathbf G$ containing $K$).   Define ``generalized'' Chern numbers 
\begin{equation}
\check c^{I}(\lambda_1,\lambda_2,\ldots,\lambda_r) = (c^{i_1}(\check E_1)\cup
c^{i_2}(\check E_2) \cup \ldots \cup c^{i_r}(\check E_r))\cap [\check D] \in
\mathbb Z \end{equation}
and
\begin{equation}
(\bar c^{I}(\lambda_1,\lambda_2,\ldots,\lambda_r) = \bar
c^{i_1}(\nabla^{\sf p}_1) \cup \bar c^{i_2}(\nabla^{\sf p}_2) \cup \ldots
\cup c^{i_r}(\nabla^{\sf p}_r)) \cap [\overline{X}] \in \mathbb C 
\end{equation}
where $\nabla^{\sf p}_j$ denotes the patched connection on $E'_j \to X$ and where
$[\overline{X}]\in H_{2n}(\overline{X};\mathbb C)$ denotes the fundamental class of the
Baily-Borel compactification. Let $v(\Gamma)\in \mathbb Q$ denote the constant which
appears in the proportionality theorem of Hirzebruch \cite{Hirzebruch}, \cite{Mumford}.

\begin{prop}\label{prop-prop}  For any choice $\lambda_1,\lambda_2,\ldots,\lambda_r$ of
representations and for any partition  $I = (i_1,i_2,\ldots,i_r)$ of $n= \dim_{\mathbb
C}(D)$ we have 
\begin{equation} \bar c^{I}(\lambda_1,\lambda_2,\ldots,\lambda_r) = v(\Gamma)\check
c^{I}(\lambda_1,\lambda_2,\ldots,\lambda_r). \end{equation} \end{prop}

\subsection{Proof}  
Each of the vectorbundles $E'_j$ has a canonical extension $\bar E'_{j,\Sigma} \to
\overline{X}_{\Sigma}.$  The same proof as in \cite{Mumford} (which is the same proof as
in \cite{Hirzebruch}) (cf. \cite{letter}) shows that the Chern classes of these extended
bundles satisfy the proportionality formula
\begin{equation}
(c^{i_1}(\bar E'_{1,\Sigma})\cup c^{i_2}(\bar E'_{2,\Sigma})\cup \ldots \cup
c^{i_r}(\bar E'_{r,\Sigma})) \cap [\overline{X}_{\Sigma}] = v(\Gamma) \check
c^{I}(\lambda_1,\lambda_2,\ldots,\lambda_r). \end{equation} 
The result now follows immediately from Theorem \ref{thm-toroidal} .  \qed 
\quash{(and the projection
formula $\tau_*(\tau^*a \cap b) = a \cap \tau_*(b)$ for any $a \in H^*(\overline{X})$ and
any $b \in H_*(\overline{X}_{\Sigma}.$)   \qed} 

\section{Proof of Theorem \ref{thm-toroidal}}\label{sec-pfthmtor}
\subsection{}\label{pf-tor0}  Fix $\epsilon \le \epsilon_0.$
We claim that the partition of unity (\ref{eqn-partition}) $\sum_{Z \le X} B_Z^{\epsilon}
=1$ on $\overline{X}$ pulls back to a smooth partition of unity on
$\overline{X}_{\Sigma}.$  Fix a pair of strata $Y,Z$ of $\overline{X}.$  We must verify
that $\tau^*B_Z^{\epsilon}$ is smooth near $\tau^{-1}(Y).$  It can be shown that the
mapping $\tau: \overline{X}_{\Sigma} \to \overline{X}$ is a complex analytic morphism
between complex analytic varieties, as is the projection $\pi_Y:T_Y(\epsilon) \to Y.$  It
follows that the composition $\pi_Y\tau: \tau^{-1}T_Y(\epsilon) \to Y$ is a complex
analytic morphism between smooth complex varieties, so it is smooth.  If the stratum $Z$
is not comparable to $Y$ or if $Z > Y$ then $B_Z^{\epsilon}$ vanishes on
$T_Y(\epsilon/2)$ hence $\tau^*B_Z^{\epsilon}$ vanishes on $\tau^{-1}(T_Y(\epsilon/2))$
so it is smooth.  If $Z \le Y$ then by (\ref{eqn-Bpi1}) and (\ref{eqn-Bpi1a}),
$\tau^*B_Z^{\epsilon}(x) = B_Z^{\epsilon}\pi_Y\tau(x)$ for all $x\in
\tau^{-1}(T_Y(\epsilon/2))$ so $\tau^*B_Z^{\epsilon}$ is smooth in this open set.
  
\subsection{}\label{pf-tor1}  We may assume that $\mathbf G$ is simple over $\mathbb Q.$
The ``boundary'' $\overline{X}_{\Sigma} - X$ of the toroidal compactification has a
distinguished covering by open sets $U_Y$, one for each stratum $Y\subset \overline{X}$
of the Baily-Borel compactification,  such that $\tau(U_Y)\subset \overline{X}$ is a
neighborhood of $Y$, and for which the restriction $\bar E'_{\Sigma}|U_Y$  arises from an
automorphy factor (\cite{HZ1} \S 3.3). 
\quash{(Harris and Zucker
use the ``canonical automorphy factor for $P$'' of \cite{Harris1}; however
one may just as well use our automorphy factor of (\ref{eqn-autP}) (cf. \S
\ref{subsec-twoauto}).)}
\begin{prop}\label{prop-extension}
For any smooth connection $\nabla_Y$ on $(E'_Y,Y)$ the parabolically induced connection
$\Phi_{XY}^*(\nabla_Y)$ (which is defined only on $E'|(U_Y \cap X)= \bar E'_{\Sigma}
|(U_Y\cap X)$) extends canonically to a smooth connection (which we denote by
$\overline{\Phi}_{XY}^*(\nabla_Y)$) on $\overline{E}'_{\Sigma}| U_Y.$ \end{prop}
\subsection{Proof of Theorem \ref{thm-toroidal}}
We may assume that $\epsilon_0 > 0$ was chosen so small that $T_Y(\epsilon_0)\subset U_Y$
for each stratum $Y < X$ of $\overline{X}.$  By \S \ref{pf-tor0} the partition of unity
which is used to construct the patched connection $\nabla^{\sf p}_X$ extends to a smooth
partition of unity on $\overline{X}_{\Sigma}.$  Hence, using Proposition
\ref{prop-extension}, 
\[
\overline{\nabla}^{\sf p}_{\Sigma} = \tau^*B_X^{\epsilon_X} + \sum_{Y<X} \tau^*
B_Y^{\epsilon X} \overline{\Phi}_{XY}^*(\nabla_Y^{\sf p})\]
is a smooth connection on $\overline{E}'_{\Sigma} \to \overline{X}_{\Sigma}$ which
coincides with the patched connection $\nabla_X^{\sf p}$ on $E' \to X.$  Therefore its
Chern forms are smooth and everywhere defined and they restrict to the Chern forms of
$\nabla^{\sf p}_X.$ It follows from Lemma \ref{lem-integrate} that each Chern class
of $\overline{\nabla}^{\sf p}_{\Sigma}$ is the pullback of the corresponding Chern class
of $\nabla^{\sf p}_X.$  \qed

The remainder of \S \ref{sec-toroidal} is devoted to the proof of Proposition
\ref{prop-extension} which is essentially proven in  \cite{HZ1} (3.3.9) (following
\cite{Harris2}).  We will now verify the details.

\subsection{}  Let us recall the construction \cite{AMRT} of the toroidal
compactification.  Fix a rational boundary component $F$ with normalizing parabolic
subgroup $P = \mathcal UG_hG_{\ell}$ and let $Y = \Gamma_h \backslash F \subset
\overline{X}$ denote the corresponding stratum in the Baily-Borel
compactification of $X = \Gamma \backslash D.$  Let $Z_F = \text{Center}
(\mathcal U)$ and $\mathfrak z =
\text{Lie}(Z_F).$  The vectorspace $\mathfrak z$ is preserved under the adjoint
action of $G_{\ell}$ and contains a unique open orbit $C_F\subset \mathfrak
z$; it is a rationally defined self adjoint homogeneous cone.  Its Satake
compactification $C_F^*$ consists of $C_F$ together with all its rational
boundary components (in the Satake topology).  The toroidal compactification
is associated to a collection $\Sigma = \{\Sigma_F\}$ of rational polyhedral
cone decompositions of the various $C_F^*$ which are compatible under $\Gamma.$

Let $\check{D}$ denote the compact dual symmetric space (so $\check{D} =
\mathbf G(\mathbb C)/K(\mathbb C)P^{-}$ in the notation of Proposition \ref{prop-Cayley};
cf.\cite{AMRT}, \cite{Satake}) and let $\beta:D \to \check{D}$ denote the
Borel embedding.  Set $D_F = Z_F(\mathbb C).\beta(D).$  The domain $D_F$ is homogeneous
under $P.Z_F(\mathbb C)$ and it admits ``Siegel coordinates''
$D(F) \cong Z_F(\mathbb C) \times \mathbb C^a \times F$ in which the subset
$\beta(D) \subset D_F$ is defined by a certain well known inequality (\cite{AMRT} p. 239,
\cite{Satake} \S III (7.4)).  Now consider the commutative diagram \cite{HZ1} (1.2.5),
reproduced in figure \ref{fig-1}.
\begin{figure}[ht]
\begin{equation*}
\begin{CD}
D @>>> \Gamma'_F\backslash D & {\subset }& D_{F,\Sigma} & @>>{\varphi_{F,\Sigma}}>
\overline{X}_{\Sigma}\\
{\cap} && {\cap} && {\cap} && \\
D_F @>>{\Gamma'_F}>  M'_F &\subset & M'_{F,\Sigma} && \\
@V{\theta_2}VV @V{\pi_2}V{T_F}V @V{\pi_{2,\Sigma}}V{T_{F,\Sigma}}V  && \\
D_F/Z_F(\mathbb C) @>>> A_F @= A_F  &&\\
@V{\theta_1}VV @V{\pi_1}VV  &&&& \\
F @>{\Gamma_h}>> Y &&&&
\end{CD}
\end{equation*}\caption{Toroidal compactification}\label{fig-1}\end{figure}

Here, $\Gamma'_F = \Gamma \cap (G_h\mathcal U)$ and $M'_F =\Gamma'_F
\backslash D_F$.  The algebraic torus $T_F = Z_F(\mathbb C)/(\Gamma \cap Z_F)$
acts on $M'_F$ with quotient $A_F$, which is in turn an abelian scheme over
$Y$.  The choice $\Sigma_F$ of polyhedral cone decomposition determines a
torus embedding $T_F \hookrightarrow T_{F,\Sigma}$ and a partial
compactification $ M'_{F,\Sigma} = M'_F \times_{T_F} T_{F,\Sigma}$
of $M'_F.$  Let $D_{F,\Sigma}$ denote the interior of the closure of
$M'_F = \Gamma'_F \backslash D$ in $M'_{F,\Sigma}.$   The quotient mapping $\Gamma'_F 
\backslash D \to X$ extends to a local isomorphism $\varphi_{F,\Sigma}: D_{F,\Sigma} 
\to \overline{X}_{\Sigma}$  (cf. \cite{AMRT} p. 250).  In other words,
$\varphi_{F,\Sigma}$ is an open analytic mapping with discrete fibers which, near the
boundary, induces an embedding $D_{F,\Sigma}/(\Gamma'_F\backslash \Gamma_P)
\hookrightarrow \overline{X}_{\Sigma}$ whose image is the neighborhood $U_Y$ referred to
in
\S \ref{pf-tor1}.  The mappings $\varphi_{F,\Sigma}$ for the
various strata $Y\subset \overline{X}$ cover the boundary of $\overline{X}_{\Sigma}.$
The composition $\theta_1\theta_2|D: D \to F$ coincides with the canonical projection
$\pi.$

\subsection{}  Each of the spaces in Figure \ref{fig-1} comes equipped with 
a vectorbundle and most of the mappings in this diagram are covered by vectorbundle
isomorphisms.  We will make the following notational convention:  If $E_1 \to M_1$ and
$E_2 \to M_2$ are (smooth) vectorbundles, $\alpha: M_1 \to M_2$ is a (smooth) mapping,
and $\Phi:E_1 \to E_2$ is a vectorbundle mapping which induces an isomorphism $E_1 \cong
\alpha^*(E_2)$, then we will write $\Phi:(E_1,M_1)\sim (E_2,M_2)$ and refer to $\Phi$ as
being a {\it vectorbundle isomorphism which covers $\alpha$}.

As in \cite{Mumford}, complexify $\lambda:K \to GL(V)$ and extend it
trivially over $P^{-}$ to obtain a representation $\tilde\lambda:K(\mathbb C)P^{-} \to
GL(V)$.  The homogeneous vectorbundle $E = G\times_K V$ has a canonical extension 
\[\check{E} = \mathbf G(\mathbb C) \times_{K(\mathbb C)P^{-}} V\] over the
compact dual symmetric space $\check{D}.$  Its restriction $E_F$ to $D_F$ is a 
$PZ_F(\mathbb C)$-homo\-gen\-eous bundle, and it passes to a vectorbundle 
$E'_F \to M'_F$ upon dividing by $\Gamma'_F.$  The restriction $E'_F|(\Gamma'_F
\backslash D)$ coincides with the vectorbundle obtained from the homogeneous vectorbundle
$E = G \times_K V$ upon dividing by $\Gamma'_F$.  We will denote this restriction also by
$(E'_F, \Gamma'_F \backslash D).$

Define $\widetilde{E} = P.Z_F(\mathbb C)  \times_{K_P.Z_F(\mathbb C)}V.$  This
vectorbundle on $D_F/Z_F(\mathbb C)$ is homogeneous under $P.Z_F(\mathbb C)$  and it
passes to a vectorbundle $E^A_F \to A_F$ upon dividing by $\Gamma'_F.$  As in \cite{HZ1}
(3.2.1) and (3.3.5), there is a  canonical vectorbundle isomorphism 
\begin{equation}\label{eqn-psi}
\psi:(E'_F, M'_F) \sim (E^A_F, A_F) \end{equation}
which covers $\pi_2.$  In fact, the isomorphism $\psi$ is obtained from the canonical
isomorphism of $P.Z_F(\mathbb C)$-homogeneous vectorbundles,
\begin{equation}
\Psi : (E_F,D_F) \sim (\widetilde{E},D_F/Z_F(\mathbb C)) \end{equation}
which covers $\theta_2$ and which is given by the quotient mapping
\[ E_F = P.Z_F(\mathbb C)\times_{K_P}V \to \widetilde{E} = P.Z_F(\mathbb C) 
\times_{K_P.Z_F(\mathbb C)}V. \]
Let $\bar E'_{\Sigma}$ denote Mumford's canonical extension of the vectorbundle $E' =
\Gamma \backslash E \to X$ to the toroidal compactification $\overline{X}_{\Sigma}$, and
let $\bar E'_{F,\Sigma}= \varphi_{F,\Sigma}^*(E'_{\Sigma})$ be its pullback to
$D_{F,\Sigma}.$  Then we have a further canonical identification $E' = \bar
E'_{F,\Sigma}|(\Gamma'_F\backslash D)=E'_F|(\Gamma'_F\backslash D).$  We also have
vectorbundles $E_h = G_h\times_{K_h}V$ on $F$ and its quotient $E'_Y \to Y = \Gamma_h
\backslash F.$

According to \cite{HZ1} (3.3.9) (which in turn relies on \cite{Harris2}), {\it the
canonical isomorphism} \ref{eqn-psi} {\it extends to a vectorbundle isomorphism}
\begin{equation}
\psi_{\Sigma}: (\bar E'_{F,\Sigma},D_{F,\Sigma}) \sim (E^A_F,A_F) \end{equation}
{\it which covers} $\pi_{2,\Sigma}.$  This is the key point in the argument:  the
isomorphism $\psi_{\Sigma}$ identifies Mumford's canonical extension (which is defined
using a growth condition on a singular connection) with a vectorbundle,
$\pi_{2,\Sigma}^*(E^A_F)$ which is defined topologically, and which is trivial on each
torus embedding $\pi_{2,\Sigma}^{-1}(a) \cong T_{F,\Sigma}.$  We will use this
isomorphism to extend the parabolically induced connection over the toroidal
compactification, because such a parabolically induced connection is also pulled up from
$E^A_F.$

As in (\ref{eqn-Paction}), define an action of $P.Z_F(\mathbb C)$ on the vectorbundle
$E_h \to F$ by
\begin{equation}
ug_hg_{\ell}z.[g'_h,v] = [g_hg'_h,\lambda_1(g_{\ell})v]\end{equation}
(where $u\in \mathcal U_P$, $g_h,g'_h \in G_h$, $g_{\ell} \in G_{\ell}$, $z\in
Z_f(\mathbb C)$, and $v\in V$).  Define a mapping
\begin{equation}
\tilde\Phi_F: PZ_F(\mathbb C) \times_{K_PZ_F(\mathbb C)}V \to G_h \times_{K_h}V
\end{equation}
by $\tilde\Phi_F([ug_hg_{\ell}z,v]) = [g_h, \lambda_1(g_{\ell})v].$  Then $\tilde\Phi_F$
is well defined, it is $P.Z_F(\mathbb C)$-invariant, and it gives a $P.Z_F(\mathbb
C)$-equivariant isomorphism of vectorbundles,
\begin{equation}\label{eqn-newphi} 
{\Phi_F}:(\widetilde{E},D_F/Z_F(\mathbb C)) \sim (E_h,F) \end{equation}
which covers $\theta_1.$  Moreover the composition $\Phi_F\Psi|(E,D)$ is precisely the
isomorphism $\Phi:(E,D) \cong \pi^*(E_h,F)$ of (\ref{eqn-Phi}).
In summary, this array of vectorbundles appears in figure \ref{fig-2}.
\begin{figure}[h!]
\begin{equation*} \begin{CD}
(E,D) &\overset{\text{mod }\Gamma'_F}{-\dashrightarrow} & (E'_F,\Gamma'_F\backslash D) 
& \overset{\text{extend}}{-\dashrightarrow} & (\bar E'_{F,\Sigma},D_{F,\Sigma}) 
& \overset{\text{extend}}{-\dashrightarrow} & 
(\bar E'_{\Sigma},\overline{X}_{\Sigma}) \\
\begin{CD} @V{\text{extend}}VV \\ (E_F,D_F) \\ @V{\Psi}V{\sim}V \end{CD} 
& -\dashrightarrow & 
\begin{CD} @V{\text{extend}}VV \\ (E'_F,M'_F) \\ @V{\psi}V{\sim}V \end{CD} 
&& @V{\psi_{\Sigma}}V{\sim}V \\
(\widetilde{E}, D_F /Z_F(\mathbb C)) & {-\dashrightarrow} &
 (E^A_F,A_F) @= (E^A_F,A_F)  \\
@V{{\Phi}}V{\sim}V @VVV &&& \\
(E_h,F) & \overset{\text{mod }\Gamma_h}{-\dashrightarrow} & (E'_Y,Y) &&&
\end{CD}   
\end{equation*}
\caption{Vectorbundles on the toroidal compactification}\label{fig-2}\end{figure}

\subsection{} 
Each of the vectorbundles in Figure \ref{fig-2} comes equipped with a connection.
Let $\nabla_Y$ be a given connection on $E'_Y \to Y$ and let $\nabla_h$ be its pullback
to $E_h \to F.$  Define $\tilde\nabla = \Phi^*_F(\nabla_h)$ to be the pullback of
$\nabla_h$ under the isomorphism (\ref{eqn-newphi}).  Then $\Psi^*(\tilde\nabla)$
is an extension of the parabolically induced connection $\Phi^*(\nabla_h)$ on $(E,D).$ 
Both connections are invariant under $\Gamma'_F.$ Let $\nabla'_F$ denote the resulting
connection on the quotient $(E'_F,M'_F)$ (where again we use the same symbol to denote
this connection as well as its restriction to $(E'_F,\Gamma'_F \backslash D)$).  We need
to show that this connection $\nabla'_F$ on $(E'_F, \Gamma'_F \backslash D)$ has a smooth
extension to a connection $\nabla'_{F,\Sigma}$ on $(\bar E'_{F,\Sigma}, D_{F,\Sigma})$
which is invariant under $\Gamma'_F \backslash \Gamma_P.$

The connection $\tilde\nabla$ passes to a connection $\nabla^A_F$ on $(E^A_F,A_F)$ such
that $ \nabla'_F = \psi^*(\nabla^A_F).$ Therefore the connection 
\[ \nabla'_{F,\Sigma} = \psi^*_{\Sigma}(\nabla^A_F) \]
on $(\bar E'_{F,\Sigma}, M'_{F,\Sigma})$ is a smooth extension of $\nabla'_F.$  
 The $\Gamma'_F\backslash \Gamma_P$-invariance of $\nabla'_{F,\Sigma}$
follows from the $P.Z_F(\mathbb C)$-invariance of $\widetilde{\nabla}.$  This completes
the proof of Proposition \ref{prop-extension}.  \qed

\section{Chern classes and constructible functions}\label{sec-constructible}
\subsection{} A {\it constructible function} $F:W \to \mathbb Z$ on a complete (complex)
algebraic variety $W$ is one which is constant on the strata of some algebraic (Whitney)
stratification of $W$.  The Euler characteristic of such a constructible function $F$ is
the sum
\begin{equation*}
\chi(W;F) = \sum_{\alpha}\chi(W_{\alpha})F(W_{\alpha})
\end{equation*}
over strata $W_{\alpha}\subset W$ along which the function $F$ is constant.
If $f:W\to W'$ is an (proper) algebraic mapping, then the pushforward of the
constructible function $F$ is the constructible function
\begin{equation}\label{eqn-pushforward}
f_*(F)(w') = \chi(f^{-1}(w');F)
\end{equation}
(for any $w'\in W'$).
According to \cite{Mac}, for each constructible function $F:W\to \mathbb Z$
it is possible to associate a unique Chern class $c_*(W;F)\in H_*(W;\mathbb Z)$ which
depends linearly on $F$, such that
$f_*c_*(W;F) =c_*(W';f_*F)$ (whenever $f:W \to W'$ is a proper morphism), and
such that $c_*(W;\mathbf 1_W) = c^*(W)\cap [W]$ if $W$ is nonsingular.
(Here, $[W]\in H_{2\dim(W)}(W;\mathbb Z)$ denotes the fundamental class of $W$.)
The {\it MacPherson-Schwartz} Chern class of $W$ is the Chern class of the constructible
function $\mathbf 1_W.$

\subsection{}
Now let $\overline{Z}$ be a nonsingular complete complex algebraic variety and
let $D = D_1 \cup D_2 \cup \ldots \cup D_m$ be a union of smooth divisors with
normal crossings in $\overline{Z}$.  Set $Z = \overline{Z} - D.$  The tangent bundle
$T_Z$ 
of $Z$ has a ``logarithmic'' extension to $\overline{Z}$,
\begin{equation*}
T_{\overline Z}(-\log D) = \text{Hom}(\Omega^1_{\overline{Z}}(\log D), 
\mathcal O_{\overline{Z}})
\end{equation*}
which is called the ``log-tangent bundle'' of $(\overline{Z},D).$  It is the vectorbundle
whose sheaf of sections near any k-fold multi-intersection $\{z_1 = z_2 = \ldots =
z_k = 0\}$ of the divisors is generated
by $z_1\frac{\partial}{\partial z_1}, z_2\frac{\partial}{\partial z_2},\ldots,
z_k\frac{\partial}{\partial z_k},z_{k+1},\ldots,z_n$ (where $n=\dim(Z)$).  The following
result was discovered independently by P. Aluffi \cite{Aluffi}.

\begin{prop}\label{prop-log}
The Chern class of the log tangent bundle is equal to the Chern class of the
constructible
function which is $1$ on $Z = \overline{Z}-D,$ that is,
\begin{equation*}
c^*(T_{\overline{Z}}(-\log D)) \cap [\overline{Z}] = c_*(\mathbf 1_Z).
\end{equation*}
\end{prop}

\subsection{Proof} 
For any subset $I\subset \{1,2,\ldots,m\}$ let $D_I = \bigcap_{i\in I}D_i$, let 
\[ D^I = D_I \cap \bigcup_{j\notin I}D_j \]
denote the ``trace'' of the divisor $D$ in $D_I$, and let $D_I^o = D_I - D^I$ denote its
complement.  The restriction of the log tangent bundle of $(\overline{Z},D)$ to
any intersection $D_I$ is (topologically) isomorphic to the direct sum of vectorbundles
\begin{equation}\label{eqn-log}
T_{\overline{Z}}(-\log D)|D_I \cong T_{D_I}(-\log(D^I)) \oplus |I|\mathbf{1}
\end{equation}
(the last symbol denoting $|I|$ copies of the trivial bundle).  (This follows from the
short exact sequence of locally free sheaves on $D_j$,
\begin{equation*}\begin{CD}
0 @>>> \Omega^1_{D_j}(\log D^{\{j\}}) @>>>
\Omega^1_{\overline{Z}}(\log D)|D_j @>>> 
 \mathcal O_{D_j} @>>> 0
\end{CD}\end{equation*}
by dualizing and induction.)  We will prove Proposition \ref{prop-log} by induction on
the number $m$ of divisors, with the case $m=0$ being trivial.  For any constructible
function $F$ on $\overline{Z}$, denote by $c(F)\in H^*(\overline{Z})$ the Poincar\'e dual
of the (homology) Chern class of $F$.  Each divisor $D_j$ carries a fundamental homology
class whose Poincar\'e dual we denote by $[D_j]\in H^2(\overline{Z}).$ The Chern
class of the line bundle $\mathcal O(D_j)$ is $1+[D_j].$  Let $\tilde c $ denote the
Chern class of the bundle $T_{\overline{Z}}(-\log D).$ If $I\subset \{1,2,\ldots,m\}$ and
if $i:D_I \to \overline{Z}$ denotes the inclusion then
\[
\tilde c\cdot[D_I] = i_*(c(T_{\overline{Z}}(-\log D))|D_I) = i_*c(T_{D_I}(-\log(D^I))
= i_*c(\mathbf 1_{D_I^o}) 
\]
by (\ref{eqn-log}) and induction.  Using \cite{Tsushima} Proposition 1.2 we see,
\begin{align*}
c(\overline{Z}) &= \tilde c \cdot \prod_i (1+[D_i]) = 
\tilde c + \tilde c \cdot \sum_{I}[D_I] \\
&= \tilde c +  \sum_I c(\mathbf 1_{D_I^o}) = \tilde c + c(\mathbf 1_D)
\end{align*}
since each point in $D$ occurs in exactly one multi-intersection of divisors.  \qed

\begin{thm}\label{thm-constructible}
Let $X = \Gamma \backslash G /K$ be a Hermitian locally symmetric space as in \S
\ref{sec-BB}, with Baily-Borel compactification $\overline{X}.$  Let $\bar
c^i(\overline{X})\in H^{2i}(\overline{X};\mathbb C)$ denote the cohomology Chern class of
the tangent bundle, constructed in Theorem \ref{thm-A}.  Then its homology image
\[ \bar c^*(\overline{X}) \cap [\overline{X}] = c_*(\mathbf 1_X) \in H_*(\overline{X};
\mathbb Z) \]
lies in integral homology and coincides with the (MacPherson) Chern class of the
constructible function which is $1$ on $X$ and is $0$ on $\overline{X} - X.$
\end{thm}

\subsection{Proof}
Let $\tau:\overline{X}_{\Sigma}\to \overline{X}$ denote a smooth toroidal resolution of
singularities, having chosen the system of polyhedral cone decompositions $\Sigma$ so
that the exceptional divisor $D$ is a union of smooth divisors with normal crossings.
Let $T_{\overline{X}_{\Sigma}}(-\log D)$ denote the log tangent bundle of
$(\overline{X}_{\Sigma},D).$  As in \cite{Mumford} Prop. 3.4, this bundle is isomorphic
to Mumford's canonical extension $\overline{T}_{X,\Sigma}$ of the tangent bundle.
Therefore
\begin{align*}
\bar c^*(\overline{X})\cap [\overline{X}] &= \tau_*(\tau^*\bar c^*(\overline{X}) \cap 
[\overline{X}_{\Sigma}]))  \\
&= \tau_*(c^*(\overline{T}_{X,\Sigma})\cap [\overline{X}_{\Sigma}])\\
&= \tau_*c_*(\mathbf 1_X)=c_*(\tau_*(\mathbf 1_X)) = c_*(\mathbf 1_X)
\end{align*}
by Theorem \ref{thm-toroidal}, Proposition \ref{prop-log} and (\ref{eqn-pushforward}).
\qed

\begin{cor}
The MacPherson-Schwartz Chern class of the Baily-Borel compactification $\overline{X}$
is given by the sum over strata $Y\subset \overline{X}$,
\[
c_*(\mathbf 1_{\overline{X}}) = c_*(\sum_{Y\subset \overline{X}}\mathbf 1_Y) = 
\sum_{Y\subset \overline{X}}i_*\bar c^*(\overline{Y})\cap [\overline{Y}]
\]
where $i:\overline{Y} \hookrightarrow \overline{X}$ is the inclusion of the closure of
$Y$ (which is also the Baily-Borel compactification of $Y$) into $\overline{X}.$  \qed
\end{cor}

\section{Cohomology of the Baily-Borel Compactification}\label{sec-cohomology}

\subsection{}\label{subsec-BK}
Let $K$ be a compact Lie group and let $EK\to BK$ be the universal principal $K$-bundle.
For any representation $\lambda:K\to GL(V)$ on a complex vectorspace $V$, let
$E_{\lambda} = EK\times_KV$ be the associated vectorbundle.  The Chern classes $c^i(E)\in
H^{2i}(BK;\mathbb C)$ of all such vectorbundles generate a subalgebra which we denote
$H^*_{\text{Chern}}(BK;\mathbb C).$  Two cases are of particular interest:  if $K=
\mathbf{U(n)}$ then $BK = \lim_{k\to\infty}G_n(\mathbb C^{n+k})$ is the infinite
Grassmann manifold and $H^*(BK;\mathbb C) = H^*_{\text{Chern}}(BK;\mathbb C)$.  In fact,
the standard representation $\lambda:\mathbf{ U(n)} \to \mathbf{GL_n}(\mathbb C)$ gives
rise to a single vectorbundle $E_{\lambda}\to BK$ such that the algebra $H^{*}(BK;\mathbb
C)$ is canonically isomorphic to the polynomial algebra in the Chern
classes $c^1(E_{\lambda}), c^2(E_{\lambda}),\ldots,c^n(E_{\lambda}).$ If $K=\mathbf{ S
O(n)}$ then $BK=\lim_{k\to\infty}G^o_n(\mathbb R^{n+k})$ is the infinite Grassmann
manifold of real oriented n-planes.  Let $\tilde\lambda:\mathbf{ S O(n)}\to \mathbf{
GL_n}(\mathbb R)$ be the standard representation with resulting vectorbundle
$E_{\tilde\lambda} \to BK,$ and let $\lambda:\mathbf{S O(n)} \to \mathbf {GL_n}(\mathbb
C)$ denote the composition of $\tilde\lambda$ with the inclusion $\mathbf{ GL_n}(\mathbb
R) \subset \mathbf{ GL_n}(\mathbb C)$.  The associated vectorbundle $E_{\lambda} =
E_{\tilde\lambda}(\mathbb C)$ is the complexification of $E_{\tilde\lambda}$.  If $n$ is
odd, then $H^*(BK;\mathbb C)$ is canonically isomorphic to the polynomial algebra
generated by the Pontrjagin classes $p^i(E_{\tilde\lambda}) =c^{2i}(E_{\lambda})\in
H^{4i}(BK;\mathbb C)$ for $i=1,2,\ldots,n.$  Hence $H^*(BK;\mathbb C) =
H^*_{\text{Chern}} (BK;\mathbb C).$   If $n$ is even then the algebra $H^*(BK;\mathbb C)$
has an additional generator, the Euler class $e = e(E_{\tilde\lambda})\in
H^{n}(BK;\mathbb C)$.  (It satisfies $e^2 = p^{n/2}$.)  If $n=2$ then $e$ is the first
Chern class of the line bundle corresponding to the representation $\mathbf{SO}(2) \cong
\mathbf U(1) \subset \mathbf{GL_1}(\mathbb C).$  

\subsection{}\label{subsec-BK2}
Now suppose that $K=K_1\times K_2\times \ldots \times K_r$ is a product of unitary
groups, odd orthogonal groups, and copies of $\mathbf{SO(2)}.$  According to the
preceding paragraph,
there are representations $\lambda_1,\ldots,\lambda_r$ of $K$ on certain complex
vectorspaces $V_1,V_2,\ldots, V_r$ so that the Chern classes of the resulting
``universal'' complex vectorbundles $E_i = EK\times_K V_i \to BK$ generate the polynomial
algebra $H^*(BK;\mathbb C) = H^*_{\text{Chern}}(BK;\mathbb C).$

\subsection{}\label{subsec-table}
Suppose that $\mathbf G$ is a semisimple algebraic group defined over $\mathbb Q$, and
that $\mathbf G(\mathbb R)^0$ acts as the identity component of the group of
automorphisms of a Hermitian symmetric space $D = G/K.$  Recall (\cite{Helgason} X \S 6,
\cite{Borelthesis}) that the irreducible components of $D$ come from the following list:

\medskip

\centerline{\begin{tabular}{ | l | l | l |}
\hline 
Type & Symmetric Space & Compact Dual \\ \hline
AIII & $\mathbf{  U (p,q)}/ \mathbf {U(p) \times U(q)}$ & $\mathbf{U(p+q)}/\mathbf{U(p)
\times U(q)}$ \\
DIII & $\mathbf{SO^*(2n)}/\mathbf{U(n)}$ & $\mathbf{SO(2n)}/\mathbf{U(n)}$ \\
BDI  & $\mathbf{SO(p,2)}/\mathbf{SO(p)}\times \mathbf{SO(2)}$ & $
\mathbf{SO(p+2)}/\mathbf{SO(p)}\times \mathbf{SO(2)}$ \\
CI   & $\mathbf{Sp(n,\mathbb R)}/\mathbf{U(n)}$ & $ \mathbf{Sp(n)}/\mathbf{U(n)} $ \\
EIII & $\mathbf{E^3_6}/\mathbf{Spin(10)}\times \mathbf{SO(2)} $ &
$\mathbf{E_6}/\mathbf{Spin(10)}\times
\mathbf{SO(2)}$ \\
EVII & $\mathbf{E^3_7}/\mathbf{E_6}\times \mathbf{SO(2)}$ &
$\mathbf{E_7}/\mathbf{E_6}\times
\mathbf{SO(2)} $ \\ \hline \end{tabular}}

\medskip

Let $X = \Gamma\backslash G/K$, with $\Gamma\subset \mathbf G(\mathbb Q)$ a neat
arithmetic group, and let $\overline{X}$ denote the Baily-Borel
compactification of $X$.  Let $\check{D} = G_u/K$ be the compact dual symmetric space,
where $G_u\subset\mathbf G(\mathbb C)$ is a compact real form containing $K$.  The
principal bundles $\Gamma\backslash G \to X$ and $G_u \to \check{D}$ are classified
by mappings $\Phi:X \to BK$ and $\Psi:\check{D} \to BK$ (respectively) which are uniquely
determined up to homotopy.  A theorem of Borel \cite{Borelthesis} states that (in this
Hermitian case) the resulting homomorphism $\Psi^{*}:H^{*}(BK;\mathbb C) \to
H^*(\check{D};\mathbb C)$ is surjective.

Suppose the irreducible factors of $D = G/K$ are of type AIII, DIII, CI, or BDI for $p$
odd or $p=2.$  The construction of $\pi$-fiber Chern forms in Section
\ref{sec-maintheorem} determines a homomorphism $\widetilde{\Phi}^{*}:H^*(BK;\mathbb C)
\to H^*(\overline{X};\mathbb C)$ by setting $\widetilde{\Phi}^*(c^i(E_j)) = \bar
c^i(E'_{j})$ (where $E_j \to BK$ is the universal vectorbundle corresponding to the
representation $\lambda_j$ of \S \ref{subsec-BK2} and $E'_{j} \to X$ is the
corresponding automorphic vectorbundle).  Let us denote the image of $\widetilde\Phi^*$
by $H^*_{\text{Chern}}(\overline{X};\mathbb C).$

\begin{thm}\label{thm-cohomology} Suppose $X=\Gamma\backslash G/K$ is a Hermitian locally
symmetric space such that
the irreducible factors of $D=G/K$ are of type AIII, DIII, CI, or BDI for $p$ odd or
$p=2$. Then the  mappings $\widetilde{\Phi}^{*}$ and $\Psi^{*}$ determine a surjection
\begin{equation}h:\label{eqn-surjection}H^*_{\text{Chern}}(\overline{X};\mathbb C) \to
H^*(\check{D};\mathbb C)\end{equation} from this subalgebra of the cohomology of the
Baily-Borel compactification, to the cohomology of the compact dual symmetric space.
Moreover, for each ``universal'' vectorbundle $E_j \to BK$ we have
\begin{equation*} h(\bar c^i(E'_{j})) = c^i(\check E_{j}) \end{equation*}
where $E'_{j} \to X$ and $\check E_{j} \to \check D$ are the associated automorphic
and homogeneous vectorbundles, respectively.
\end{thm}

\subsection{Proof} 
Define the mapping $h:H^*_{\text{Chern}}(\overline{X};\mathbb C) \to
H^*(\check{D};\mathbb C)$ by $h\tilde{\Phi}^*(c) = {\Psi}^*(c)$ for any $c\in
H^*(BK;\mathbb C)$.  If this is well defined, it is surjective by Borel's theorem.  To
show it is well defined, let us suppose that $\tilde{\Phi}^*(c)=0.$  We must show that
$\Psi^*(c)=0,$  so we assume the contrary.

Let $x=\Psi^*(c)\in H^i(\check{D};\mathbb C).$  By Poincar\'e duality, there exists a
complementary class $y\in H^{2n-i}(\check{D};\mathbb C)$ so that $(x\cup y )\cap
[\check{D}] \ne 0$ (where $n = \dim_{\mathbb C}(\check{D})$).  Then $y$ has a lift, $d\in
H^{2n-i}(BK;\mathbb C)$ with $\Psi^*(d)=y.$  Let us write $K=K_1K_2\ldots K_r$ for the
decomposition of $K$ into irreducible factors.  By \S \ref{subsec-BK2} the polynomial
algebra $H^*(BK;\mathbb C)$ is generated by the Chern classes of the universal
vectorbundles  $E_1,E_2,\ldots,E_r$ corresponding to representations $\lambda_i:K_i \to
GL(V_i).$ Hence, both $c$ and $d$ are polynomials in the Chern
classes of the vectorbundles $E_{1}, E_{2}, \ldots, E_{r}.$   Hence
$(\widetilde\Phi^*(c)\cup\widetilde\Phi^*(d))\cap
[\overline{X}] \in \mathbb C$ is a sum of ``generalized'' Chern numbers which, by the
Proposition \ref{prop-prop}, coincides with the corresponding sum of ``generalized''
Chern numbers for the compact dual symmetric space, $v(\Gamma)(x\cup y)\cap [\check{D}]
\ne 0.$ This implies that $\tilde{\Phi}^*(c)\ne 0$ which is a contradiction.  \qed

\subsection{Remarks}  We do not know whether the surjection (\ref{eqn-surjection}) has a
canonical splitting.  However, the intersection cohomology $IH^*(\overline{X};\mathbb
C)$ contains, in a canonical way, a copy of the cohomology $H^*(\check{D};\mathbb C)$ of
the compact dual symmetric space.  By the Zucker conjecture (\cite{Looijenga} and
\cite{SaperStern}), the intersection
cohomology may be identified with the $L^2$ cohomology of $X$ which, in turn may be
identified with the relative Lie algebra cohomology $H^*(\mathfrak g,
K;L^2(\Gamma\backslash G))$.  But $L^2(\Gamma\backslash G)$ contains a copy of the
trivial representation $\mathbf 1$ (the constant functions), whose cohomology
$H^*(\mathfrak g,K;\mathbf 1) \cong H^*(\check{D};\mathbb C)$ is the
cohomology of the compact dual symmetric space.  We sketch a proof that the following
diagram commutes.
\begin{equation*}\begin{CD}
H^{2k}_{\text{Chern}}(\overline{X};\mathbb C) @>>j> IH^{2k}(\overline{X};\mathbb C) \\
@VhVV @AAiA \\
H^{2k}(\check{D};\mathbb C) @= H^{2k}(\check D; \mathbb C)
\end{CD}\end{equation*}

If $E' \to X$ and $\check E \to \check D$ are vectorbundles arising from the same
representation $\lambda$ of $K$ then the class $j(\bar c^k(E'))$ is represented by the
differential form $\sigma^k(\nabla^{\sf p}_X)$ which is $\pi$-fiber, hence bounded, hence
$L^2.$  The class $i(c^k(\check E))$ is represented by the differential form
$\sigma^k(\nabla^{\text{Nom}}_X)$ which is ``invariant'' (meaning that its pullback to
$D$ is invariant), hence $L^2.$  The intersection cohomology of $\overline{X}$ embeds
into the ordinary cohomology of any toroidal resolution $\overline{X}_{\Sigma}.$  But
when these two differential forms are considered on $\overline{X}_{\Sigma}$, they both
represent the same cohomology class, $c^k(\overline{E}_{\Sigma})$ (using Theorem 
\ref{thm-toroidal} and \cite{Mumford}).  Alternatively, one may deform the connection
$\nabla^{\sf p}_X$ to $\nabla^{\text{Nom}}_X,$ obtaining a differential form $\Psi \in
\mathcal A^{2k-1}(X)$ such that $d\Psi = \sigma^k(\nabla^{\sf p}_X) -
\sigma^k(\nabla^{\text{Nom}}_X)$, and check that $\Psi$ is $L^2.$

In the case BDI, the compact dual is $\check{D} = SO(p+2)/SO(p)\times
SO(2)$.  Suppose $p$ is even.  Its cohomology
$H^*(\check{D};\mathbb C)$ has a basis $\{ 1, c_1,
c_1^2, \ldots, c_1^{p-1},e\}$ where $c_1$ is the Chern class of the complexification of
the line bundle arising from the standard representation of $SO(2)$ and where $e$ is the
Euler class of the vectorbundle arising from the standard representation of $SO(p),$
\cite{BorelHirzebruch} \S 16.5.  All
these classes lift canonically to $IH^*(\overline{X};\mathbb C)$ and $c_1^j$ lifts
further to the (ordinary) cohomology of the Baily-Borel compactification.  However
(except in the case $p=2$) we do not know whether $e$  also lifts further to
$H^*(\overline{X}; \mathbb C).$

\end{document}